\documentclass[11pt,notitlepage]{article}
\usepackage[dvips]{graphicx}
\usepackage{afterpage,amsmath,amsthm,amssymb,epsfig,fullpage,pictex,verbatim,wasysym}

\newtheorem{theorem}{Theorem}

\newtheorem{proposition}[theorem]{Proposition}

\author{Matt Davis\footnote{This material is based upon work supported by the National Science
Foundation under Grant No. DMS-0839966.}\\
Department of Mathematics\\
Harvey Mudd College\\
davis@math.hmc.edu}
\title{Representations of Rank Two Affine Hecke Algebras at Roots of Unity}

\begin{document}

\maketitle

\begin{abstract}
In this paper, we will fully describe the representations of the crystallographic rank two affine Hecke algebras using elementary methods, for all possible values of $q$. The focus is on the case when $q$ is a root of unity of small order.
\end{abstract}

\section{Introduction}

The affine Hecke algebra $\widetilde{H}$ was introduced by Iwahori and Matsumoto (\cite{I-M}). Knowing the representations of $\widetilde{H}$ gives a substantial amount of information about the representations of a closely related $\mathfrak{p}$-adic group. The definition of $\widetilde{H}$ involves a parameter $q$ which can have a large effect on the structure of the algebra. In this paper, we will fully describe the representations of the affine Hecke algebras of type $C_{2}$ and $G_{2}$, for all possible values of $q$. The methods are essentially those introduced in \cite{Ra1}, with the modifications required to deal with $q$ being a root of unity.

The representations in type $A$ were described in the non-root of unity case by Zelevinsky, in terms of combinatorial objects called multisegments (see \cite{BZ} and \cite{Z}). In the root of unity case, these representations are indexed by the aperiodic multisegments (see the appendix of \cite{LTV} for an argument relying on the results of \cite{Lu}). The representations of $\widetilde{H}$ in all types have been classified geometrically by Kazhdan and Lusztig (\cite{KL}) in the non-root of unity case, and studied in the root of unity case by Grojnowski (\cite{G}) and N.Xi (\cite{Xi}, \cite{Xi2}), among others. In the root of unity case, Grojnowksi gives a simple description of a geometric indexing set ($\cite{G}$, Theorem 2) only in type $A$. However, Theorem 1 of \cite{G} does not apply in all cases (see the remark on p. 524 of \cite{Xi}). And, to this author, at least, it is not obvious how to turn the statement of Theorem 1 into Theorem 2. One hopes that a better understanding of the representations of $\widetilde{H}$ in some small cases will help clarify these issues.

We begin by defining the affine Hecke algebra $\widetilde{H}$, and recalling some basic facts about the representations of $\widetilde{H}$. We will make extensive use of $\mathbb{C}[X]$, a large commutative subalgebra of $\widetilde{H}$, and \textit{weights}, elements of $\textrm{Hom}(\mathbb{C}[X],\mathbb{C})$, which describe the simple reprsentations of $\mathbb{C}[X]$. A $\widetilde{H}$ module $M$ can be described in part by which weights appear, i.e. which simple $\mathbb{C}[X]$ modules are composition factors of it. The most important construction we will use is that of the \textit{principal series module} $M(t)$ which can be constructed from any weight $t \in T$, since every simple $\widetilde{H}$-module is a quotient of $M(t)$ for an appropriate choice of a weight $t$ (Theorem \ref{PrinQuot}c). We also recall from \cite{Ra1} several facts needed to analyze the modules $M(t)$, with some adaptations as necessary to deal with the root of unity case.

The main goal of the paper is to describe a way of visualizing and describing the composition factors of $M(t)$ directly from the combinatorial data of the weight $t$. This can be done with particular pictures based on the root system underlying $\widetilde{H}$. The following are examples in the Type $A_{2}$ case.

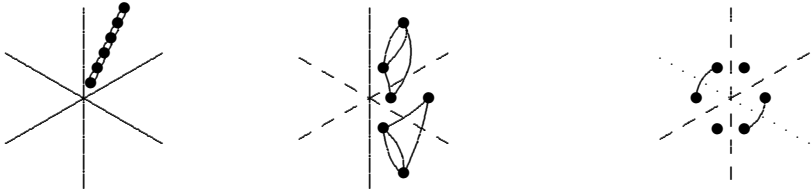
\[ \beginpicture
\setcoordinatesystem units <0.2cm,0.2cm>
\setplotarea x from -27 to 25, y from -5 to 6
\linethickness=3pt
\plot -23 6 -23 -6 /
\plot -17.804 3 -28.196 -3 /
\plot -17.804 -3 -28.196 3 /
\put{$\bullet$} at -21.2 4
\put{$\bullet$} at -21.65 3
\put{$\bullet$} at -22.1 2
\put{$\bullet$} at -20.75 5
\put{$\bullet$} at -20.3 6
\put{$\bullet$} at -22.55 1
\plot -4 6 -4 -6 /
\put{$\bullet$} at -2.6 0
\put{$\bullet$} at -0.1 0
\put{$\bullet$} at -3.1 2
\put{$\bullet$} at -1.75 5
\put{$\bullet$} at -3.1 -2
\put{$\bullet$} at -1.75 -5
\put{$\bullet$} at 17.7 0
\put{$\bullet$} at 22.3 0
\put{$\bullet$} at 19.1 2
\put{$\bullet$} at 20.9 2
\put{$\bullet$} at 19.1 -2.1
\put{$\bullet$} at 20.9 -2.1
\setquadratic
\plot -22.55 1 -22.55 1.5 -22.1 2 /
\plot -22.1 2  -22.1 2.5 -21.65 3 /
\plot -21.65 3 -21.65 3.5 -21.2 4 /
\plot -21.2 4 -21.2 4.5 -20.75 5 /
\plot -20.75 5 -20.75 5.5 -20.3 6 /
\plot -22.55 1 -22.1 1.5 -22.1 2 /
\plot -22.1 2  -21.65 2.5 -21.65 3 /
\plot -21.65 3 -21.2 3.5 -21.2 4 /
\plot -21.2 4 -20.75 4.5 -20.75 5 /
\plot -20.75 5 -20.3 5.5 -20.3 6 /
\plot -3.1 2 -2.6 4 -1.75 5 /
\plot -3.1 2 -1.7 3.8 -1.75 5 /
\plot -3.1 -2 -2.6 -4 -1.75 -5 /
\plot -3.1 -2 -2 -3.4 -1.75 -5 /
\plot -3.1 2 -2.75 1 -2.6 0 /
\plot -1.75 5 -1.3 2.5 -2.6 0 /
\plot -1.75 -5 -0.7 -2.5 -0.1 0 /
\plot -3.1 -2 -1.6 -1.2 -0.1 0 /
\setquadratic
\circulararc 60 degrees from 19.1 2.1 center at 20 0
\circulararc 60 degrees from 20.9 -2.1 center at 20 0
\setlinear
\setdashes
\plot 25.196 3 14.804 -3 /
\plot 20 6 20 -6 /
\plot 1.196 3 -9.196 -3 /
\plot 1.196 -3 -9.196 3 /
\setdots
\plot 25.196 -3 14.804 3 /
\endpicture \]

The lines in this picture represent the hyperplanes perpendicular to the roots $\alpha$, and are drawn as solid, shaded, or dotted based on the values of the weight $t$ on the elements $X^{\alpha} \in \widetilde{H}$. Each dot in the picture represents one dimension of the module $M(t)$, and dots are connected if a single composition factor of $M(t)$ contains both of these basis elements. The general goal is to determine a few rules that determine which of these lines should be drawn. That is, we hope to find a few algebraic statements that describe how $M(t)$ breaks down into composition factors which can be translated into these pictures. Essentially, Theorem $\ref{calibsame}$b, Theorem $\ref{2dims}$, and Lemma $\ref{eq:lemma}$ below are sufficient to complete the classification in the rank two cases, for all values of $q$. These pictures provide a very straightforward way of determining the composition factors of $M(t)$, without relying on heavy computations. One also hopes that a complete classification of the rank two crystallographic cases will facilitate a greater understanding of the representation theory of $\widetilde{H}$ in all types.

\section{Definitions}
\label{sec:definitions}

In this section, we introduce the needed definitions and several preliminary results about the affine Hecke algebra. Proofs of most previously known results will not be given.

\vspace{.1in}

\noindent \textbf{The Affine Hecke Algebra.} \label{subsectHecke} Let $R$ be a root system in $\mathbb{R}^{n}$ with simple roots $\alpha_{1}, \ldots \alpha_{n}$. Let $R^{+}$ be the set of positive roots and $R^{-}$ the set of negative roots. We define the \textit{rank} of $R$ to be the number of simple roots $n$.

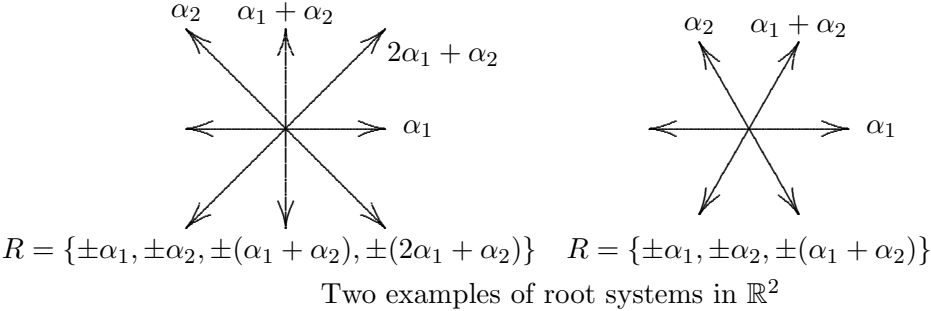
\[ \beginpicture
\setcoordinatesystem units <0.22cm,0.22cm>
\setplotarea x from -25 to 25, y from -5 to 6
\linethickness=5pt
\setlinear
\thicklines
\arrow <10 pt> [.2,.6] from -16 0 to -16 6
\arrow <10 pt> [.2,.6] from -16 0 to -16 -6
\arrow <10 pt> [.2,.6] from -16 0 to -10 -6
\arrow <10 pt> [.2,.6] from -16 0 to -22 6
\arrow <10 pt> [.2,.6] from -16 0 to -22 0
\arrow <10 pt> [.2,.6] from -16 0 to -10 0
\arrow <10 pt> [.2,.6] from -16 0 to -10 6
\arrow <10 pt> [.2,.6] from -16 0 to -22 -6
\arrow <10 pt> [.2,.6] from 12 0 to 15 5.196
\arrow <10 pt> [.2,.6] from 12 0 to 9 -5.196
\arrow <10 pt> [.2,.6] from 12 0 to 15 -5.196
\arrow <10 pt> [.2,.6] from 12 0 to 9 5.196
\arrow <10 pt> [.2,.6] from 12 0 to 18 0
\arrow <10 pt> [.2,.6] from 12 0 to 6 0
\put{Two examples of root systems in $\mathbb{R}^{2}$} at 0 -10
\put{$\alpha_{1}$} at -8 0
\put{$\alpha_{1}$} at 20 0
\put{$\alpha_{2}$} at 9 6.196
\put{$\alpha_{2}$} at -22 7
\put{$\alpha_{1} + \alpha_{2}$} at -16 7
\put{$\alpha_{1} + \alpha_{2}$} at 15 6.196
\put{$2\alpha_{1} + \alpha_{2}$} at -6.5 4.5
\put{$R = \{ \pm \alpha_{1}, \pm \alpha_{2}, \pm (\alpha_{1} + \alpha_{2}), \pm (2\alpha_{1} + \alpha_{2}) \}$} at -17 -7.4
\put{$R = \{ \pm \alpha_{1}, \pm \alpha_{2}, \pm (\alpha_{1} + \alpha_{2}) \}$} at 12 -7.4
\endpicture\]

The reflection through $H_{\alpha}$ will be denoted by $s_{\alpha}$, or $s_{i}$ for the reflection through $H_{\alpha_{i}}$. If $\pi/m_{ij}$ is the angle between $H_{\alpha_{i}}$ and $H_{\alpha_{j}}$, then $m_{ij} \in \{2,3,4,6\}$ for $1 \leq i,j \leq n$, and the \textit{Weyl Group} $W_{0}$ has presentation \[W_{0} = \left\langle s_{1}, \ldots s_{n} \, | \, s_{i}^{2} = 1 \,\, , \,\, \underbrace{s_{i}s_{j}s_{i} \ldots }_{m_{ij} \textrm{ factors}} = \underbrace{s_{j}s_{i}s_{j} \ldots }_{m_{ij} \textrm{ factors }}, \,\, \textrm{ for } 1 \leq i,j \leq n \right\rangle.\]

Let $P$ be the \textit{weight lattice}, spanned by the elements $\omega_{i}$ satisfying \[ \langle \omega_{i}, \alpha_{j} \rangle = \delta_{ij} \cdot \frac{1}{2} \langle \alpha_{j},\alpha_{j} \rangle,\] for $\alpha_{i}$ and $\alpha_{j}$ simple roots. Let $Q$ be the lattice spanned by the simple roots $\alpha_{i}$. Let \begin{equation} X = \{X^{\lambda} \, | \, \lambda \in P\}, \,\, \textrm{ with } X^{\lambda} \cdot X^{\mu} = X^{\lambda + \mu} \quad \textrm{ for } \lambda, \mu \in P. \end{equation} Then $W_{0}$ acts on $X$ by \[ w \cdot X^{\lambda} = X^{w\cdot \lambda},\] and this action extends linearly to an action of $W_{0}$ on the group algebra $\mathbb{C}[X]$.

The \textit{affine Hecke algebra} $\widetilde{H}$ is the $\mathbb{C}$-algebra generated by $\{T_{i} \, | \, i \in I \}$ and $\{ X^{\lambda} \, | \, \lambda \in P \}$, where $\mathbb{C}[X]$ is a subalgebra of $\widetilde{H}$, and subject to the relations

\begin{equation} \label{eq:quadrel} T_{i}^{2} = (q - q^{-1})T_{i} + 1,  \quad \textrm{ for } i = 1,2, \ldots n, \end{equation} \vspace{.001in}
\begin{equation} \label{eq:heckerel} \underbrace{T_{i}T_{j}T_{i} \ldots }_{m_{ij} \textrm{ factors}} = \underbrace{ T_{j}T_{i}T_{j}\ldots }_{m_{ij} \textrm{ factors}} \quad \textrm{ for }i \neq j, \,\, \textrm{ and } \end{equation}
\begin{equation} \label{eq:relation} X^{\lambda}T_{s_{i}} = T_{s_{i}}X^{s_{i}\cdot \lambda} + (q -
q^{-1})\frac{X^{\lambda} - X^{s_{i}\lambda}}{1 - X^{-\alpha_{i}}}, \quad \textrm{ for } \lambda \in P, 1 \leq i \leq n.\end{equation} The \textit{rank} of $\tilde{H}$ is defined to be the rank of the underlying root system $R$. For $w \in W_{0}$, let \[ T_{w} = T_{i_{1}}T_{i_{2}}\ldots T_{i_{k}}\] for a reduced word $w = s_{i_{1}}s_{i_{2}} \ldots s_{i_{k}}$ in $W_{0}$. Then $\{ T_{w}X^{\lambda} \, | \, w \in W_{0}, \lambda \in P\}$ is a $\mathbb{C}$-basis for $\widetilde{H}$.

\vspace{.1in}

\noindent \textbf{Weights.} Let $T = \textrm{ Hom}(X, \mathbb{C}^{\times})$ be the set of group homomorphisms from $X$ to $\mathbb{C}^{\times}$. Then $T$ is an abelian group with $W_{0}$-action given by \[ w \cdot t(X^{\lambda}) = t(X^{w^{-1}\cdot \lambda}) \qquad \textrm{ for } t \in T, w \in W_{0}, \lambda \in P. \] An element of $T$ is called a \textit{weight}. For a weight $t$, the subgroup of $W_{0}$ that fixes $t$ under this action is generated by $\{ s_{i} \, | \, t(X^{\alpha_{i}}) = 1 \}$. (This relies on the fact that we chose $P$ rather than $Q$ to build $\widetilde{H}$. See \cite{St}, 3.15, 4.2, and 5.3).

For any finite-dimensional $\widetilde{H}$-module $M$, define the
{\it $t$-weight space} and the {\it generalized $t$-weight space} of $M$ by

\renewcommand{\baselinestretch}{1} \normalsize

\[M_{t} = \{m \in M \, |\, X^{\lambda}\cdot m = t(X^{\lambda})m \textrm{
for all } X^{\lambda} \in X \} \textrm{, and }\]

\[M_{t}^{\textrm{gen}} = \{m \in M \, |\, \textrm{ for all } X^{\lambda}
\in X, (X^{\lambda} - t(X^{\lambda}))^{k} m = 0 \textrm{ for some } k \in \mathbb{Z}_{>0} \},\]
respectively. Then $M = \bigoplus_{t \in T}
M_{t}^{\textrm{gen}}$ is a decomposition of $M$ into Jordan blocks for
the action of $\mathbb{C}[X]$. An element $t \in T$ is a {\it weight of $M$} if $M_{t}^{\textrm{gen}} \neq 0$.

\vspace{.1in}

\noindent \textbf{Induced Modules and Intertwining Operators.} If $I \subseteq \{1,\ldots n\}$, define $W_{I} = \langle s_{i} \, | \, i \in I \rangle$ and \[ \widetilde{H}_{I} = \{ T_{w}X^{\lambda} \, | \, \lambda \in P, w \in W_{I}\}.\] For $t \in T$ such that
$t(X^{\alpha_{i}}) = q^{2}$ for $i \in I$, then define $\mathbb{C}v_{t}$ to be the one dimensional
$\widetilde{H}_{I}$-module spanned by $v_{t}$, with $\widetilde{H}_{I}$ action
given by

\[ T_{i}\cdot v_{t} = qv_{t} \quad \textrm{ and } \quad X^{\lambda}\cdot v_{t} =
t(X^{\lambda})v_{t}, \quad \textrm{ for } X^{\lambda} \in X.\]

\begin{proposition} \label{eq:ind} \label{weightbasis} \label{PrinQuot} (\cite{Ra1}, Lemma 1.17)
Let $\mathbb{C}v_{t}$ be defined as above, and let $M =$ Ind$_{\widetilde{H}_{I}}^{\widetilde{H}} \mathbb{C}v_{t}$. Let $W_{I} = \langle s_{i} \, | \, i \in I \rangle$, and let $W_{0}/W_{I}$ be a set of minimal length coset representatives of $W_{I}$-cosets in $W_{0}$.

(a) Then the weights of $M$ are $\{wt \, | \, w \in W_{0}/W_{I}\}$, and \[ \textrm{dim}(M^{\textrm{gen}}_{wt}) = (\# \textrm{ of } v \in W_{0}/W_{I} \textrm{ with } vt = wt).\]

(b) There is a basis of $M$ consisting of elements of the form \[ m_{w} = T_{w} v_{t} + \sum_{u < w, u \in W_{0}/W_{I}} p_{w,u} T_{u}v_{t},\] for $w \in W_{0}/W_{I}$, such that $m_{w} \in M_{wt}$.

(c) If $t$ is a weight of an irreducible $\widetilde{H}$-module $N$ and $I = \emptyset$, then $N$ is a quotient of $M$. In fact, if $v \in N$ is a non-zero vector in $N_{t}$, then

\renewcommand{\baselinestretch}{1} \normalsize
\[\phi:M \rightarrow N\] \[v_{t} \mapsto v\]
extends to a surjective $\widetilde{H}$-module homomorphism.

\end{proposition}

\vspace{.05in}

In particular, if $I = \emptyset$, then we call \[M(t) = \widetilde{H} \otimes_{\mathbb{C}[X]} \mathbb{C}v_{t} = \textrm{span}\{ T_{w}v_{t} \, | \, w \in W_{0}\}\] the \textit{principal series module} for $t$.

Part $c$ of this lemma implies that the weights of a single simple finite-dimensional module $M$ lie in a single orbit $Wt$. We call this orbit (and, by abuse of terminology, any element of the orbit) the \textit{central character} of $M$. In fact, $\widetilde{H}$ has finite dimension over its center, and thus all simple $\widetilde{H}$-modules are finite-dimensional (see \cite{Ra1}, section 2.3). Thus, this theorem tells us that understanding the composition factors of all the principal series modules $M(t)$ is sufficient for understanding all the simple $\tilde{H}$-modules.

For a weight $t$ with $t(X^{\alpha_{i}}) \neq 1$ and a $\widetilde{H}$-module $M$, define a $\mathbb{C}$-linear operator $\tau_{i}: M_{t}^{\textrm{gen}} \rightarrow M$ by \renewcommand{\baselinestretch}{1} \normalsize

\[\tau_{i}(m) = (T_{i} - \frac{q-q^{-1}}{1-X^{-\alpha_{i}}})\cdot m. \label{taudef} \]

\begin{theorem} \label{tauthm} (\cite{Ra1}, Prop. 1.18)

\renewcommand{\baselinestretch}{1} \normalsize
\begin{enumerate}
\item[(a)] $1-X^{-\alpha_{i}}$ is invertible as an operator on $M_{t}^{\textrm{gen}}$, so that $\tau_{i}:M^{\textrm{gen}}_{t} \rightarrow M$ is well-defined.

\item[(b)] As operators on $M_{t}^{\textrm{gen}}$, $X^{\lambda}\tau_{i} = \tau_{i}X^{s_{i}\lambda}$ for all $X^{\lambda} \in X$, so that $\tau_{i}(M_{t}^{\textrm{gen}}) \subseteq M_{s_{i}t}^{\textrm{gen}}$.

\item[(c)] As operators on $M_{t}^{\textrm{gen}}$, $\displaystyle{\tau_{i}\tau_{i} = \frac{(q-q^{-1}X^{\alpha_{i}})(q-q^{-1}X^{-\alpha_{i}})}{(1-X^{\alpha_{i}})(1-X^{-\alpha_{i}})}}$

\item[(d)] The maps $\tau_{i}:M_{t}^{\textrm{gen}} \rightarrow M_{s_{i}t}^{\textrm{gen}}$ and $\tau_{i}:M_{s_{i}t}^{\textrm{gen}} \rightarrow M_{t}^{\textrm{gen}}$ are both invertible if and only if $t(X^{\alpha_{i}}) \neq q^{\pm2}$.

\item[(e)] If $i \neq j$ and $m_{ij}$ is defined as in ($\ref{eq:heckerel}$), then $\underbrace{\tau_{i}\tau_{j}\tau_{i}\ldots}_{m_{ij} \textrm{ factors}} = \underbrace{ \tau_{j}\tau_{i}\tau_{j} \ldots}_{m_{ij} \textrm{factors}} $, whenever both sides are well-defined operators. \end{enumerate} \end{theorem}

For $t \in T$, the \textit{calibration graph} of $t$ is the graph with vertices labeled by the elements of the orbit $W_{0}t$ and edges $(wt, s_{i}wt)$ if $(wt)(X^{\alpha_{i}}) \neq q^{\pm 2}$. The $\tau$ operators are used to prove the following.

\begin{theorem} \label{PrinOrbit}

\begin{itemize}
\item[(a)] (\cite{Rog}, Proposition 2.3) If $w \in W_{0}$ and $t \in T$ then $M(t)$ and $M(wt)$ have the same
composition factors.

\item[(b)] \label{calibsame} (\cite{Ra1}, Prop. 1.6)
Let $M$ be a finite dimensional $\widetilde{H}$-module, and let $t$ and $wt$ be two elements of $W_{0}t$ in the same connected component of the calibration graph for $t$. Then

\[ \textrm{dim}(M_{t}^{\textrm{gen}}) =
\textrm{dim}(M_{wt}^{\textrm{gen}}).\]

\item[(c)] \label{eq:Kato} (\cite{K1}) $M(t)$ is irreducible if and only if $P(t) := \{ \alpha \in R^{+} \, | \, t(X^{\alpha}) = q^{\pm 2} \} = \emptyset$.
\end{itemize}

\end{theorem}

\vspace{.1in}

\noindent \textbf{The Structure of Modules.} Theorem \ref{PrinOrbit} shows that the connected components of the calibration graph encode certain sets of weights whose corresponding weight spaces $M_{t}^{\textrm{gen}}$ must have the same dimension in any irreducible $\tilde{H}$-module $M$. These ideas lead us to the following propositions which will be fundamental in our later classification.

\begin{proposition} \label{2dims}
Let $M$ be an irreducible 2-dimensional $\widetilde{H}$-module and assume $q^{2} \neq 1$.

(a) If $M$ has two different weight spaces $M_{t}$ and $M_{t'}$, then $t' = s_{i}t$ for some $i$, and $t(X^{\alpha_{i}}) \neq q^{\pm 2}$ or $1$, but $t(X^{\alpha_{j}}) = q^{\pm 2}$ and $s_{i}t(X^{\alpha_{j}}) = q^{ \pm 2}$ for $j \neq i$. Moreover, there is a unique 2-dimensional module (up to isomorphism) containing these two weight spaces.

(b) If $M$ has only one weight space $M_{t}^{\textrm{gen}}$, then $t(X^{\alpha_{i}}) = 1$ for some $i$, and for $j \neq i$, $t(X^{\alpha_{j}}) = q^{2}$, and either $\langle \alpha_{j}, \alpha_{i}^{\vee} \rangle = 0$ or $q^{2}= -1$ and it is not the case that $\langle \alpha_{i}, \alpha_{j}^{\vee} \rangle = -1$ and $\langle \alpha_{j}, \alpha_{i}^{\vee} \rangle = -2$.

\end{proposition}

\noindent
\textit{Proof.} (a) If $t(X^{\alpha_{i}}) = 1$ for some $i$, then consider $M$ as a $\widetilde{H}_{\{i\}}$-module. By Kato's criterion ( Theorem $\ref{eq:Kato}$), the fact that $q^{2} \neq 1$, and Theorem $\ref{PrinQuot}$c, there is only one irreducible $\widetilde{H}_{\{i\}}$-module $N$ with central character $t$, where $t(X^{\alpha_{i}}) = 1$. This module is 2-dimensional with dim $N_{t}^{\textrm{gen}} = 2$. Thus $M \cong N$ as $\widetilde{H}_{\{1\}}$-modules and $t = t'$.

Assume $M$ has two different weight spaces $M_{t}$ and $M_{t'}$. Then since $M$ is irreducible, some $\tau_{i}$ must be non-zero on $M_{t}$, and $t' = s_{i}t$. Then $\tau_{i}$ must also be non-zero on $M_{s_{i}t}$, and $t(X^{\alpha_{i}}) \neq q^{\pm 2}$. Since $M_{s_{j}t} = 0 = M_{s_{j}t'}$ for $j \neq i$, $t(X^{\alpha_{j}}) = q^{\pm 2}$ and $s_{i}t(X^{\alpha_{j}}) = q^{\pm 2}$. The weight structure determines the action of $\mathbb{C}[X]$ on $M$, and since we know how the operators $\tau_{i}$ act, the actions of the $T_{i}$ are determined as well, so that the module structure of $M$ is determined by its weight structure.

(b) Assume $M$ consists of one generalized weight space $M_{t}^{\textrm{gen}}$, with $v_{t} \in M_{t}$. If all the operators $\tau_{i}$ were defined on $M_{t}$, then $\tau_{i}(v_{t}) = 0$ for all $i$. Hence $T_{i}v_{t} \in \mathbb{C}v_{t}$ for all $i$ and $v_{t}$ would span a submodule of $M$, a contradiction. Thus some $\tau_{i}$ is not well-defined and $t(X^{\alpha_{i}}) = 1$.

If $t(X^{\alpha_{j}}) = 1$ for any $j \neq i$, then $t(X^{\beta}) = 1$ for any $\beta$ in the span of $\alpha_{i}$ and $\alpha_{j}$. Then $M$ is irreducible as a module over $\widetilde{H}_{\{i,j\}}$, the submodule of $\widetilde{H}$ generated by $\mathbb{C}[X]$, $T_{i}$ and $T_{j}$, and must have dimension at least as large as the number of roots in the root subsystem generated by $\alpha_{i}$ and $\alpha_{j}$, which is a contradiction since this number will be greater than 2. Then $t(X^{\alpha_{j}}) \neq 1$ for all $j \neq i$.

Then consider $M$ as an $\widetilde{H}_{\{i\}}$ module, which is irreducible by Theorem $\ref{eq:Kato}$c. The action of $T_{i}$ on the basis $\{v_{t},T_{i}v_{t}\}$ is given by
\renewcommand{\baselinestretch}{1} \normalsize
\[ M(T_{i}) = \begin{bmatrix} 0 & 1 \\ 1 & q-q^{-1} \end{bmatrix}.\]
The action of $\mathbb{C}[X]$ on this subspace is given by \[X^{\lambda}\cdot v = t(X^{\lambda})v \textrm{ and } X^{\lambda} \cdot T_{i}v = \left(T_{i}X^{s_{i}\lambda} + (q-q^{-1})\frac{X^{\lambda}-X^{s_{i}\lambda}}{1-X^{-\alpha_{i}}}\right)v.\]
Specifically, \[ M(X^{\alpha_{j}}) = t(X^{\alpha_{j}}) \begin{bmatrix} 1 & (q-q^{-1})\langle \alpha_{j}, \alpha_{i}^{\vee} \rangle \\ 0 & 1 \end{bmatrix}.\]
Then since $\tau_{j}$ (which is well-defined since $t(X^{\alpha_{j}}) \neq 1$) must be the zero map on $M_{t}^{\textrm{gen}}$, we have $M(T_{j}) = M(\frac{q-q^{-1}}{1-X^{-\alpha_{j}}})$, and
\[M(T_{j})  = (q-q^{-1})\left(\frac{1}{1-t(X^{-\alpha_{j}})}\right) \left[ \begin{array}{cc} 1 & \frac{(q-q^{-1})t(X^{-\alpha_{j}})}{1-t(X^{-\alpha_{j}})}\langle - \alpha_{j}, \alpha_{i}^{\vee} \rangle \\ 0 & 1 \end{array} \right].\]

However, since $T_{j}$ satisfies the relation $T_{j}^{2} = (q-q^{-1})T_{j} + 1$, or $(T_{j} - q)(T_{j}+q^{-1}) = 0$, $M(T_{j})$ must have eigenvalues $q$ or $-q^{-1}$. Hence, since $M(T_{j})$ has $\frac{q-q^{-1}}{1-t(X^{-\alpha_{j}})}$ as each diagonal entry, we have $\frac{q-q^{-1}}{1-t(X^{-\alpha_{j}})} = q$ or $-q^{-1}$, so that $t(X^{-\alpha_{j}}) = q^{\pm 2}$.

If $q^{2} \neq -1$, so that $q \neq -q^{-1}$, then either $M(T_{j})-qI$ or $M(T_{j}) + q^{-1}I$ is invertible, so that the other must actually be zero and so the off-diagonal term must be zero. The only way this can occur is if $\langle \alpha_{j}, \alpha_{i}^{\vee} \rangle = 0. $

If $q^{2} = -1$ then
\renewcommand{\baselinestretch}{1} \normalsize
\[M(T_{j}) = \begin{bmatrix} q & \langle \alpha_{j} ,\alpha_{i}^{\vee} \rangle \\ 0 & q \end{bmatrix}.\]
Then $M(T_{i})$ and $M(T_{j})$ must satisfy the same braid relation as $T_{i}$ and $T_{j}$, which is determined by the type of root system spanned by $\alpha_{i}$ and $\alpha_{j}$. A check of the possible root systems ($A_{1} \times A_{1}$, $A_{2}$, $C_{2}$, and $G_{2}$) shows that the braid relation is satisfied unless $\langle \alpha_{i}, \alpha_{j}^{\vee} \rangle = -1$ and $\langle \alpha_{j}, \alpha_{i}^{\vee} \rangle = -2$. $\square$

\begin{theorem} \label{eq:lemma} (\cite{Ra1}, Lemma 1.19)
Assume $q^{2} \neq 1$. Let $t \in T$ such that $t(X^{\alpha_{i}}) = 1$ and suppose that $M$ is a $\widetilde{H}(q)$-module such that $M_{t}^{\textrm{gen}} \neq 0$. Let $W_{t}$ be the stabilizer of $t$ under the action of $W_{0}$ on $T$. Assume that $\overline{w} \in W_{0}/W_{t}$ is such that $t$ and $\overline{w}t$ are in the same connected component of the calibration graph for $t$, and let $w$ be a minimal length coset representative for $\overline{w}$.Then
\begin{itemize}
\item[(a)] dim$(M_{wt}^{\textrm{gen}}) \geq 2$ and dim$M_{wt}^{\textrm{gen}} > $ dim$M_{wt}$.

\item[(b)] If $M_{s_{j}wt}^{\textrm{gen}} = 0$ then $(\overline{w}t)(X^{\alpha_{j}}) = q^{\pm 2}$ and if, in addition, $q^{2} \neq -1$, then $\langle w^{-1}\alpha_{j},\alpha_{i}^{\vee} \rangle = 0.$
\end{itemize}
\end{theorem}

\vspace{.1in}

\noindent \textbf{Visualizing Modules.} \label{locregsect} For $t \in T$, define \[ Z(t) = \{\alpha \in R^{+} \, | \, t(X^{\alpha}) = 1\} \textrm{ and } P(t) = \{ \alpha \in R^{+}\, | \, t(X^{\alpha}) = q^{\pm 2}\}.\] Notice that $|Z(t)|$ and $|P(t)|$ are constant on orbits $W_{0}t$, since the action of $W_{0}$ permutes the multiset $\{ t(X^{\alpha}) \, | \, \alpha \in R \}$.

The $\tau$ operators and the sets $Z(t)$ and $P(t)$ provide extensive information about the structure (and sometimes the composition factors) of $M(t)$. Let $H_{\alpha}$ be the hyperplane fixed by $s_{\alpha}$ for $\alpha \in R$. A \textit{chamber} is a connected component of $\mathbb{R}^{n} \setminus \cup_{\alpha \in R^{+}} H_{\alpha}$, and $W_{0}$ acts faithfully and simply transitively on the set of chambers. Choose a fundamental chamber $C$ and define the positive side of a hyperplane $H_{\alpha}$ to be the side on which $C$ lies. The map \begin{eqnarray} \label{chamberbij} \{ \textrm{chambers} \} & \leftrightarrow & W_{0} \\ w^{-1}C & \mapsto & w \nonumber \end{eqnarray} is a bijection. In Type $A_{2}$, the picture is:

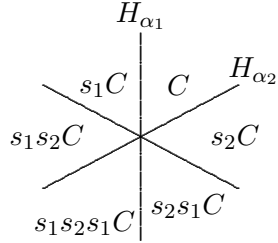
\[ \beginpicture
\setcoordinatesystem units <0.25cm,0.23cm>
\setplotarea x from -25 to 25, y from -5 to 5
\linethickness=3pt
\plot 0 6 0 -6 /
\plot 5.196 3 -5.196 -3 /
\plot 5.196 -3 -5.196 3 /
\put{$C$} at 2 3
\put{$s_{1}C$} at -2 3
\put{$s_{2}C$} at 5 0
\put{$s_{1}s_{2}C$} at -5 0
\put{$s_{2}s_{1}C$} at 2.5 -4
\put{$s_{1}s_{2}s_{1}C$} at -3 -5
\put{$H_{\alpha_{1}}$} at 0 7
\put{$H_{\alpha_{2}}$} at 6 3.8
\endpicture .\]

By 3.15, 4.2, and 5.3 of Steinberg, (\cite{St}) the stabilizer of $t$ is \[ W_{t} = \langle s_{\alpha} | \alpha \in Z(t) \rangle.\]  Thus, if $W_{0}/W_{t}$ is a set of minimal length coset representatives of $W_{t}$-cosets in $W_{0}$, then \begin{equation} \begin{matrix}  \label{chamberbij2} W_{0}/W_{t} & \leftrightarrow & W_{0}t & \leftrightarrow & \{ \textrm{ chambers on the positive side of all } H_{\alpha}, \alpha \in Z(t) \} \\ w & \mapsto & wt & \mapsto  & w^{-1}C \quad, \textrm{ for } w \in W_{0}/W_{t} \end{matrix} \end{equation} are bijections. Again using Type $A_{2}$ as an example, each $W_{0}$-orbit in $T$ has a representative such that the bijection \eqref{chamberbij2} is illustrated by one of the following pictures.

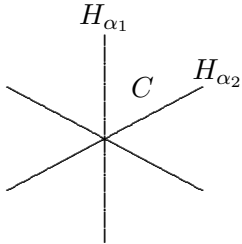
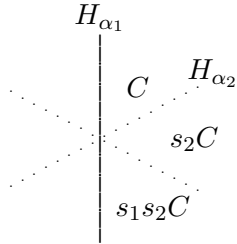
\[ \beginpicture
\setcoordinatesystem units <0.25cm,0.23cm>
\setplotarea x from -25 to 25, y from -5 to 5
\linethickness=3pt
\plot -23 6 -23 -6 /
\plot -17.804 3 -28.196 -3 /
\plot -17.804 -3 -28.196 3 /
\put{$C$} at -21 3
\put{$H_{\alpha_{1}}$} at -23 7
\put{$H_{\alpha_{2}}$} at -17 3.8
\put{$Z(t) = R^{+} $} at -23 -9
\put{$W_{0}/W_{t} = \{ 1 \}$} at -23 -11
\plot -4 6 -4 -6 /
\put{$C$} at -2 3
\put{$s_{2}C$} at 1 0
\put{$s_{1}s_{2}C$} at -1.2 -4
\put{$H_{\alpha_{1}}$} at -4 7
\put{$H_{\alpha_{2}}$} at 2 3.8
\put{$Z(t) = \{ \alpha_{1} \}$} at -2 -9
\put{$W_{0}/W_{t} = \{ 1, s_{2}, s_{1}s_{2}\}$} at -2 -11
\put{$C$} at 22 3
\put{$s_{1}C$} at 18 3
\put{$s_{2}C$} at 25 0
\put{$s_{2}s_{1}C$} at 15 0
\put{$s_{1}s_{2}C$} at 22.5 -4
\put{$s_{1}s_{2}s_{1}C$} at 17 -5
\put{$H_{\alpha_{1}}$} at 20 7
\put{$H_{\alpha_{2}}$} at 26 3.8
\put{$Z(t) = \emptyset $} at 20 -9
\put{$W_{0}/W_{t} = W_{0}$} at 20 -11
\setlinear
\setdots
\plot 20 6 20 -6 /
\plot 25.196 3 14.804 -3 /
\plot 25.196 -3 14.804 3 /
\plot 1.196 3 -9.196 -3 /
\plot 1.196 -3 -9.196 3 /
\endpicture \]

Recall that the elements of the orbit $W_{0}t$ are the vertices of the calibration graph. The hyperplanes $H_{\alpha}$ for $\alpha \in P(t)$ (which are drawn as dashed lines) divide the chambers into subsets corresponding to the components of the calibration graph. The bijection $\eqref{chamberbij2}$ shows that the weights of $M(t)$ are in bijection with the chambers on the positive side of the $H_{\alpha}$ with $\alpha \in Z(t)$. The picture of these chambers can be used to visualize the composition factors of $M(t)$, which is the main goal of this paper.

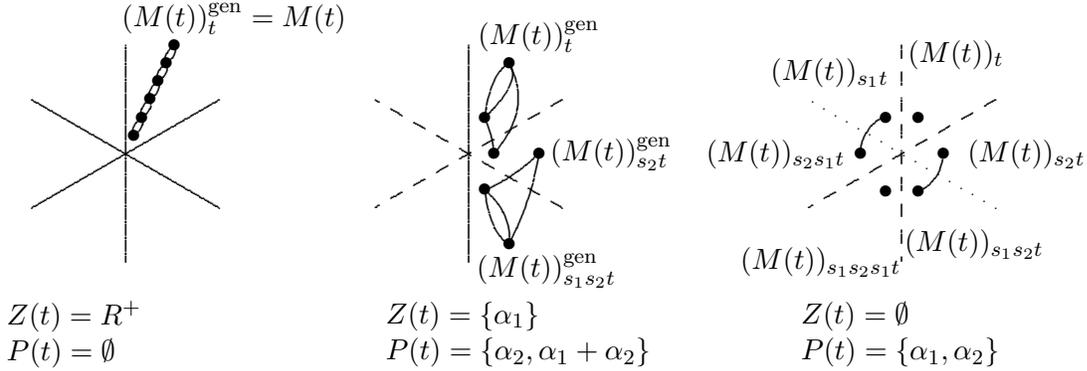
\[ \beginpicture
\setcoordinatesystem units <0.24cm,0.24cm>
\setplotarea x from -27 to 25, y from -5 to 6
\linethickness=3pt
\plot -23 6 -23 -6 /
\plot -17.804 3 -28.196 -3 /
\plot -17.804 -3 -28.196 3 /
\put{$\bullet$} at -21.2 4
\put{$\bullet$} at -21.65 3
\put{$\bullet$} at -22.1 2
\put{$\bullet$} at -20.75 5
\put{$\bullet$} at -20.3 6
\put{$\bullet$} at -22.55 1
\put{$(M(t))_{t}^{\textrm{gen}} = M(t)$} at -17 7.5
\put{\parbox{2in}{$Z(t) = R^{+} \\ P(t) = \emptyset$}} at -19 -10
\plot -4 6 -4 -6 /
\put{$\bullet$} at -2.6 0
\put{$\bullet$} at -0.1 0
\put{$\bullet$} at -3.1 2
\put{$\bullet$} at -1.75 5
\put{$\bullet$} at -3.1 -2
\put{$\bullet$} at -1.75 -5
\put{$(M(t))_{t}^{\textrm{gen}}$} at -0.1 6.5
\put{$(M(t))_{s_{2}t}^{\textrm{gen}}$} at 3.9 0
\put{$(M(t))_{s_{1}s_{2}t}^{\textrm{gen}}$} at 0.3 -6.5
\put{\parbox{2in}{$Z(t) = \{\alpha_{1}\} \\ P(t) = \{\alpha_{2}, \alpha_{1} + \alpha_{2}\}$}} at 2 -10
\put{$\bullet$} at 17.7 0
\put{$\bullet$} at 22.3 0
\put{$\bullet$} at 19.1 2
\put{$\bullet$} at 20.9 2
\put{$\bullet$} at 19.1 -2.1
\put{$\bullet$} at 20.9 -2.1
\put{$(M(t))_{t}$} at 23 5.5
\put{$(M(t))_{s_{2}t}$} at 26.9 0
\put{$(M(t))_{s_{1}s_{2}t}$} at 24 -5
\put{$(M(t))_{s_{1}t}$} at 16 4.5
\put{$(M(t))_{s_{2}s_{1}t}$} at 13 0
\put{$(M(t))_{s_{1}s_{2}s_{1}t}$} at 15.5 -6
\put{\parbox{2in}{$Z(t) = \emptyset \\ P(t) = \{\alpha_{1}, \alpha_{2}\}$}} at 25 -10
\setquadratic
\plot -22.55 1 -22.55 1.5 -22.1 2 /
\plot -22.1 2  -22.1 2.5 -21.65 3 /
\plot -21.65 3 -21.65 3.5 -21.2 4 /
\plot -21.2 4 -21.2 4.5 -20.75 5 /
\plot -20.75 5 -20.75 5.5 -20.3 6 /
\plot -22.55 1 -22.1 1.5 -22.1 2 /
\plot -22.1 2  -21.65 2.5 -21.65 3 /
\plot -21.65 3 -21.2 3.5 -21.2 4 /
\plot -21.2 4 -20.75 4.5 -20.75 5 /
\plot -20.75 5 -20.3 5.5 -20.3 6 /
\plot -3.1 2 -2.6 4 -1.75 5 /
\plot -3.1 2 -1.7 3.8 -1.75 5 /
\plot -3.1 -2 -2.6 -4 -1.75 -5 /
\plot -3.1 -2 -2 -3.4 -1.75 -5 /
\plot -3.1 2 -2.75 1 -2.6 0 /
\plot -1.75 5 -1.3 2.5 -2.6 0 /
\plot -1.75 -5 -0.7 -2.5 -0.1 0 /
\plot -3.1 -2 -1.6 -1.2 -0.1 0 /
\setquadratic
\circulararc 60 degrees from 19.1 2.1 center at 20 0
\circulararc 60 degrees from 20.9 -2.1 center at 20 0
\setlinear
\setdashes
\plot 25.196 3 14.804 -3 /
\plot 20 6 20 -6 /
\plot 1.196 3 -9.196 -3 /
\plot 1.196 -3 -9.196 3 /
\setdots
\plot 25.196 -3 14.804 3 /
\endpicture \]

Each dot in a chamber represents one dimension of the generalized weight space corresponding to that chamber under the bijection \eqref{chamberbij2}. Viewing the dots as vertices of a graph, the connected components of the graph correspond to irreducible modules with central character $t$. In the first picture, $M(t)$ is irreducible, with dim $M(t)_{t}^{\textrm{gen}} = 6$, while in the second picture, $M(t)$ has two composition factors, each 3-dimensional. One has dim $M_{t}^{\textrm{gen}} = 2$ and dim $M_{s_{2}t} = 1$, while the other has dim$M_{s_{1}s_{2}t}^{\textrm{gen}} = 2$ and dim $M_{s_{2}t} = 1$. In the third picture, $M(t)$ has four composition factors - two 1-dimensional modules with weights $t$ and $s_{1}s_{2}s_{1}t$, and two 2-dimensional modules, one with dim $M_{s_{1}t} = 1 = $ dim $M_{s_{2}s_{1}t}$, and one with dim $M_{s_{2}t} = 1 = $ dim $M_{s_{1}s_{2}t}$. These claims are proven by Theorem \ref{calibsame}b and Lemma \ref{eq:lemma}, and these pictures are a good way to visualize those results.

\vspace{.1in}

\noindent \textbf{Calibrated Modules and Weights.} A weight $t$ is defined to be \emph{regular} if $W_{t}$, the
subgroup of the Weyl group that fixes $t$, is trivial. Then a weight $t$ is regular if and only if $Z(t) = \emptyset$.

A representation $M$ is \textit{calibrated} if
\[M_{t}^{\textrm{gen}} = M_{t}\] for all weights $t$, i.e. the subalgebra $\mathbb{C}[X] \subseteq \widetilde{H}$ acts diagonally on $M$.

\begin{proposition} (\cite{Ra1}, Proposition 1.10) \label{calibfund}
(a) If $\, q^{2} \neq 1$, an irreducible $\widetilde{H}$-module is calibrated if and only if

\renewcommand{\baselinestretch}{1} \normalsize
\[ \textrm{dim} (M_{t}^{\textrm{gen}}) = 1\] \renewcommand{\baselinestretch}{1} \normalsize
for all weights $t$ of $M$.

(b) If $M$ is an $\widetilde{H}$-module with regular central character,
then $M$ is calibrated.
\end{proposition}

When $q^{2} = 1$, all irreducible modules are calibrated, as will be shown in section \ref{clifford}.

\vspace{.1in}

\noindent \textbf{Calibrated Modules with regular Central Character.} \begin{theorem} (See \cite{Ra1}, Prop. 1.11) \label{eq:locreg} Assume $q^{2} \neq 1$. Let $t$ be a regular central character, and let $C$ be a component of the calibration graph. Define

\[ \widetilde{H}^{(t,C)} = \mathbb{C}\textrm{-span}\{v_{w} | wt \in C \}\]

Then the vector space $\widetilde{H}^{(t,C)}$ is an irreducible calibrated
$\widetilde{H}$-module with action

\renewcommand{\baselinestretch}{1} \normalsize
\[ X^{\lambda} \cdot v_{w} = (wt)(X^{\lambda})v_{w} \textrm{ for }
X^{\lambda} \in X, w \in W_{0} \textrm{, and} \]

\[ T_{i} \cdot v_{w} = (T_{i})_{w}v_{w} + (q^{-1} +
(T_{i})_{w})v_{s_{i}w} \textrm{ for } 1 \leq i \leq n , w \in W_{0}, \]

where $(T_{i})_{w} = \frac{q - q^{-1}}{1-wt(X^{-\alpha_{i}})}$, and
$v_{s_{i}w} = 0$ if $s_{i}wt \notin C$.
\end{theorem}

Note that since $t$ is a regular central character, $wt(X^{\alpha_{i}}) \neq 1$ for $w \in \mathcal{F}^{(t,C)}$ and $i = 1, \ldots n$. Hence $(T_{i})_{w}$ is well-defined for $w \in \mathcal{F}^{(t,J)}$. The most difficult part of this theorem is checking that the given $\widetilde{H}$-module structure satisfies the braid relation. Since $t$ is assumed to be regular, this essentially follows from the braid relation on the $\tau_{i}$. (See \cite{Ra1} for details.) In fact, more is true.

\begin{theorem} (See \cite{Ra1}, Prop. 1.11)
Assume $q^{2} \neq 1$, and let $M$ be an irreducible $\widetilde{H}$-module with regular central character $t$. ($M$ is therefore calibrated). Then if $wt$ is a weight of $M$, let $C$ be the component of the calibration graph containing $wt$. Then the weights of $M$ are exactly the vertices in $C$. In addition, $M$ is isomorphic to the module $\widetilde{H}^{(t,C)}$ given in Theorem $\ref{eq:locreg}$.
\end{theorem}

\noindent \textit{Proof.} By theorem $\ref{calibsame}$b, the weights of $M$ are exactly $C_{1} \cup C_{2} \cup \ldots \cup C_{n}$ for some components $C_{i}$ of the calibration graph. We prove that if $t$ and $s_{i}t$ are weight spaces of $M$, then $\tau_{i}: M_{t} \rightarrow M_{s_{i}t}$ is a bijection.

To see this, let $v_{t}$ be a nonzero vector in $M_{t}$. Then since $M$ is irreducible, there is a series of $\tau$ operators so that \[ \tau_{i_{1}}\tau_{i_{2}}\ldots \tau_{i_{m}}v_{t} = cv_{s_{i}t}, (c \in \mathbb{C}, c \neq 0) \] a non-zero element of $M_{s_{i}t}$, and we may assume that $m$ is minimal such that this is true. Then $\tau_{i}\tau_{i_{1}} \ldots \tau_{i_{m}}v_{t} \in M_{t}$. Since $t$ is regular, $s_{i}s_{i_{1}} \ldots s_{i_{m}} = 1$, so that this word in $W$ is not reduced. Thus, since the $\tau_{i}$ satisfy the same braid relation as the $s_{i}$, and $\tau_{i}^{2}$ is a multiple of the identity, we can ``reduce'' $\tau_{i_{1}}\ldots \tau_{i_{m}}$ exactly when we can reduce $s_{i_{1}}\ldots s_{i_{m}}$. Thus, by the minimality of $m$, we must have $s_{i_{1}}\ldots s_{i_{m}} = s_{i}$, and $m=1$.

Then $\tau_{i}:M_{t} \rightarrow M_{s_{i}t}$ is a bijection since $v_{s_{i}t}$ is non-zero. Thus, if $t$ is a weight of $M$ and $t(X^{\alpha_{i}}) = q^{\pm 2}$, $s_{i}t$ is not a weight of $M$. Then the weights of $M$ are from a single component of the calibration graph. Then $M$ and $\widetilde{H}^{(t,C)}$ have the same weight spaces, and are both composition factors of $M(t)$. Since every weight space of $M(t)$ is 1-dimensional, different composition factors of $M(t)$ have distinct weight spaces. Thus $M \cong \widetilde{H}^{(t,C)}$. $\square$

\vspace{.1in}

\noindent \textbf{Clifford Theory when $q^{2} = 1$.} \label{clifford} Let $q^{2} = 1$. Then we can identify the subalgebra $H$ spanned by $\{ T_{w} \, | \, w \in W_{0} \}$ with $\mathbb{C}[W_{0}]$, so that \[ \widetilde{H} = \textrm{span}\{ wX^{\lambda} \, | \, w \in W_{0}, \lambda \in P \}.\]

Let $M$ be a finite dimensional simple $\widetilde{H}$-module and let $t \in T$ such that $M_{t} \neq 0$. Let $W_{t}$ be the stabilizer of $t$ in $W_{0}$. As vector spaces, $M_{t} \cong M_{wt}$ via the map $m \mapsto wm$, and \[ M = \bigoplus_{w \in W_{0}/W_{t}} M_{wt},\] since $M$ is simple and the right-hand side is a submodule of $M$. (This implies that all $\widetilde{H}$ modules are calibrated.)

\begin{theorem} (See also \cite{K2}.)
Let $q^{2} = 1$ and let $M$ be an irreducible $\widetilde{H}$-module. Let $t \in T$ be such that $M_{t} \neq 0$. Define $\widetilde{H}_{t} = \{ wX^{\lambda} \, | \, w \in W_{t}, \lambda \in P \}$. \\
Then \begin{enumerate} \item[(a)] $M_{t}$ is an irreducible $W_{t}$-module. \item[(b)] $M_{t}$ is an $\widetilde{H}_{W_{t}}$-module and \[M \cong \widetilde{H}  \otimes_{\widetilde{H}_{W_{t}}} M_{t}.\] \end{enumerate}
\end{theorem}

\noindent \textit{Proof.} (a) Let $q^{2} = 1$ and let $M$ be an irreducible $\widetilde{H}$-module. Let $t \in T$ be such that $M_{t} \neq 0$. If $w \in W_{t}$ and $m \in M_{t}$, then $w \cdot M_{t} = M_{wt}  = M_{t}$, so that $M_{t}$ is a $W_{t}$-module.

Let $N$ be a nontrivial $W_{t}$-submodule of $M_{t}$. If $I$ is a set of coset representatives of $W_{0} / W_{t}$, then $M = \bigoplus_{w \in I} wN$ since the right-hand side is a submodule of $M$ and $M$ is simple. Thus $N = M_{t}$ and $M_{t}$ is a simple $W_{t}$-module.

(b) Since $\mathbb{C}[X] M_{t} = M_{t}$, $M_{t}$ is an $\widetilde{H}_{W_{t}}$-module. The set $\{ T_{w} \, | \, w \in I \}$ is a basis for $\widetilde{H}$ over $\widetilde{H}_{W_{t}}$ and then \renewcommand{\baselinestretch}{1} \normalsize \begin{eqnarray*} \phi: \widetilde{H}  \otimes_{\widetilde{H}_{W_{t}}} M_{t} & \rightarrow & M \\ T_{w} \otimes m & \mapsto & T_{w} \cdot m \textrm{,  for } w \in I \end{eqnarray*} is an $\widetilde{H}$-module homomorphism. Since $M$ is simple, the map is surjective. Since both modules have dimension $|W:W_{t}| \cdot (\textrm{dim}(M_{t}))$, the map is injective as well. $\square$

\section{Type $A_{1}$} \label{A1start}

We begin with the type $A_{1}$ affine Hecke algebra, restating definitions explicitly to give a clear example of the constructions above, as well as our approach to classifying simple $\widetilde{H}$-modules.

\vspace{.1in}

\noindent \textbf{The Affine Hecke Algebra} The type $A_{1}$ affine Hecke algebra is built on the root data of SL$_{2}$. Let \[ P = \mathbb{Z}\omega_{1} \quad \textrm{and} \quad X = \{X^{k\omega_{1}} | k \in \mathbb{Z} \}\] so that $X$ is the group generated by $X^{\omega_{1}}$ and is isomorphic to $P$. The Weyl group is $W_{0} = \{1, s_{1} \}$ with $s_{1}^{2} = 1$, and setting $ s_{1}X^{\omega_{1}} = X^{-\omega_{1}}$ defines an action of $W_{0}$ on $X$. Let $q \in \mathbb{C}^{\times}$. The \textit{affine Hecke algebra} of type $A_{1}$ is

\[ \widetilde{H} = \mathbb{C}\textrm{-span}\{X^{k\omega_{1}},T_{1}X^{k\omega_{1}} \,|\,\, k \in \mathbb{Z} \} ,\label{A1affhecke}\] with relations
\begin{equation} T_{1}^{2} = (q-q^{-1})T_{1} + 1,  \quad X^{k\omega_{1}} X^{m\omega_{1}} = X^{(k+m)\omega_{1}}, \quad X^{\omega_{1}}T_{1} = T_{1}X^{-\omega_{1}} + (q-q^{-1})X^{\omega_{1}}. \label{A1rels} \end{equation}

A \textit{weight} $t$ is an element of \[T =  \textrm{Hom}_{\mathbb{C}\textrm{-alg}}(\mathbb{C}[X],\mathbb{C}) = \textrm{Hom}_{\textrm{Gp}}(X, \mathbb{C}^{\times}) \label{TDef}.\] We let $t_{z}$ denote the weight given by $t_{z}(X^{\omega_{1}}) = z$. This definition is a special case of the definition given in section \ref{subsectHecke}, using the root system $R = \{ \pm \alpha_{1} \}$, where $\alpha_{1} = 2\omega_{1}$.

\begin{proposition}\label{eq:A11dims}
There are four 1-dimensional $\widetilde{H}$ representations

\vspace{.2in}
$\begin{array}{cccc} L_{q}: \widetilde{H} \rightarrow \mathbb{C} & \qquad L_{-q}: \widetilde{H} \rightarrow \mathbb{C} & \qquad L_{-q^{-1}}: \widetilde{H} \rightarrow \mathbb{C} & \qquad L_{q^{-1}}: \widetilde{H} \rightarrow \mathbb{C} \\
 X^{\omega_{1}} \mapsto q & \qquad X^{\omega_{1}} \mapsto -q & \qquad X^{\omega_{1}} \mapsto -q^{-1} & \qquad X^{\omega_{1}} \mapsto q^{-1} \\ T_{1} \mapsto q   & \qquad T_{1} \mapsto q & \qquad T_{1} \mapsto -q^{-1} & \qquad T_{1} \mapsto -q^{-1} \end{array}$

\vspace{.2in}

If $q^{2} = -1$, then $q = -q^{-1}$, so that $L_{q} = L_{-q^{-1}}$ and $L_{-q} = L_{q^{-1}}$.

\end{proposition}

\noindent \textit{Proof.} A straightforward check of the relations \eqref{A1rels} shows that the maps $L_{\pm q^{\pm 1}}$ are homomorphisms. Assume $M$ is a 1-dimensional $\widetilde{H}$-module, with $M = \textrm{span}\{v\}$. By $\eqref{A1rels}$,\[T_{1}v = qv \,\, \textrm{  or  } \,\, T_{1}v = -q^{-1}v.\] Then if $T_{1}v = qv$, then again by \eqref{A1rels}, $X^{2\omega_{1}}v = q^{2}v$ and either $M \cong L_{q}$ or $M \cong L_{-q}$. Similarly, if $T_{1}v = -q^{-1}v,$ then $X^{2\omega_{1}}v = q^{-2}v$ and either $M \cong L_{q^{-1}}$ or $M \cong L_{-q^{-1}}. \\ \square$

\vspace{.2in}

Let $t \in T$ and let $\mathbb{C}_{t} = \textrm{span}\{v_{t}\}$ be the one-dimensional $\mathbb{C}[X]$-module given by

\[ X^{\lambda}v_{t} = t(X^{\lambda})v_{t}.\] Then the principal series module is

\[ M(t) = \textrm{Ind}_{\mathbb{C}[X]}^{\widetilde{H}} \mathbb{C}_{t} = \widetilde{H} \otimes_{\mathbb{C}[X]} \mathbb{C}_{t} = \textrm{span}\{v_{t},T_{1}v_{t}\},\] so that $\{ v_{t}, T_{1}v_{t}\}$ is a basis for $M(t)$.

\renewcommand{\baselinestretch}{1} \normalsize

\begin{proposition} \label{A1thm1}

\begin{enumerate}
\item[(a)] If $t(X^{\omega_{1}}) \neq \pm q^{\pm 1}$, then $M(t)$ is irreducible.
\item[(b)] If $\, t(X^{\omega_{1}}) = \pm q^{\pm 1}$, $M(t)$ has composition series $M(t) \supseteq M_{1} \supseteq 0$.

If $\, t(X^{\omega_{1}}) = \pm q$, then
 \[ M_{1} = \textrm{  span}\{T_{1}v_{t} - qv_{t}\} \cong L_{\pm q^{-1}} \textrm{  and } \quad M(t) / M_{1} \cong  L_{\pm q}. \]  If $\, t(X^{\omega_{1}}) = \pm q^{-1}$,
 \[ M_{1}  =  \textrm{  span}\{T_{1}v_{t} + q^{-1}v_{t}\} \cong L_{\pm q} \textrm{  and } \quad M(t) / M_{1} \cong L_{\pm q^{-1}}. \]

\item[(c)] If $q^{2} = 1$ and $t(X^{\omega_{1}}) = \pm q$, then $M(t) = M^{+} \oplus M^{-}$, where \[ M^{+} = \textrm{span}\{T_{1}v_{t} - qv_{t}\}, \, \, \textrm{ and } \quad M^{-} = \textrm{span}\{T_{1}v_{t}+q^{-1}v_{t}\}.\]
\end{enumerate}
\end{proposition}

\renewcommand{\baselinestretch}{1} \normalsize

\noindent \textit{Proof.} Part a is given by Kato's criterion (Theorem \ref{eq:Kato}), which also shows that $M(t)$ must be reducible if $t(X^{\omega_{1}}) = \pm q^{\pm 1}$. In parts b and c, it is easy to check from the defining relations for $\widetilde{H}$ that the given vectors span submodules. In part b, $M(t)$ cannot be a direct sum since it would have to be the direct sum of its two weight spaces. But $M(t)_{t}$ is spanned by $v_{t}$ and is not a submodule. $\square$

Note that if $q$ is a primitive fourth root of unity, the two cases in part (b) of this theorem coincide. To visualize this classification, identify $\{ t_{q^{x}} \,|\, x \in \mathbb{R} \}$ with the real line. In this picture, the set \[ H_{\alpha} = \{x \, |\, t_{q^{x}}(X^{\alpha_{1}}) = 1 \} \] is marked with a solid line, while \[ H_{\alpha + \delta} = \{x \, |\, t_{q^{x}}(X^{\alpha_{1}}) = q^{2}\}
\,\, \textrm{ and } \,\, H_{\alpha  -\delta } = \{ x \, | \, t_{q^{x}}(X^{\alpha_{1}}) = q^{-2} \} \] are denoted by dashed lines.

\[\beginpicture
\setcoordinatesystem units <0.25cm,0.25cm>
\setplotarea x from -25 to 25, y from -5 to 5
\linethickness=3pt
\plot -20 0 20 0 /
\plot 0 2 0 -2 /
\put{$t_{q^{0}}$} at 0 -3
\put{$t_{q^{2}}$} at 10.5 -2
\put{$H_{\alpha}$} at 0 3
\put{$H_{\alpha + \delta}$} at 5 3
\put{$H_{\alpha - \delta}$} at -5 3
\setdashes
\plot -5 2 -5 -2 /
\plot 5 2 5 -2 /
\put{$\bullet$} at 15 0
\put{$\bullet$} at -5 0
\put{$\bullet$} at 0 0
\put{$\bullet$} at 5 0
\put{$t_{q}$} at 5.8 -3
\put{$t_{q^{-1}}$} at -4.2 -3
\put{$t_{q^{3}}$} at 15.8 -2
\put{$\bullet$} at 10 0
\put{Characters $t_{q^{x}}$, generic $q$.} at 0 -7
\endpicture \]

If $q$ is a primitve $2\ell$th root of unity then $\{t_{q^{x}} | x \in \mathbb{R} \}$ is identified with $\mathbb{R} / 2\ell \mathbb{Z}$ and $H_{\alpha} = \{ k\ell \, | \, k \in \mathbb{Z} \}$. The following is the specific case $\ell = 2$, so that $t_{1} = t_{q^{2}} = \ldots$. The periodicity is evident in the picture.

\[\beginpicture
\setcoordinatesystem units <0.25cm,0.25cm>
\setplotarea x from -25 to 25, y from -5 to 5
\linethickness=3pt
\plot -22 0 22 0 /
\plot 0 2 0 -2 /
\plot 20 2 20 -2 /
\plot 10 2 10 -2 /
\plot -10 2 -10 -2 /
\plot -20 2 -20 -2 /
\put{$t_{q^{0}}$} at 0 -3
\put{$t_{q^{2}}$} at 10.5 -3
\put{$t_{q^{0}}$} at 20.5 -3
\put{$t_{q^{3}}$} at 15.5 -3
\put{$t_{q^{2}}$} at -9.5 -3
\put{$t_{q^{3}}$} at -4.5 -3
\put{$t_{q^{1}}$} at -14.5 -3
\put{$t_{q^{0}}$} at -19.5 -3
\put{$H_{\alpha}$} at 0 3
\put{$H_{\alpha}$} at 10 3
\put{$H_{\alpha}$} at 20 3
\put{$H_{\alpha}$} at -10 3
\put{$H_{\alpha}$} at -20 3
\put{$H_{\alpha \pm \delta}$} at 5 3
\put{$H_{\alpha \pm \delta}$} at 15 3
\put{$H_{\alpha \pm \delta}$} at -5 3
\put{$H_{\alpha \pm \delta}$} at -15 3
\setdashes
\plot -5 2 -5 -2 /
\plot 5 2 5 -2 /
\plot 15 2 15 -2 /
\plot -15 2 -15 -2 /
\put{$\bullet$} at 0 0
\put{$\bullet$} at 15 0
\put{$\bullet$} at -5 0
\put{$\bullet$} at -10 0
\put{$\bullet$} at -15 0
\put{$\bullet$} at -20 0
\put{$\bullet$} at 5 0
\put{$t_{q}$} at 5.8 -3
\put{$\bullet$} at 10 0
\put{$\bullet$} at 20 0
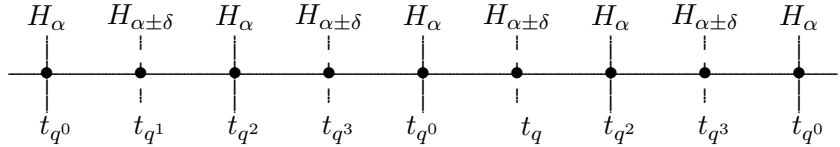
\put{Characters $t_{q^{x}}$, $q^{4} =1$.} at 0 -7
\endpicture \]

\noindent
If $q^{2} = 1$, then $t_{1} = t_{q} = t_{q^{2}} = \ldots$.

\[\beginpicture
\setcoordinatesystem units <0.25cm,0.25cm>
\setplotarea x from -25 to 25, y from -5 to 5
\linethickness=3pt
\plot -22 0 22 0 /
\plot 0 2 0 -2 /
\plot 20 2 20 -2 /
\plot 10 2 10 -2 /
\plot -10 2 -10 -2 /
\plot -20 2 -20 -2 /
\plot -5 2 -5 -2 /
\plot 5 2 5 -2 /
\plot 15 2 15 -2 /
\plot -15 2 -15 -2 /
\put{$t_{q^{0}}$} at 0 -3
\put{$t_{q^{0}}$} at 10.5 -3
\put{$t_{q^{0}}$} at 20.5 -3
\put{$t_{q}$} at 15.5 -3
\put{$t_{q^{0}}$} at -9.5 -3
\put{$t_{q}$} at -4.5 -3
\put{$t_{q}$} at -14.5 -3
\put{$t_{q^{0}}$} at -19.5 -3
\put{$H_{\alpha}$} at 0 3
\put{$H_{\alpha}$} at 10 3
\put{$H_{\alpha}$} at 20 3
\put{$H_{\alpha}$} at -10 3
\put{$H_{\alpha}$} at -20 3
\put{$H_{\alpha}$} at 5 3
\put{$H_{\alpha}$} at 15 3
\put{$H_{\alpha}$} at -5 3
\put{$H_{\alpha}$} at -15 3
\setdashes
\put{$\bullet$} at 20 0
\put{$\bullet$} at 0 0
\put{$\bullet$} at 15 0
\put{$\bullet$} at -5 0
\put{$\bullet$} at -10 0
\put{$\bullet$} at -15 0
\put{$\bullet$} at -20 0
\put{$\bullet$} at 5 0
\put{$t_{q}$} at 5.8 -3
\put{$\bullet$} at 10 0
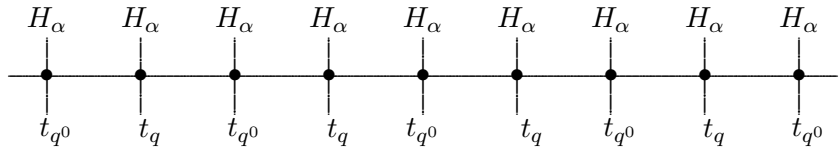
\put{Characters $t_{q^{x}}$, $q^{2} =1$.} at 0 -7
\endpicture \]

The structure of $M(t)$ can be seen using the $\tau$ operator described above (equation \eqref{taudef}). If $t(X^{\alpha_{1}}) = 1$, then $\tau_{1}$ is not defined. If $t(X^{\alpha_{1}}) = q^{\pm 2} \neq 1$, then $\tau_{1}$ is non-zero on $M(t)_{t}$, but is zero on $M(t)_{s_{1}t}$. This accounts for the structure of $M(t)$ described above (Proposition \ref{A1thm1} (b)). If $t(X^{\omega_{1}}) \neq q^{2}$ or 1, then $\tau_{1}$ is invertible on $M(t)$. The following pictures show the chambers around $t$ as described in $\eqref{locregsect}$, which give a picture of $M(t)$.

Each neighborhood is divided into two chambers, representing the two weights in the $W_{0}$-orbit of the central character. Each chamber $wC$ contains a number of large dots equal to the dimension of $M(t)_{wt}^{\textrm{gen}}$, the generalized weight space with weight $wt$, and lines connecting the dots show that the corresponding basis vectors are in the same composition factor. A solid line is drawn where $\tau_{1}$ is not defined, as in the following picture.

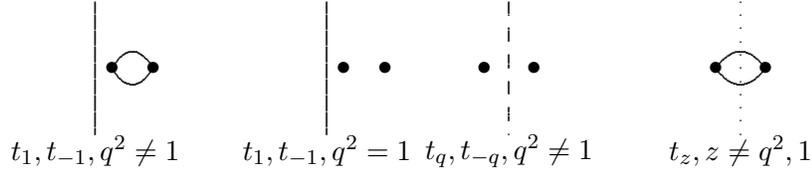
\[ \beginpicture
\setcoordinatesystem units <0.22cm,0.22cm>
\setplotarea x from -25 to 25, y from -7 to 7
\linethickness=3pt
\plot -22 4 -22 -4 /
\put{$\bullet$} at -21 0
\put{$\bullet$} at -18.5 0
\put{$t_{1},t_{-1}, q^{2} \neq 1$} at -22 -5
\plot -8 4 -8 -4 /
\put{$\bullet$} at -7 0
\put{$\bullet$} at -4.5 0
\put{$t_{1},t_{-1}, q^{2} = 1 $} at -8 -5
\setquadratic
\plot -21 0 -19.75 1 -18.5 0 /
\plot -21 0 -19.75 -1 -18.5 0 /
\linethickness=3pt
\put{$\bullet$} at 4.5 0
\put{$\bullet$} at 1.5 0
\put{$\bullet$} at 15.5 0
\put{$\bullet$} at 18.5 0
\put{$t_{q},t_{-q}, q^{2} \neq 1 $} at 3 -5
\put{$t_{z}, z \neq q^{2},1$} at 17 -5
\setquadratic
\plot 15.5 0 17 1 18.5 0 /
\plot 15.5 0 17 -1 18.5 0 /
\setlinear
\setdashes
\plot 3 4 3 -4 /
\setdots
\plot 17 4 17 -4 /
\endpicture \]

The first two pictures represent weights with $t(X^{\alpha_{1}}) = 1$, so that $s_{1}t = t$. Hence the two dots representing the weight basis vectors are drawn in the same chamber. If $q^{2} \neq 1$, $M(t)$ is irreducible, so the basis vectors are connected.
The next two pictures show $M(t)$ for the other possible weights $t$. A dashed line denotes a weight $t$ with $t(X^{\alpha_{1}}) = q^{\pm 2}$, so that the operator $\tau_{1}$ on the weight spaces of $M(t)$ is well-defined, but not invertible in both directions. This means that the corresponding weight spaces will be in different composition factors of $M(t)$. A dotted line denotes a weight $t$ with $t(X^{\alpha_{1}}) \neq q^{2}$ or 1. In this case $\tau_{1}$ is invertible, forcing the weight spaces to be in the same composition factor of $M(t)$. Accordingly, the dots representing them are connected.

Notice also the connection to the drawings of central characters above. The pictures of the $M(t_{q^{x}})$ are a picture of a small open neighborhood around the point $t_{q^{x}}$ in the picture of the characters. This correspondence is not as clear in this case as in others, since it is the smallest example of the affine Hecke algebra, but the essential ingredients are present. The weights of the $M(t)$ are all displayed, as are the actions of the $\tau$ operators that determine the composition factors of $M(t)$. The complete classification of $\widetilde{H}$ modules is summarized in the following tables.

\begin{table}[htpb]
\centering
$\begin{array}{| c | c | c | c |}
\hline \, & \multicolumn{3}{|c|}{\textrm{Dims. of Irreds.}} \\
\hline t & q^{4} \neq 1 & q^{2} = -1 & q = -1 \\
\hline t_{1}  & 2 & 2 & 1,1  \\
\hline t_{-1}  & 2 & 2 & 1,1 \\
\hline t_{q}   & 1,1 & 1,1 & N/A \\
\hline t_{-q}  & 1,1 & N/A & N/A \\
\hline t_{z}, z  \neq  \pm 1 \textrm{ or } \pm q  & 2 & 2 & 2 \\
\hline
\end{array}$
\caption{Table of possible central characters in Type $A_{1}$, with $q$ generic.}
\end{table}

It should be noted that for any value of $q$ with $q^{2} \neq \pm 1$, (that is, $q^{4} \neq 1$) the representation theory of $\widetilde{H}$ can be described in terms of $q$ only. If $q^{2} = \pm 1$, then the representation theory of $\widetilde{H}$ does not fit that same description. This fact can be seen through a number of different lenses. It is a reflection of the fact that the sets $P(t)$ and $Z(t)$ for all possible central characters $t$ can be described solely in terms of $q$. In the local region pictures, this is reflected in the fact that the hyperplanes $H_{\alpha}$ and $H_{\alpha \pm \delta}$ are distinct \textit{unless} $q^{2} = \pm 1$. When these hyperplanes coincide, the sets $P(t)$ and $Z(t)$ change for characters on those hyperplanes.

\section{Type $A_{2}$}

The type $A_{2}$ root system is $R = \{ \pm \alpha_{1}, \pm \alpha_{2}, \pm(\alpha_{1} + \alpha_{2})\},$ where $\langle \alpha_{1}, \alpha_{2}^{\vee} \rangle = -1 = \langle \alpha_{2}, \alpha_{1}^{\vee} \rangle$, with Weyl group $W_{0} \cong S_{3}$. The simple roots are $\alpha_{1}$ and $\alpha_{2}$, and $\alpha_{1} + \alpha_{2}$ is the only other positive root.

\[ \beginpicture
\setcoordinatesystem units <0.4cm,0.4cm>
\setplotarea x from -15 to 15, y from -5 to 6
\linethickness=5pt
\setlinear
\thicklines
\arrow <10 pt> [0.2,0.67] from 0 0 to -6 0
\arrow <10 pt> [0.2,0.67] from 0 0 to 6 0
\arrow <10 pt> [0.2,0.67] from 0 0 to -3 -5.196
\arrow <10 pt> [0.2,0.67] from 0 0 to -3 5.196
\arrow <10 pt> [0.2,0.67] from 0 0 to 3 -5.196
\arrow <10 pt> [0.2,0.67] from 0 0 to 3 5.196
\arrow <10 pt> [0.2,0.67] from 0 0 to 3 1.732
\arrow <10 pt> [0.2,0.67] from 0 0 to 0 3.468
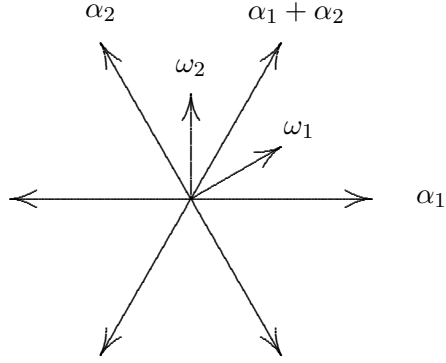
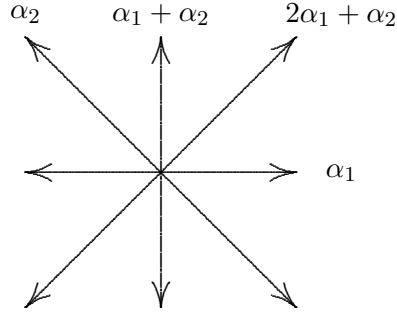
\put{The type $A_{2}$ root system} at 0 -7
\put{$\alpha_{1}$} at 8 0
\put{$\alpha_{2}$} at -3 6.196
\put{$\alpha_{1} + \alpha_{2}$} at 3.5 6.196
\put{$\omega_{1}$} at 3.6 2.3
\put{$\omega_{2}$} at 0 4.368
\endpicture\]

In this picture, $s_{i}$ is reflection through the hyperplane perpendicular to $\alpha_{i}$. The fundamental weights satisfy $$ \begin{matrix} \omega_{1} = \frac{1}{3}(2\alpha_{1} + \alpha_{2}) \\ \omega_{2} = \frac{1}{3}(2\alpha_{2}+\alpha_{1}) \end{matrix} \qquad \begin{matrix} \alpha_{1} =  2\omega_{1} - \omega_{2} \\ \alpha_{2} = 2\omega_{2} - \omega_{1} \end{matrix}.$$ Let \[P = \mathbb{Z}\textrm{-span}\{\omega_{1},\omega_{2}\} \textrm{ and } Q = \mathbb{Z}\textrm{-span}(R) \] be the \textit{weight lattice} and \textit{root lattice} of $R$, respectively.

The affine Hecke algebra $\widetilde{H}$ is defined as in section \ref{subsectHecke}. Let \[\mathbb{C}[Q] = \{ X^{\lambda} \, | \, \lambda \in Q\} \textrm{ and } T_{Q} = \textrm{Hom}_{\mathbb{C}\textrm{-alg}}(\mathbb{C}[Q],\mathbb{C}).\] Define \[t_{z,w}: \mathbb{C}[Q] \rightarrow \mathbb{C} \textrm{ by  } t_{z,w}(X^{\alpha_{1}})= z \textrm{ and } t_{z,w}(X^{\alpha_{2}}) = w .\] For each $t_{z,w} \in T_{Q}$, there are 3 elements $t \in T$ with $t|_{Q} = t_{z,w}$, determined by \[ t(X^{\omega_{1}})^{3} = z^{2}w \quad \textrm{ and } \quad t(X^{\omega_{2}}) = t(X^{-\omega_{1}})\cdot zw.\] The dimension of the simple modules with central character $t$ and the submodule structure of $M(t)$ depends only on $t|_{Q}$. Thus we begin by examining the $W_{0}$-orbits in $T_{Q}$. For a generic weight $t$, $P(t)$ and $Z(t)$ are empty so that $M(t)$ is irreducible by Theorem \ref{eq:Kato}. Thus we examine only the non-generic orbits.

\begin{proposition} \label{A2chars}
If $t \in T_{Q}$, and $P(t) \cup Z(t) \neq \emptyset$, then $t$ is in the $W_{0}$-orbit of one of the following weights:

\[ t_{1,1}, t_{1,q^{2}}, t_{q^{2},1}, t_{q^{2},q^{2}}, \{ t_{1,z} \, | \, z\in \mathbb{C}^{\times} z \neq 1, q^{\pm 2}\},\] \[ \textrm{ or } \{ t_{q^{2},z} \, | \, z\in \mathbb{C}^{\times} z \neq 1, q^{\pm 2}, q^{-4} \}\] \end{proposition}

\noindent \textit{Proof.} The proof consists of exhausting all possible cases for which configurations of roots appear in $Z(t)$ and $P(t)$, up to the action of $W_{0}$.

\noindent
\underline{Case 1:} If $Z(t)$ contains two positive roots, then it must contain the third. This implies $t = t_{1,1}$.

\noindent
\underline{Case 2:} If $Z(t)$ contains only one root, by applying an element of $W_{0}$, assume that it is $\alpha_{1}$. Then $t(X^{\alpha_{2}}) = t(X^{\alpha_{1} + \alpha_{2}})$, so either $P(t) = \emptyset$ or $P(t) = \{\alpha_{2}, \alpha_{1} + \alpha_{2}\}$. The first central character is $t_{1,z}$ for some $z \neq 1$ or $q^{\pm 2}$. (If $z =1$ or $z = q^{\pm 2}$, either $P(t)$ or $Z(t)$ would be larger.) For the second case, there are two potential choices for the orbit, arising from choosing $t(X^{\alpha_{2}}) = q^{2}$ or $q^{-2}$. However, $t_{1,q^{-2}}$ is in the same orbit as $t_{q^{2},1}$.

\noindent
\underline{Case 3:} Now assume that $Z(t) = \emptyset$. If $P(t)$ is not empty, assume that $\alpha_{1} \in P(t)$ and $t(X^{\alpha_{1}}) = q^{2}$. Then $t(X^{\alpha_{2}}) \neq q^{-2}$ by assumption on $Z(t)$. Then it is possible that $\alpha_{2} \in P(t)$, in which case $t = t_{q^{2},q^{2}}$. If $\alpha_{1} + \alpha_{2} \in P(t)$, then $t(X^{\alpha_{2}}) = q^{-4}$ and $t = t_{q^{2},q^{-4}} = s_{2}s_{1}t_{q^{2},q^{2}}$. Otherwise, $t = t_{q^{2},z}$ for some $z \neq 1, q^{\pm 2}, q^{-4}$. \\ $\square$

\noindent \textbf{Remark.} There is some redundancy present in the list of central characters above for specific values of $q$. If $q^{2} = -1$, then $t_{1,q^{2}}$, $t_{q^{2},1}$, and $t_{q^{2},q^{2}}$ are all in the same $W_{0}$-orbit. Also in this case, $t_{q^{2},z} = t_{-1,z} = s_{1}t_{-1,-z}$. If $q^{2} = 1$, then $t_{1,1} = t_{q^{2},1} = t_{1,q^{2}} = t_{q^{2},q^{2}}$, and $t_{1,z} = t_{q^{2},z}$. Also note that for every generic weight $t_{z,w}$, there are six weights in its $W_{0}$-orbits, all of which are generic.

It is helpful to draw a picture of the weights $\{t_{q^{x},q^{y}} | x,y \in \mathbb{R} \}$ for various values of $q$. Solid lines in these pictures show sets of the form \[ H_{\alpha} = \{ t \in T_{Q} \, | \, t(X^{\alpha}) = 1 \}, \] for $\alpha \in R^{+}$, while dashed lines denote sets of the form \[ H_{\alpha \pm \delta} = \{ t \in T_{Q} \, | \, t(X^{\alpha}) = q^{\pm 2} \}, \] for $\alpha \in R^{+}$. The weight $t_{q^{x},q^{y}}$ is the point $x$ units away from $H_{\alpha_{1}}$ and $y$ units away from $H_{\alpha_{2}}$. The action of $W_{0}$ is visible in this picture, as $s_{i}$ is given by reflection in the hyperplane $H_{\alpha_{i}}$.

\vspace{.3in}
\[ \beginpicture
\setcoordinatesystem units <1.5cm,1.5cm>
\setplotarea x from -4 to 4, y from -2.5 to 2.5
\linethickness=5pt
\setlinear
\plot -3.5 -1.2 -3.5 2.5 /
\plot -1.5 1.16 -5 -0.87 /
\plot -5 0.87 -1.5 -1.16 /
\setdashes
\plot -2.7 -1.2 -2.7 2.5 /
\plot -1.5 2.06 -5 0.04 /
\plot -4.93 1.74 -1.5 -0.29 /
\put{$\bullet$} at -3.5 .9
\put{$\bullet$} at -3.5 2
\put{$\bullet$} at -3.5 0
\put{$\bullet$} at -2.7 0.45
\put{$\bullet$} at -2.7 1.35
\put{$\bullet$} at -2.7 1.8
\put{$\bullet$} at -3.2 0.4
\put{$t_{1,1}$} at -3.9 0
\put{$t_{1,z}$} at -3.8 2
\put{$t_{1,q^{2}}$} at -3.9 .9
\put{$t_{q^{2},q^{2}}$} at -2.3 1.3
\put{$t_{q^{2},z}$} at -2.3 1.9
\put{$t_{q^{2},1}$} at -2.3 0.45
\put{$t_{z,w}$} at -3.2 .58
\put{$H_{\alpha_{1}}$} at -3.5 2.7
\put{$H_{\alpha_{1} + \delta}$} at -2.5 2.7
\put{$H_{\alpha_{2}}$} at -1.3 1.32
\put{$H_{\alpha_{2} + \delta}$} at -1.1 1.92
\put{$H_{\alpha_{1}+\alpha_{2}}$} at -5.3 1.03
\put{$H_{\alpha_{1}+\alpha_{2}+\delta}$} at -5.2 1.9
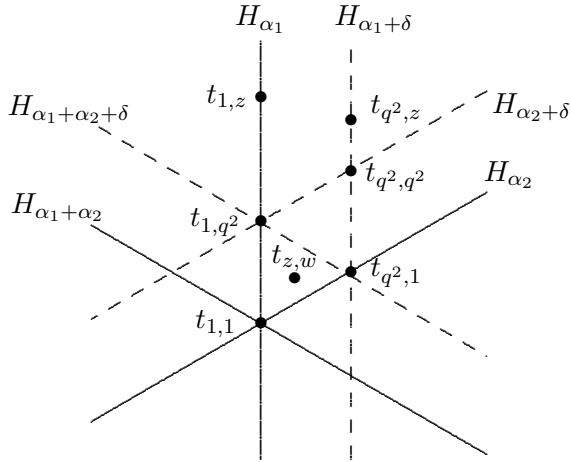
\put{Central characters with general $q$.} at -3.5 -1.5
\setlinear
\setsolid
\plot 2 -1.2 2 2.5 /
\plot 4 1.16 0.5 -0.87 /
\plot 0.5 0.87 4 -1.16 /
\plot 2.465 2.49  4 1.58 /
\setdashes
\plot 2.8 -1.2 2.8 2.5 /
\plot 3.6 -1.2 3.6 2.5 /
\plot 4 2.06 0.5 0.04 /
\plot 0.5 0.94 3.268 2.56 /
\plot 0.5 1.77 4 -0.3 /
\plot 1.29 2.24 4 0.64 /
\put{$\bullet$} at 2 .9
\put{$\bullet$} at 2 1.5
\put{$\bullet$} at 2 0
\put{$\bullet$} at 2.8 0.45
\put{$\bullet$} at 2.8 1.35
\put{$\bullet$} at 2.8 1.8
\put{$\bullet$} at 2.3 0.4
\put{$t_{1,1}$} at 1.6 0
\put{$t_{1,z}$} at 1.7 1.5
\put{$t_{1,q^{2}}$} at 1.6 .9
\put{$t_{q^{2},q^{2}}$} at 3.3 1.35
\put{$t_{q^{2},z}$} at 3.1 1.88
\put{$t_{q^{2},1}$} at 3.2 0.45
\put{$t_{z,w}$} at 2.3 .58
\put{$H_{\alpha_{1}}$} at 2.2 2.7
\put{$H_{\alpha_{1}+\delta}$} at 3.1 2.7
\put{$H_{\alpha_{1}-\delta}$} at 4 2.7
\put{$H_{\alpha_{2}}$} at 0.3 -0.44
\put{$H_{\alpha_{2}+\delta}$} at 0.3 0.46
\put{$H_{\alpha_{2}-\delta}$} at 0.3 1.26
\put{$H_{\alpha_{1} + \alpha_{2}}$} at 4.3 -0.87
\put{$H_{\alpha_{1} + \alpha_{2}+\delta}$} at 4.3 -0.02
\put{$H_{\alpha_{1} + \alpha_{2}-\delta}$} at 4.3 0.83
\put{$H_{\alpha_{1} + \alpha_{2}}$} at 4.3 1.68
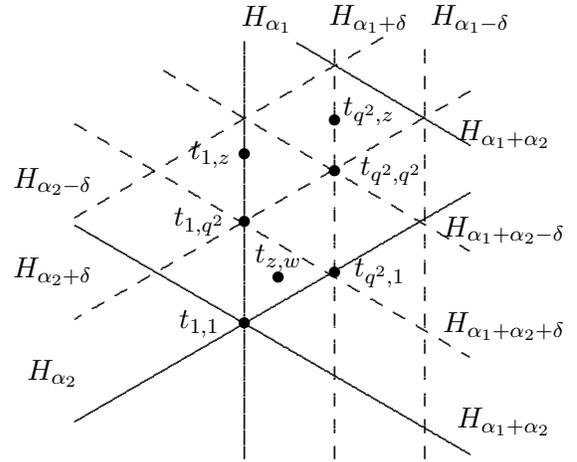
\put{Characters with $q^{2}$ a third root of unity.} at 2.5 -1.5
\endpicture\]

\vspace{.3in}
\[ \beginpicture
\setcoordinatesystem units <1.5cm,1.5cm>
\setplotarea x from -4 to 4, y from -2.5 to 2.5
\linethickness=5pt
\setlinear
\plot -5 -1 -5 2.5 /
\plot -3 1.16 -6 -0.58 /
\plot -6 0.58 -3 -1.16 /
\plot -3.045 -1 -3.045 2.5 /
\plot -5.37 2.49 -2.772 0.9675 /
\plot -6 1.66 -4.622 2.49 /
\setdashes
\plot -4 -1 -4 2.5 /
\plot -3 2.285 -6 .555 /
\plot -6 1.73 -2.772 -0.16 /
\put{$\bullet$} at -5 1.1275
\put{$\bullet$} at -5 1.8
\put{$\bullet$} at -5 0
\put{$\bullet$} at -4 2.25
\put{$\bullet$} at -4.7 0.4
\put{$t_{1,1}$} at -5.4 0
\put{$t_{1,z}$} at -5.3 1.8
\put{$t_{1,q^{2}}$} at -5.4 1.1
\put{$t_{q^{2},z}$} at -3.6 2.25
\put{$t_{z,w}$} at -4.7 .58
\put{$H_{\alpha_{1}}$} at -4.8 2.7
\put{$H_{\alpha_{1} \pm \delta}$} at -4 2.7
\put{$H_{\alpha_{1}}$} at -3.2 2.7
\put{$H_{\alpha_{2}}$} at -2.6 1.32
\put{$H_{\alpha_{2} \pm \delta}$} at -2.6 2.12
\put{$H_{\alpha_{1}+\alpha_{2}}$} at -6.3 0.74
\put{$H_{\alpha_{1}+\alpha_{2}\pm\delta}$} at -2.4 0.22
\put{$H_{\alpha_{1}+\alpha_{2}}$} at -5.7 2.7
\put{$H_{\alpha_{2}}$} at -6.3 1.66
\put{Central characters with $q^{2} = -1$.} at -4.5 -1.6
\setlinear
\setsolid
\plot 2 -1.3 2 2 /
\plot 3.2 0.812 1 -0.58 /
\plot 1 0.58 3.2 -0.812 /
\plot 2.9 -1.3 2.9 2 /
\plot 3.2 1.937 1 0.555 /
\plot 1 1.7238 3.2 0.434 /
\put{$\bullet$} at 2 1.8
\put{$\bullet$} at 2 0
\put{$\bullet$} at 2.3 0.4
\put{$t_{1,1}$} at 1.6 0
\put{$t_{1,z}$} at 1.7 1.8
\put{$t_{z,w}$} at 2.3 .58
\put{$H_{\alpha_{1}}$} at 2.2 2.2
\put{$H_{\alpha_{1}}$} at 3 2.2
\put{$H_{\alpha_{2}}$} at 3.5 0.74
\put{$H_{\alpha_{2}}$} at 3.5 1.74
\put{$H_{\alpha_{1}+\alpha_{2}}$} at 0.5 0.74
\put{$H_{\alpha_{1}+\alpha_{2}}$} at 0.6 1.82
\put{Central characters with $q^{2} = 1$.} at 1.8 -1.5
\endpicture\]


\vspace{.1in}

\noindent \textbf{Analysis of the characters.} \label{A2analysis} \begin{proposition} \label{eq:A2ones}
Fix $\epsilon = e^{2\pi i/3}$. The 1-dimensional representations of $\widetilde{H}$ are

\renewcommand{\baselinestretch}{1} \normalsize

$\begin{array}{ccc} \qquad \qquad L_{q^{2},\, q^{2}}: \widetilde{H} \rightarrow \mathbb{C} \qquad \,\,\,\,\,\,\,\, & L_{\epsilon q^{2},\, \epsilon^{2} q^{2}}: \widetilde{H} \rightarrow \mathbb{C} \qquad \,\,\,\,\, & L_{\epsilon^{2} q^{2},\, \epsilon q^{2}}: \widetilde{H} \rightarrow \mathbb{C} \qquad \\ T_{1} \mapsto q & T_{1} \mapsto q & T_{1} \mapsto q \\ T_{2} \mapsto q & T_{2} \mapsto q & T_{2} \mapsto q \\ X^{\omega_{1}} \mapsto q^{2} & X^{\omega_{1}} \mapsto \epsilon q^{2} & X^{\omega_{1}} \mapsto \epsilon^{2} q^{2} \\ X^{\omega_{2}} \mapsto q^{2} & X^{\omega_{2}} \mapsto \epsilon^{2} q^{2} & X^{\omega_{1}} \mapsto \epsilon q^{2} \end{array}$

$\begin{array}{ccc} \qquad \qquad L_{q^{-2},\, q^{-2}}: \widetilde{H} \rightarrow \mathbb{C} \qquad & L_{\epsilon q^{-2},\, \epsilon^{2} q^{-2}}: \widetilde{H} \rightarrow \mathbb{C} \qquad & L_{\epsilon^{2} q^{-2},\, \epsilon q^{-2}}: \widetilde{H} \rightarrow \mathbb{C} \qquad \\ T_{1} \mapsto -q^{-1} & T_{1} \mapsto -q^{-1} & T_{1} \mapsto -q^{-1} \\ T_{2} \mapsto -q^{-1} & T_{2} \mapsto -q^{-1} & T_{2} \mapsto -q^{-1} \\ X^{\omega_{1}} \mapsto q^{-2} & X^{\omega_{1}} \mapsto \epsilon q^{-2} & X^{\omega_{1}} \mapsto \epsilon^{2} q^{-2} \\ X^{\omega_{2}} \mapsto q^{-2} & X^{\omega_{2}} \mapsto \epsilon^{2} q^{-2} & X^{\omega_{1}} \mapsto \epsilon q^{-2} \end{array}$

\renewcommand{\baselinestretch}{1} \normalsize

\end{proposition}

\noindent \textit{Proof.} A straightforward check shows that the maps above respect the defining relations for $\widetilde{H}$, so that the maps are homomorphisms.

Let $\mathbb{C}v$ be any 1-dimensional $\widetilde{H}$-module. By equations $\eqref{eq:quadrel}$ and \eqref{eq:heckerel}, \[T_{i}v = qv \quad \textrm{or} \quad T_{i}v = -q^{-1}v \quad \textrm{ and } \quad T_{1}v=T_{2}v.\]

Once the actions of the $T_{i}$ are known, the relation $\eqref{eq:relation}$ forces the actions of the $X^{\alpha_{i}}$. This then defines the action of the $X^{\omega_{i}}$ up to a third root of unity. \\ $\square$

\noindent \textbf{Principal Series Modules}

Let $t \in T$. The principal series module is \[ M(t) = \textrm{Ind}_{\mathbb{C}[X]}^{\widetilde{H}} \mathbb{C}_{t} = \widetilde{H} \otimes_{\mathbb{C}[X]} \mathbb{C}_{t},\] where $\mathbb{C}_{t}$ is the one-dimensional $\mathbb{C}[X]$-module given by
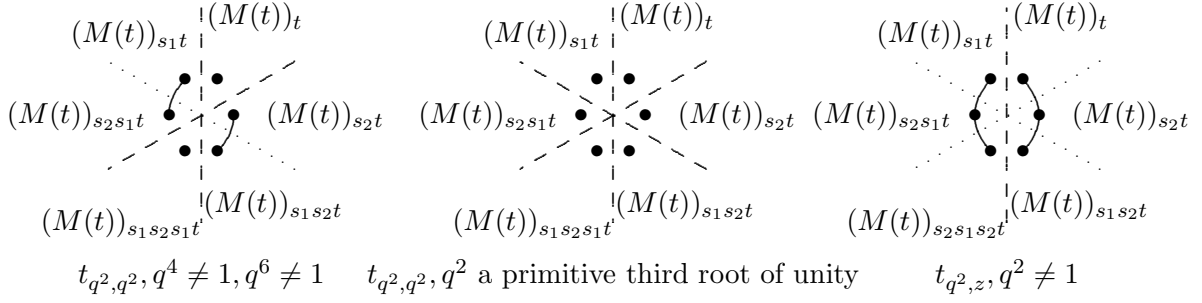
\[ \mathbb{C}_{t} = \textrm{span}\{v_{t}\} \quad \textrm{ and } \quad X^{\lambda}v_{t} = t(X^{\lambda})v_{t}.\]
By \eqref{PrinQuot}c, every irreducible $\widetilde{H}$ module is a quotient of some principal series module $M(t)$. Thus,  finding all the composition factors of $M(t)$ for all central characters $t$ will find all the irreducible $\widetilde{H}$-modules.

A weight $t|_{Q}$ corresponds to a point in the root lattice $Q$ as described above. The composition structure of $M(t)$ is largely determined by the structure of the operators $\tau_{i}$, which can be encoded in the picture of a small neighborhood of the point $t$ in $Q$. The following pictures show a small neighborhood of $t$, including the hyperplanes $H_{\alpha_{i}}$ that include the point $t$. Each picture shows 6 chambers around $t$, (bounded by the affine hyperplanes through $t$), representing the 6 (not necessarily distinct) elements of $W_{0}t$. The solid lines are hyperplanes $H_{\alpha_{i}}$ with $t(X^{\alpha_{i}}) = 1$, so that $s_{\alpha_{i}} \in W_{t}$. The chambers between any two successive solid lines are in bijection with the distinct elements of the orbit $W_{0}t$, which are in bijection with the cosets $W_{0}/W_{t}$ and thus the weight spaces of $M(t)$. These weight spaces are $|W_{t}|$-dimensional, so each chamber has $|W_{t}|$ points drawn in it.

The dotted lines represent hyperplanes $H_{\alpha_{i}}$ such that $t(X^{\alpha_{i}}) = q^{\pm 2}$. The operators $\tau_{i}:M_{t} \rightarrow M_{s_{i}t}$ and $\tau_{i}: M_{s_{i}t} \rightarrow M_{t}$ exist and are invertible exactly when $t(X^{\alpha_{i}}) \neq 1$ or $q^{2}$, i.e. exactly when the corresponding chambers in our picture do not have a hyperplane between them. Thus the hyperplanes in the picture of the neighborhood of $t$ show which $\tau$ operators are invertible and which are not. In most cases, this is enough to determine the exact composition factors of $M(t)$. Assume for now that $q^{2} \neq 1$.

\underline{Case 1:} $P(t)$ empty. By $\ref{eq:Kato}$, if $P(t) = \{ \alpha \in R^{+} | t(X^{\alpha_{1}}) = 1 \}$ is empty, then $M(t)$ is irreducible and is the only irreducible module with central character $t$. This case includes the central characters $t_{1,1}$, $t_{1,z}$, and $t_{z,w}$ for generic $z,w$ - that is, any $z$ and $w$ for which $P(t_{z,w}) = Z(t_{z,w}) = \emptyset$.

These pictures also show the weight space structure of $M(t)$. For $t = t_{1,1}$, all of $M(t)$ is in the $t$ weight space, so all the dots lie in the same chamber. For $t = t_{1,z}$, $s_{1}t = t$, so that the weights of $M(t)$ are $t$, $s_{2}t$, and $s_{1}s_{2}t$, and each weight space is 2-dimensional. The weight $ t = t_{z,w}$ is regular, so there are six different weights in $M(t_{z,w})$.

\[ \beginpicture
\setcoordinatesystem units <0.25cm,0.25cm>
\setplotarea x from -25 to 25, y from -5 to 6
\linethickness=3pt
\plot -23 6 -23 -6 /
\plot -17.804 3 -28.196 -3 /
\plot -17.804 -3 -28.196 3 /
\put{$\bullet$} at -21.2 4
\put{$\bullet$} at -21.65 3
\put{$\bullet$} at -22.1 2
\put{$\bullet$} at -20.75 5
\put{$\bullet$} at -20.3 6
\put{$\bullet$} at -22.55 1
\put{$(M(t))_{t}^{\textrm{gen}} = M(t)$} at -17 7.5
\put{$t_{1,1}, \, q^{2} \neq 1$} at -23 -9
\plot -4 6 -4 -6 /
\put{$\bullet$} at -2.6 0
\put{$\bullet$} at -0.1 0
\put{$\bullet$} at -3.1 2
\put{$\bullet$} at -1.75 5
\put{$\bullet$} at -3.1 -2
\put{$\bullet$} at -1.75 -5
\put{$(M(t))_{t}^{\textrm{gen}}$} at -0.1 6.5
\put{$(M(t))_{s_{2}t}^{\textrm{gen}}$} at 3.9 0
\put{$(M(t))_{s_{1}s_{2}t}^{\textrm{gen}}$} at 0.3 -6.5
\put{$t_{1,z}, \, q^{2} \neq 1$} at -2 -9
\put{$\bullet$} at 17.7 0
\put{$\bullet$} at 22.3 0
\put{$\bullet$} at 19.1 2
\put{$\bullet$} at 20.9 2
\put{$\bullet$} at 19.1 -2.1
\put{$\bullet$} at 20.9 -2.1
\put{$(M(t))_{t}$} at 23 5.5
\put{$(M(t))_{s_{2}t}$} at 26.9 0
\put{$(M(t))_{s_{1}s_{2}t}$} at 24 -5
\put{$(M(t))_{s_{1}t}$} at 16 4.5
\put{$(M(t))_{s_{2}s_{1}t}$} at 13 0
\put{$(M(t))_{s_{1}s_{2}s_{1}t}$} at 15.5 -6
\put{$t_{z,w}$} at 20 -9
\setquadratic
\plot -22.55 1 -22.55 1.5 -22.1 2 /
\plot -22.1 2  -22.1 2.5 -21.65 3 /
\plot -21.65 3 -21.65 3.5 -21.2 4 /
\plot -21.2 4 -21.2 4.5 -20.75 5 /
\plot -20.75 5 -20.75 5.5 -20.3 6 /
\plot -22.55 1 -22.1 1.5 -22.1 2 /
\plot -22.1 2  -21.65 2.5 -21.65 3 /
\plot -21.65 3 -21.2 3.5 -21.2 4 /
\plot -21.2 4 -20.75 4.5 -20.75 5 /
\plot -20.75 5 -20.3 5.5 -20.3 6 /
\plot -2.6 0 -1.35 0.8 -0.1 0 /
\plot -2.6 0 -1.35 -0.8 -0.1 0 /
\plot -3.1 2 -2.6 4 -1.75 5 /
\plot -3.1 2 -1.7 3.8 -1.75 5 /
\plot -3.1 -2 -2.6 -4 -1.75 -5 /
\plot -3.1 -2 -1.7 -3.8 -1.75 -5 /
\plot -3.1 2 -2.6 0 -3.1 -2 /
\plot -1.75 5 -0.1 0 -1.75 -5 /
\setquadratic
\circulararc 360 degrees from 22.3 0 center at 20 0
\setlinear
\setdots
\plot 20 6 20 -6 /
\plot 25.196 3 14.804 -3 /
\plot 25.196 -3 14.804 3 /
\plot 1.196 3 -9.196 -3 /
\plot 1.196 -3 -9.196 3 /
\endpicture \]

\underline{Case 2:} $Z(t) = \emptyset$, $P(t) \neq \emptyset$. This case includes the central characters $t_{q^{2},q^{2}}$, and $t_{q^{2},z}$. If $Z(t)$ is empty, then $M(t)$ is calibrated and the irreducible modules with central character $t$ are in one-to-one correspondence with the components of the calibration graph. These components are the same as the sets of chambers between two successive dashed hyperplanes in the picture of $t$. The segments in the picture represent the fact the the $\tau$ operators between two successive weight vectors are invertible, and thus these weight vectors must lie in the same composition factor of $M(t)$.

\[ \beginpicture
\setcoordinatesystem units <0.238cm,0.238cm>
\setplotarea x from -25 to 25, y from -5 to 6
\linethickness=3pt
\put{$\bullet$} at -21.2 0
\put{$\bullet$} at -24.8 0
\put{$\bullet$} at -22.1 2
\put{$\bullet$} at -23.9 2
\put{$\bullet$} at -22.1 -2
\put{$\bullet$} at -23.9 -2
\put{$(M(t))_{t}$} at -20 5.5
\put{$(M(t))_{s_{2}t}$} at -16.1 0
\put{$(M(t))_{s_{1}s_{2}t}$} at -19 -5
\put{$(M(t))_{s_{1}t}$} at -27 4.5
\put{$(M(t))_{s_{2}s_{1}t}$} at -30 0
\put{$(M(t))_{s_{1}s_{2}s_{1}t}$} at -27.5 -6
\put{$t_{q^{2},q^{2}}, q^{4} \neq 1, q^{6} \neq 1$} at -23 -9
\put{$\bullet$} at 1.8 0
\put{$\bullet$} at -1.8 0
\put{$\bullet$} at 0.9 2
\put{$\bullet$} at -0.9 2
\put{$\bullet$} at 0.9 -2
\put{$\bullet$} at -0.9 -2
\put{$(M(t))_{t}$} at 3 5.5
\put{$(M(t))_{s_{2}t}$} at 6.9 0
\put{$(M(t))_{s_{1}s_{2}t}$} at 4 -5
\put{$(M(t))_{s_{1}t}$} at -4 4.5
\put{$(M(t))_{s_{2}s_{1}t}$} at -7 0
\put{$(M(t))_{s_{1}s_{2}s_{1}t}$} at -4.5 -6
\put{$t_{q^{2},q^{2}}, q^{2}$ a primitive third root of unity} at 0 -9
\put{$\bullet$} at 23.8 0
\put{$\bullet$} at 20.2 0
\put{$\bullet$} at 22.9 2
\put{$\bullet$} at 21.1 2
\put{$\bullet$} at 22.9 -2
\put{$\bullet$} at 21.1 -2
\put{$(M(t))_{t}$} at 25 5.5
\put{$(M(t))_{s_{2}t}$} at 28.9 0
\put{$(M(t))_{s_{1}s_{2}t}$} at 26 -5
\put{$(M(t))_{s_{1}t}$} at 18 4.5
\put{$(M(t))_{s_{2}s_{1}t}$} at 15 0
\put{$(M(t))_{s_{2}s_{1}s_{2}t}$} at 17.5 -6
\put{$t_{q^{2},z}, q^{2} \neq 1$} at 22 -9
\setquadratic
\plot -24.8 0 -24.6 1 -23.9 2 /
\plot -21.2 0 -21.4 -1 -22.1 -2 /
\plot 22.9 2 23.8 0 22.9 -2 /
\plot 21.1 2 20.2 0 21.1 -2 /
\setlinear
\setdashes
\plot -17.804 3 -28.196 -3 /
\plot -23 6 -23 -6 /
\plot 22 -6 22 6 /
\plot 0 6 0 -6 /
\plot 5.196 3 -5.196 -3 /
\plot 5.196 -3 -5.196 3 /
\setdots
\plot -17.804 -3 -28.196 3 /
\plot 27.196 3 16.804 -3 /
\plot 27.196 -3 16.804 3 /
\endpicture \]

\underline{Case 3:} $P(t) \neq \emptyset,$ $Z(t) \neq \emptyset$. The only central characters with both $Z(t)$ and $P(t)$ non-empty are $t_{1,1}$ and $t_{1,z}$ when $q^{2} = 1$, and $t_{q^{2},1}$ and $t_{1,q^{2}}$ in all cases. If $q^{4} = 1$, then $t_{q^{2},1}$ and $t_{1,q^{2}}$ are in the same orbit, and are in the same orbit as $t_{q^{2},q^{2}}$. If $q^{2} = 1$, then $t_{q^{2},1} = t_{1,q^{2}} = t_{1,1}$.

If $q^{2} \neq \pm 1$ and $t|_{Q} = t_{1,q^{2}}$ or $t_{q^{2},1}$, then Lemma $\ref{eq:lemma}$ shows that $M(t)$ has two 3-dimensional composition factors. When $q^{2} = -1$ and $t|_{Q} = t_{1,q^{2}}$, Theorem $\ref{2dims}$ shows the two dimensional weight space $M(t)_{t}^{\textrm{gen}}$ makes up an entire composition factor, as does $M(t)_{s_{1}s_{2}t}^{\textrm{gen}}$. The remaining composition factors are two copies of the 1-dimensional module with weight $s_{2}t$.

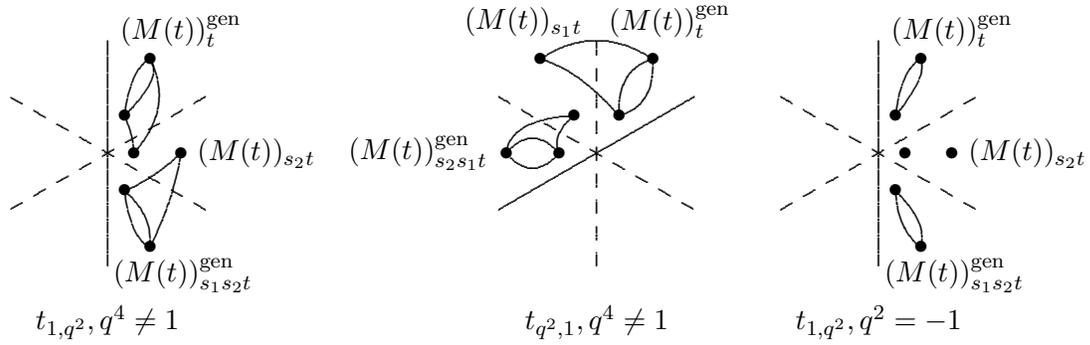
\[ \beginpicture
\setcoordinatesystem units <0.25cm,0.25cm>
\setplotarea x from -14 to 40, y from -5 to 6
\linethickness=3pt
\plot -12 6 -12 -6 /
\put{$\bullet$} at -10.6 0
\put{$\bullet$} at -8.1 0
\put{$\bullet$} at -11.1 2
\put{$\bullet$} at -9.75 5
\put{$\bullet$} at -11.1 -2
\put{$\bullet$} at -9.75 -5
\put{$(M(t))_{t}^{\textrm{gen}}$} at -8.1 6.5
\put{$(M(t))_{s_{2}t}$} at -4.1 0
\put{$(M(t))_{s_{1}s_{2}t}^{\textrm{gen}}$} at -8 -6.5
\put{$t_{1,q^{2}}, q^{4} \neq 1$} at -12 -9
\plot 19.196 3 8.804 -3 /
\put{$\bullet$} at 12 0
\put{$\bullet$} at 9.2 0
\put{$\bullet$} at 15.2 2
\put{$\bullet$} at 17 5
\put{$\bullet$} at 12.8 2
\put{$\bullet$} at 11 5
\put{$(M(t))_{t}^{\textrm{gen}}$} at 17.9 7
\put{$(M(t))_{s_{1}t}$} at 10.1 7
\put{$(M(t))_{s_{2}s_{1}t}^{\textrm{gen}}$} at 4.5 0
\put{$t_{q^{2},1}, q^{4} \neq 1$} at 14 -9
\setquadratic
\plot -11.1 2 -10.6 4 -9.75 5 /
\plot -11.1 2 -9.7 3.8 -9.75 5 /
\plot -11.1 -2 -10.6 -4 -9.75 -5 /
\plot -11.1 -2 -10 -3.4 -9.75 -5 /
\plot -11.1 2 -10.75 1 -10.6 0 /
\plot -9.75 5 -9.3 2.5 -10.6 0 /
\plot -9.75 -5 -8.7 -2.5 -8.1 0 /
\plot -11.1 -2 -9.6 -1.2 -8.1 0 /
\plot 9.2 0 10.6 0.8 12 0 /
\plot 9.2 0 10.6 -0.8 12 0 /
\plot 15.2 2 16.6 3 17 5 /
\plot 15.2 2 15.6 4 17 5 /
\plot 15.2 2 13.05 4 11 5 /
\plot 17 5 14 6 11 5 /
\plot 12 0 12 1 12.8 2 /
\plot 9.2 0 10.2 1.4 12.8 2 /
\setlinear
\setdashes
\plot -6.804 3 -17.196 -3 /
\plot -6.804 -3 -17.196 3 /
\plot 8.804 3 19.196 -3 /
\plot 14 -6 14 6 /
\setsolid
\plot 29 6 29 -6 /
\put{$\bullet$} at 30.4 0
\put{$\bullet$} at 32.9 0
\put{$\bullet$} at 29.9 2
\put{$\bullet$} at 31.25 5
\put{$\bullet$} at 29.9 -2
\put{$\bullet$} at 31.25 -5
\put{$(M(t))_{t}^{\textrm{gen}}$} at 32.9 6.5
\put{$(M(t))_{s_{2}t}$} at 36.9 0
\put{$(M(t))_{s_{1}s_{2}t}^{\textrm{gen}}$} at 33 -6.5
\put{$t_{1,q^{2}}, q^{2} = -1$} at 29 -9
\setquadratic
\plot 29.9 2 30.4 4 31.25 5 /
\plot 29.9 2 31.1 3.8 31.25 5 /
\plot 29.9 -2 31.1 -3.6 31.25 -5 /
\plot 29.9 -2 30.3 -3.8 31.25 -5 /
\setlinear
\setdashes
\plot 23.804 -3 34.196 3 /
\plot 23.804 3 34.196 -3 /
\endpicture \]

Explicitly, let $\mathbb{C}_{q^{2},1}$ be the 1-dimensional $\widetilde{H}_{\{1\}}$-module spanned by $v_{t}$ and let $\mathbb{C}_{1,q^{-2}}$ be the 1-dimensional $\widetilde{H}_{\{2\}}$-module spanned by $v_{s_{2}s_{1}t}$, given by \[ X^{\lambda}v_{t}  = t(X^{\lambda})v_{t}, \quad \textrm{ and } \quad T_{1}v_{t} = qv_{t},  \quad \textrm{ and }\] \[ X^{\lambda}v_{s_{2}s_{1}t} = (s_{2}s_{1}t)(X^{\lambda})v_{s_{2}s_{1}t}  \quad  \textrm{ and } \quad T_{1}v_{s_{2}s_{1}t} = -q^{-1}v_{s_{2}s_{1}t}.\] Then \[M = \widetilde{H} \otimes_{\widetilde{H}_{\{1\}}} \mathbb{C}_{q^{2},1} \quad \textrm{ and } \quad N = \widetilde{H} \otimes_{\widetilde{H}_{\{2\}}} \mathbb{C}_{1,q^{-2}}\] are 3-dimensional $\widetilde{H}$-modules with central character $t_{q^{2},1}$.

\begin{proposition} \label{A2fourthstandard} Let $M = \widetilde{H} \otimes_{\widetilde{H}_{\{1\}}} \mathbb{C}_{q^{2},1} \quad \textrm{ and } \quad N = \widetilde{H} \otimes_{\widetilde{H}_{\{2\}}} \mathbb{C}_{1,q^{-2}}$.

(a) If $q^{4} \neq 1$ then $M$ and $N$ are irreducible.

(b) If $q^{2} = -1$ then $M_{s_{1}t}$ is an irreducible submodule of $M$ and $N_{s_{1}t}$ is an irreducible submodule of $N$. The quotients $N/N_{s_{1}t}$ and $M/M_{s_{1}t}$ are irreducible.

\end{proposition}

\noindent \textit{Proof.} (a) Assume $q^{4} \neq 1$. If either $M$ or $N$ were reducible, it would have a 1-dimensional submodule or quotient, which cannot happen since the 1-dimensional modules have central character $t_{q^{2},q^{2}}$. Thus both $M$ and $N$ are reducible.

(b) If $t|_{Q} = t_{q^{2},1}$, then the action of $\tau_{1}$ is non-zero on $M_{t}^{\textrm{gen}}$ by theorem $\ref{weightbasis}$, and $M_{t}^{\textrm{gen}}$ is not a submodule of $M$. But $M$ is not irreducible, and the only possible remaining submodule is $M_{s_{1}t}$. A similar argument shows the result for $N$ as well. $\square$

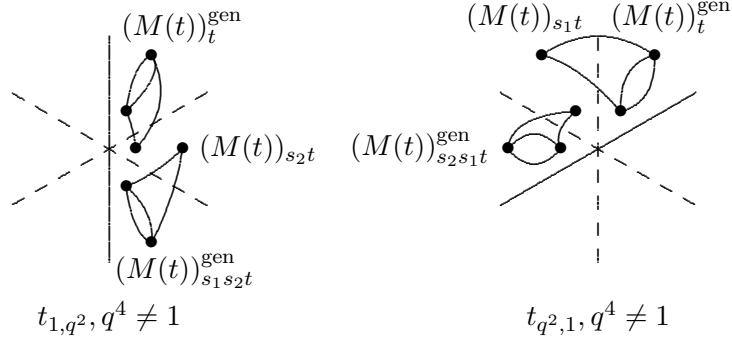
\[ \beginpicture
\setcoordinatesystem units <0.25cm,0.25cm>
\setplotarea x from -25 to 25, y from -5 to 6
\linethickness=3pt
\plot -12 6 -12 -6 /
\put{$\bullet$} at -10.6 0
\put{$\bullet$} at -8.1 0
\put{$\bullet$} at -11.1 2
\put{$\bullet$} at -9.75 5
\put{$\bullet$} at -11.1 -2
\put{$\bullet$} at -9.75 -5
\put{$(M(t))_{t}^{\textrm{gen}}$} at -8.1 6.5
\put{$(M(t))_{s_{2}t}$} at -4.1 0
\put{$(M(t))_{s_{1}s_{2}t}^{\textrm{gen}}$} at -8 -6.5
\put{$t_{1,q^{2}}, q^{4} \neq 1$} at -12 -9
\plot 19.196 3 8.804 -3 /
\put{$\bullet$} at 12 0
\put{$\bullet$} at 9.2 0
\put{$\bullet$} at 15.2 2
\put{$\bullet$} at 17 5
\put{$\bullet$} at 12.8 2
\put{$\bullet$} at 11 5
\put{$(M(t))_{t}^{\textrm{gen}}$} at 17.9 7
\put{$(M(t))_{s_{1}t}$} at 10.1 7
\put{$(M(t))_{s_{2}s_{1}t}^{\textrm{gen}}$} at 4.5 0
\put{$t_{q^{2},1}, q^{4} \neq 1$} at 14 -9
\setquadratic
\plot -11.1 2 -10.6 4 -9.75 5 /
\plot -11.1 2 -9.7 3.8 -9.75 5 /
\plot -11.1 -2 -10.6 -4 -9.75 -5 /
\plot -11.1 -2 -10 -3.4 -9.75 -5 /
\plot -11.1 2 -10.75 1 -10.6 0 /
\plot -9.75 5 -9.3 2.5 -10.6 0 /
\plot -9.75 -5 -8.7 -2.5 -8.1 0 /
\plot -11.1 -2 -9.6 -1.2 -8.1 0 /
\plot 9.2 0 10.6 0.8 12 0 /
\plot 9.2 0 10.6 -0.8 12 0 /
\plot 15.2 2 16.6 3 17 5 /
\plot 15.2 2 15.6 4 17 5 /
\plot 15.2 2 13.05 4 11 5 /
\plot 17 5 14 6 11 5 /
\plot 12 0 12 1 12.8 2 /
\plot 9.2 0 10.2 1.4 12.8 2 /
\setlinear
\setdashes
\plot -6.804 3 -17.196 -3 /
\plot -6.804 -3 -17.196 3 /
\plot 8.804 3 19.196 -3 /
\plot 14 -6 14 6 /
\endpicture \]

The modules with central character $t$ such that $t|_{Q} = t_{1,q^{2}}$ can be constructed in an entirely analoguous fashion, for $q^{2} \neq 1$.

If $q^{2} = 1$, then the results of section \ref{clifford} suffice to classify the representations of $\widetilde{H}$ with central characters $t_{1,1}$ and $t_{1,z}$ for $z \neq q^{\pm 2}$. Specifically, if $t|_{Q} = t_{1,1}$, then $W_{t} = W_{0}$ and a $\widetilde{H}$-module is merely a $W_{0}$-module (via the isomorphism $H \cong \mathbb{C}[W_{0}]$) on which $X^{\lambda} \in \mathbb{C}[X]$ acts by the scalar $t(X^{\lambda})$. In fact, $M(t)$ considered as a $W_{0}$ module is the regular $W_{0}$ module. If $t|_{Q} = t_{1,z}$ for $z \neq 1$, then $W_{t} = \{ 1, s_{1} \}$. Since $W_{t}$ has two 1-dimensional irreducible representations, $\widetilde{H}$ has two irreducible 3-dimensional representations obtained by inducing up from $\widetilde{H}_{\{1\}}$.

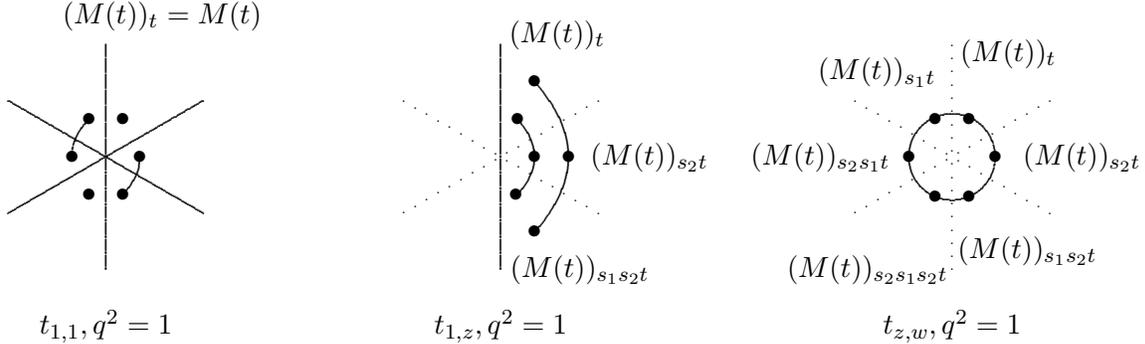
\[ \beginpicture
\setcoordinatesystem units <0.25cm,0.25cm>
\setplotarea x from -25 to 25, y from -5 to 6
\linethickness=3pt
\put{$\bullet$} at -21.2 0
\put{$\bullet$} at -24.8 0
\put{$\bullet$} at -22.1 2
\put{$\bullet$} at -23.9 2
\put{$\bullet$} at -22.1 -2
\put{$\bullet$} at -23.9 -2
\put{$t_{1,1}, q^{2} = 1$} at -23 -9
\put{$\bullet$} at -0.2 0
\put{$\bullet$} at 1.6 0
\put{$\bullet$} at -1.1 2
\put{$\bullet$} at -0.2 4
\put{$\bullet$} at -1.2 -2
\put{$\bullet$} at -0.2 -4
\put{$t_{1,z}, q^{2} = 1$} at -2 -9
\put{$\bullet$} at 19.7 0
\put{$\bullet$} at 24.3 0
\put{$\bullet$} at 21.1 2
\put{$\bullet$} at 22.9 2
\put{$\bullet$} at 21.1 -2.1
\put{$\bullet$} at 22.9 -2.1
\put{$t_{z,w}, q^{2} = 1$} at 22 -9
\put{$(M(t))_{t} = M(t)$} at -20 7.5
\put{$(M(t))_{t}$} at 1 6.5
\put{$(M(t))_{s_{2}t}$} at 5.9 0
\put{$(M(t))_{s_{1}s_{2}t}$} at 2.2 -6
\put{$(M(t))_{t}$} at 25 5.5
\put{$(M(t))_{s_{2}t}$} at 28.9 0
\put{$(M(t))_{s_{1}s_{2}t}$} at 26 -5
\put{$(M(t))_{s_{1}t}$} at 18 4.5
\put{$(M(t))_{s_{2}s_{1}t}$} at 15 0
\put{$(M(t))_{s_{2}s_{1}s_{2}t}$} at 17.5 -6
\setquadratic
\plot -24.8 0 -24.6 1 -23.9 2 /
\plot -21.2 0 -21.4 -1 -22.1 -2 /
\plot -1.1 2 -0.2 0 -1.1 -2 /
\plot -0.2 4 1.6 0 -0.2 -4 /
\circulararc 360 degrees from 24.3 0 center at 22 0
\setlinear
\plot -17.804 3 -28.196 -3 /
\plot -23 6 -23 -6 /
\plot -2 -6 -2 6 /
\plot -17.804 -3 -28.196 3 /
\setdots
\plot 27.196 3 16.804 -3 /
\plot 27.196 -3 16.804 3 /
\plot 22 -6 22 6 /
\plot 3.196 3 -7.196 -3 /
\plot 3.196 -3 -7.196 3 /
\endpicture \]

The following tables summarize the classification. It should be noted that for any value of $q$ with $q^{2} \neq \pm 1$ and $q^{2}$ not a primitive third root of unity, the representation theory of $\widetilde{H}$ can be described in terms of $q$ only. If $q^{2}$ is a primitive root of unity of order 3 or less, then the representation theory of $\widetilde{H}$ does not fit that same description. This fact can be seen through a number of different lenses. It is a reflection of the fact that the sets $P(t)$ and $Z(t)$ for all possible central characters $t$ can be described solely in terms of $q$. In the pictures of the torus, this is reflected in the fact that the hyperplanes $H_{\alpha}$ and $H_{\alpha \pm \delta}$ are distinct \textit{unless} $q^{2}$ is a root of unity of order 3 or less. When these hyperplanes coincide, the sets $P(t)$ and $Z(t)$ change for characters on those hyperplanes.
\begin{table}[htpb]
\centering
$\begin{array}{| c | c | c | c | c |}
\hline \, & \multicolumn{4}{|c|}{\textrm{Dims. of Irreds.}} \\
\hline t & q^{6} \neq 1, q^{4} \neq 1 & q^{6} = 1 & q^{2} = -1 & q = -1 \\
\hline t_{1,1} & 6 & 6 & 6 & 1,1,2 \\
\hline t_{1,z} & 6 & 6 & 6 & 3,3 \\
\hline t_{1,q^{2}} & 3,3 & 3,3 & 1,2,2 & \textrm{N/A} \\
\hline t_{q^{2},1} & 3,3 & 3,3 & \textrm{N/A} & \textrm{N/A} \\
\hline t_{q^{2},q^{2}} & 1,1,2,2 & 1,1,1,1,1,1 & \textrm{N/A} & \textrm{N/A} \\
\hline t_{q^{2},z} & 3,3 & 3,3 & 3,3 & \textrm{N/A} \\
\hline t_{z,w} & 6 & 6 & 6 & 6 \\ \hline \end{array}$
\caption{Table of possible central characters in Type $A_{2}$, for varying values of $q$.}
\end{table}

\section{Type $C_{2}$}

The type $C_{2}$ root system is $R = \{ \pm \alpha_{1}, \pm \alpha_{2}, \pm (\alpha_{1} + \alpha_{2}), \pm(2\alpha_{1} + \alpha_{2}) \},$ where $\langle \alpha_{1}, \alpha_{2}^{\vee} \rangle = -1$ and $\langle \alpha_{2}, \alpha_{1}^{\vee} \rangle = -2$. Then the Weyl group is $W_{0} = \langle s_{1}, s_{2} \, | \, s_{1}^{2} = s_{2}^{2} = 1, s_{1}s_{2}s_{1}s_{2} = s_{2}s_{1}s_{2}s_{1} \rangle$, which is isomorphic to the dihedral group of order 8. The simple roots are $\alpha_{1}$ and $\alpha_{2}$, with additional positive roots $\alpha_{1} + \alpha_{2}$ and $2\alpha_{1} + \alpha_{2}$. Then the action of $W_{0}$ on $R$ can be seen in the following picture, where $s_{i}$ acts by reflection through $H_{\alpha_{i}}$,. the hyperplane perpendicular to $\alpha_{i}$.

\[ \beginpicture
\setcoordinatesystem units <0.3cm,0.3cm>
\setplotarea x from -25 to 25, y from -5 to 6
\linethickness=5pt
\setlinear
\thicklines
\arrow <10 pt> [.2,.67] from 0 0 to 0 6
\arrow <10 pt> [.2,.67] from 0 0 to 0 -6
\arrow <10 pt> [.2,.67] from 0 0 to 6 6
\arrow <10 pt> [.2,.67] from 0 0 to 6 -6
\arrow <10 pt> [.2,.67] from 0 0 to -6 6
\arrow <10 pt> [.2,.67] from 0 0 to -6 -6
\arrow <10 pt> [.2,.67] from 0 0 to 6 0
\arrow <10 pt> [.2,.67] from 0 0 to -6 0
\put{The type $C_{2}$ root system} at 0 -10
\put{$\alpha_{1}$} at 8 0
\put{$\alpha_{2}$} at -6 7
\put{$\alpha_{1} + \alpha_{2}$} at 0 7
\put{$2\alpha_{1} + \alpha_{2}$} at 8 7
\endpicture\]

The fundamental weights satisfy \[ \begin{matrix} \omega_{1} = \alpha_{1} + \frac{1}{2}\alpha_{2} & \quad & \alpha_{1} = 2\omega_{1} - \omega_{2} \\ \omega_{2} = \alpha_{1} + \alpha_{2} & \quad & \alpha_{2} = 2\omega_{2} - 2\omega_{1}. \end{matrix}\] Let \[P = \mathbb{Z}\textrm{-span}\{\omega_{1},\omega_{2}\} \textrm{ and } Q = \mathbb{Z}\textrm{-span}(R) \] be the weight lattice and root lattice of $R$, respectively.

Then the affine Hecke algebra $\widetilde{H}$ is defined as in section \ref{subsectHecke}. Let \[\mathbb{C}[Q] = \{ X^{\lambda} \, | \, \lambda \in Q\} \textrm{ and } T_{Q} = \textrm{Hom}_{\mathbb{C}\textrm{-alg}}(\mathbb{C}[Q],\mathbb{C}).\] Define
\[ t_{z,w}: \mathbb{C}[Q]: \rightarrow \mathbb{C} \quad \textrm{ by } t_{z,w}(X^{\alpha_{1}}) = z \textrm{ and } t_{z,w}(X^{\alpha_{2}}) = w .\]  Pictorially, a weight $t_{q^{x},q^{y}} \in T_{Q}$, for $x,y \in \mathbb{R}$, is identified with with the point $x$ units from the hyperplane $H_{\alpha_{1}}$ and $y$ units from the hyperplane $H_{\alpha_{2}}$. Then \[ H_{\alpha} = \{ t \in T_{Q} \, | \, t(X^{\alpha}) = 1\}, \] and we define \[ H_{\alpha \pm \delta} = \{ t \in T_{Q} \, | \, t(X^{\alpha}) = q^{\pm 2} \}.\]

\[ \beginpicture
\setcoordinatesystem units <0.3cm,0.3cm>
\setplotarea x from -25 to 25, y from -5 to 6
\linethickness=5pt
\setlinear
\thicklines
\plot 0 -8 0 8 /
\plot -8 -8 8 8 /
\plot -8 0 8 0 /
\plot -8 8 8 -8 /
\setdashes
\plot 2 -8 2 8 /
\plot-8 -6 6 8 /
\plot -8 2 8 2 /
\plot -6 8 8 -6 /
\put{The torus $T_{Q}$} at 0 -12
\put{$H_{\alpha_{1}}$} at -1.2 9
\put{$H_{\alpha_{2}}$} at 9 9
\put{$H_{\alpha_{1} + \alpha_{2}}$} at 10.5 0
\put{$H_{2\alpha_{1} + \alpha_{2}}$} at -8.8 9
\put{$H_{\alpha_{1} + \delta}$} at 2.5 -8.8
\put{$H_{\alpha_{2}+\delta}$} at -8.5 -4
\put{$H_{\alpha_{1} + \alpha_{2} + \delta}$} at -11 2
\put{$H_{2\alpha_{1} + \alpha_{2} +\delta}$} at 11 -6
\endpicture\]

Note that $W_{0}$ acts on the points in this picture by letting $s_{i}$ denote reflection through the hyperplane $H_{\alpha_{i}}$, for $i=1$ or $2$.

For all weights $t_{z,w} \in T_{Q}$, there are 2 elements $t \in T$ with $t|_{Q} = t_{z,w}$, determined by \[ t(X^{\omega_{1}})^{2} = z^{2}w \quad \textrm{ and } \quad t(X^{\omega_{2}}) = zw.\] We denote these two elements as $t_{z,w,1}$ and $t_{z,w,2}$. Which particular weight $t_{z,w,i}$ is which is unimportant since we will always be examining them together. And in fact, most of the time, we will only refer to the restricted weight $t_{z,w}$, since the dimension of the modules with central character $t$ depends only on $t|_{Q}$. One important remark, though, is that if $t(X^{\alpha_{1}}) = -1$, then the two weights $t$ with $t|_{Q} = t_{-1,w}$ are in the same $W_{0}$-orbit and represent the same central character. To see this, let $t|_{Q} = t_{-1,w}$. Then $t(X^{\omega_{1}}) = w^{1/2}$, and \[s_{1}t(X^{\omega_{1}}) = t(X^{\omega_{2}-\omega_{1}}) = -w/t(X^{\omega_{1}}) = - t(X^{\omega_{1}}),\] while \[ s_{1}t(X^{\omega_{2}}) = t(X^{\omega_{2}}).\]

The dimension of the modules with central character $t$ and the submodule structure of $M(t)$ depends only on $t|_{Q}$. Thus we begin by examining the $W_{0}$-orbits in $T_{Q}$. The structure of the modules with weight $t$ depends virtually exclusively on $P(t) = \{ \alpha \in R^{+} \, | \, t(X^{\alpha})=q^{\pm 2} \}$ and $Z(t) = \{ \alpha \in R^{+} \, | \, t(X^{\alpha}) = 1 \}$. For a generic weight $t$, $P(t)$ and $Z(t)$ are empty, so we examine only the non-generic orbits.

\begin{proposition} \label{B2chars}
Let $q$ be generic. If $t \in T_{Q}$, and $P(t) \cup Z(t) \neq \emptyset$, then $t$ is in the $W_{0}$-orbit of one of the following weights:

\[ t_{1,1}, t_{-1,1}, t_{1,q^{2}}, t_{q^{2},1}, t_{\pm q, 1}, t_{q^{2},q^{2}}, t_{-1,q^{2}}, \{ t_{1,z} \, | \, z \neq 1, q^{\pm 2} \}, \{t_{z,1} \, | \, z \neq \pm 1, q^{\pm 2}, \pm q^{\pm 1} \}, \] \[ \{ t_{q^{2},z} \, | \, z \neq 1, q^{\pm 2}, q^{-4}, q^{-6} \}, \textrm{ or } \{ t_{z,q^{2}} \, | \, z \neq \pm 1, q^{\pm 2}, -q^{-2}, q^{-4}, \pm q^{-1} \}.\]

\end{proposition}

\noindent \textit{Proof.} The proof consists of exhausting the possible cases for the roots that appear in the sets $Z(t)$ and $P(t)$, up to the action of $W_{0}$. These cases are summarized in pictures after the proof, with the points denoting representatives of the orbits in those cases. In the following, we refer to $\alpha_{1}$ and $\alpha_{1}+\alpha_{2}$ as ``short'' roots, and $\alpha_{2}$ and $2\alpha_{1}+\alpha_{2}$ as ``long'' roots.

\underline{Case 1: $|Z(t)| \geq 2$.} If $Z(t)$ contains a short root and any other root, then $t = t_{1,1}$. If $Z(t)$ contains two long roots, then $t = t_{-1,1}$.

\vspace{.05in}

\underline{Case 2: $|Z(t)| = 1$.} If $Z(t)$ contains exactly one root, we may assume it is either $\alpha_{1}$ or $\alpha_{2}$. If $t(X^{\alpha_{1}}) = 1$, then $t(X^{\alpha_{2}}) = t(X^{\alpha_{1}+\alpha_{2}}) = t(X^{2\alpha_{1}+\alpha_{2}})$. Thus either $P(t) = \{ \alpha_{2}, \alpha_{1}+\alpha_{2}, 2\alpha_{1}+\alpha_{2}\}$, or $P(t) = \emptyset$. That is, $t$ is in the orbit of $t_{1,q^{2}}$ or $t_{1,z}$ for some $z \neq q^{\pm 2}, 1$. If $t(X^{\alpha_{2}}) = 1$, then $t(X^{\alpha_{1}}) = t(X^{\alpha_{1}+\alpha_{2}})$. Then either $P(t) = \{ \alpha_{1}, \alpha_{1}+\alpha_{2}\}$, $P(t) = \{ 2\alpha_{1}+\alpha_{2}\}$, or $P(t) = \emptyset$. These are the orbits of $t_{q^{2},1}$, $t_{\pm q,1}$, and $t_{z,1}$, respectively, where $z \neq q^{\pm 2}, \pm q^{\pm 1}, \textrm{ or } 1$.

\vspace{.05in}

\underline{Case 3: $Z(t) = \emptyset$.} Now assume that $Z(t)$ is empty, and $P(t)$ is not empty. First assume that $P(t)$ contains at least one short root. We can apply an element of $w$ to assume that $\alpha_{1} \in P(t)$ and $t(X^{\alpha_{1}}) = q^{2}$. Then if $\alpha_{2} \in P(t)$, we must have $t(X^{\alpha_{2}}) = q^{2}$ or else $Z(t)$ would be non-empty. Thus $t = t_{q^{2},q^{2}}$. If $\alpha_{1} + \alpha_{2} \in P(t)$, then $t(X^{\alpha_{1}+\alpha_{2}}) = q^{\pm 2}$, so that either $t(X^{\alpha_{2}}) = 1$ or $t(X^{2\alpha_{1}+\alpha_{2}}) = 1$. If $2\alpha_{1} + \alpha_{2} \in P(t)$, then $t(X^{2\alpha_{1}+\alpha_{2}}) = q^{-2}$ or else $\alpha_{1}+\alpha_{2} \in Z(t)$. Hence $t(X^{\alpha_{2}}) = q^{-6}$. But then $s_{2}s_{1}s_{2}t = t_{q^{2},q^{2}}$. If $P(t) = \{ \alpha_{1}\}$, then $t = t_{q^{2},z}$ for some $z \neq 1, q^{\pm 2}, q^{-4}, q^{-6}$.

Now, assume that $P(t)$ contains a long root but no short roots. Then we may assume that $t(X^{\alpha_{2}}) = q^{2}$. If $t(X^{2\alpha_{1}+\alpha_{2}}) = q^{2}$ then $t(X^{\alpha_{1}}) = -1$. If $t(X^{2\alpha_{1} + \alpha_{2}}) = q^{-2}$ then $t(X^{\alpha_{1}}) = -q^{-2}$. However, $s_{1}s_{2}s_{1}t_{-q^{-2},q^{2}} = t_{-1,q^{2}}$. Thus $t$ is in the same orbit as $t_{-1,q^{2}}$ or $t_{z,q^{2}}$ for $z \neq \pm 1, q^{\pm 2}, q^{-4}, -q^{-2}, \pm q^{-1}$. \\ $\square$

\[\beginpicture
\setcoordinatesystem units <0.3cm,0.3cm>
\setplotarea x from -25 to 25, y from -3 to 12
\setlinear
\setplotsymbol({{\scriptsize .}})
\thicklines
\plot -19 -4 -19 15 /
\plot -23 0 -4 0 /
\plot -23 -4 -4 15 /
\plot -23 15 -4 -4 /
\plot -23 4 -15 -4 /
\plot -23 11 -4 11 /
\plot -8 -4 -8 15 /
\plot 6 -4 6 15 /
\plot 2 0 21 0 /
\plot 2 -4 21 15 /
\plot 2 15 21 -4 /
\plot 2 4 10 -4 /
\plot 2 11 21 11 /
\plot 17 -4 17 15 /
\setdashes
\plot -23 6 -13 -4 /
\plot -23 13 -6 -4 /
\plot -17 15 -17 -4 /
\plot -23 -2 -6 15 /
\plot -23 2 -4 2 /
\plot 2 2 21 2 /
\plot 2 6 12 -4 /
\plot 2 13 19 -4 /
\plot 8 15 8 -4 /
\plot 2 -2 19 15 /
\put{$\newmoon$} at -19 0
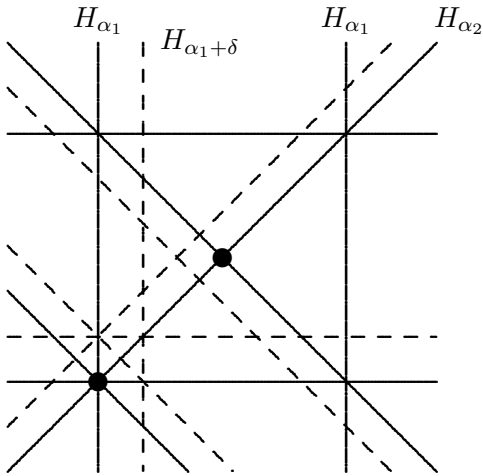
\put{Case 1: $|Z(t)| \geq 2$, so that $t$} at -15 -6
\put{lies on at least two hyperplanes $H_{\alpha}$.} at -14.5 -7.5
\put{$\newmoon$} at -13.5 5.5
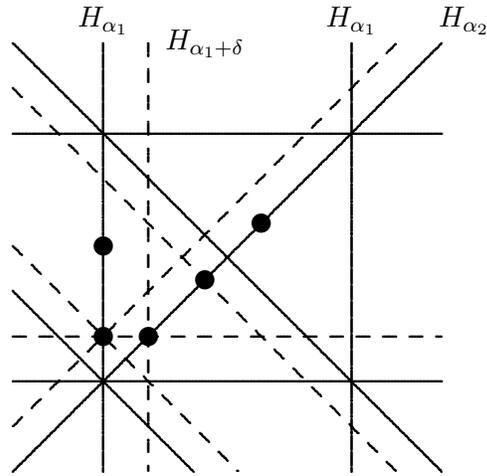
\put{Case 2: $|Z(t)| = 1$, so that $t$} at 10 -6
\put{lies on exactly one hyperplane $H_{\alpha}$.} at 10.5 -7.5
\put{$\newmoon$} at 6 2
\put{$\newmoon$} at 6 6
\put{$\newmoon$} at 13 7
\put{$\newmoon$} at 8 2
\put{$\newmoon$} at 10.5 4.5
\put{$H_{\alpha_{1}}$} at -19 16
\put{$H_{\alpha_{1}}$} at 6 16
\put{$H_{\alpha_{1}}$} at -8 16
\put{$H_{\alpha_{1}}$} at 17 16
\put{$H_{\alpha_{2}}$} at -3 16
\put{$H_{\alpha_{2}}$} at 22 16
\put{$H_{\alpha_{1} + \delta}$} at -14.5 15
\put{$H_{\alpha_{1} + \delta}$} at 10.5 15
\endpicture\]

\[\beginpicture
\setcoordinatesystem units <0.3cm,0.3cm>
\setplotarea x from -25 to 25, y from -3 to 12
\setlinear
\setplotsymbol({{\scriptsize .}})
\thicklines
\plot -19 -4 -19 15 /
\plot -23 0 -4 0 /
\plot -23 -4 -4 15 /
\plot -23 15 -4 -4 /
\plot -23 4 -15 -4 /
\plot -23 11 -4 11 /
\plot -8 -4 -8 15 /
\plot 6 -4 6 15 /
\plot 2 0 21 0 /
\plot 2 -4 21 15 /
\plot 2 15 21 -4 /
\plot 2 4 10 -4 /
\plot 2 11 21 11 /
\plot 17 -4 17 15 /
\setdashes
\plot -23 6 -13 -4 /
\plot -23 13 -6 -4 /
\plot -17 15 -17 -4 /
\plot -23 -2 -6 15 /
\plot -23 2 -4 2 /
\plot 2 2 21 2 /
\plot 2 6 12 -4 /
\plot 2 13 19 -4 /
\plot 8 15 8 -4 /
\plot 2 -2 19 15 /
\put{$\newmoon$} at -17 6
\put{$\newmoon$} at -17 4
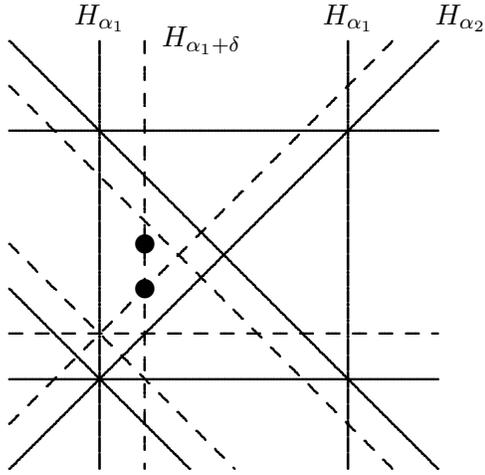
\put{Case 3: $Z(t) = \emptyset, \alpha_{1} \in P(t) $, so that $t$} at -15 -6
\put{lies on $H_{\alpha_{1} \pm \delta}$, but not on any $H_{\alpha}$.} at -15 -7.5
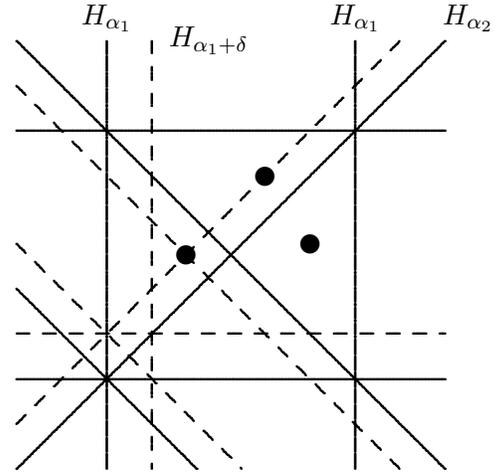
\put{Case 3: $Z(t) = \emptyset, \alpha_{1} \notin P(t) $, so that $t$} at 10 -6
\put{does not lie on $H_{\alpha_{1} \pm \delta}$ or any $H_{\alpha}$.} at 10 -7.5
\put{$\newmoon$} at 13 9
\put{$\newmoon$} at 9.5 5.5
\put{$\newmoon$} at 15 6
\put{$H_{\alpha_{1}}$} at -19 16
\put{$H_{\alpha_{1}}$} at 6 16
\put{$H_{\alpha_{1}}$} at -8 16
\put{$H_{\alpha_{1}}$} at 17 16
\put{$H_{\alpha_{2}}$} at -3 16
\put{$H_{\alpha_{2}}$} at 22 16
\put{$H_{\alpha_{1} + \delta}$} at -14.5 15
\put{$H_{\alpha_{1} + \delta}$} at 10.5 15
\endpicture\]

$\quad \square$

\setlength{\unitlength}{1in}

\noindent \textbf{Remark.} If $q^{2}$ is a root of unity of order less than or equal to 4, there is redundancy in the list of characters given above. Essentially, this is a result of the periodicity in $T_{Q}$ when $q^{2}$ is an $\ell$th root of unity. If $q^{2}$ is a primitive fourth root of unity, then $t_{q^{2},q^{2}} = s_{1}s_{2}s_{1}t_{-1,q^{2}}$.

If $q^{2}$ is a primitive third root of unity, $t_{q^{2},q^{2}} = s_{2}s_{1}s_{2}t_{q^{2},1}$. Also, one of $t_{q,1}$ and $t_{-q,1}$ is equal to $t_{q^{-2},1}$ and is in the same orbit as $t_{q^{2},1}$. (Which one it is depends on whether $q^{3} = 1$ or $-1$. In either case, $t_{q^{2}+1,1}$ is in a different orbit than $t_{q^{2},1}$, so $t_{q^{2}+1,1}$ is our preferred notation for this character.)

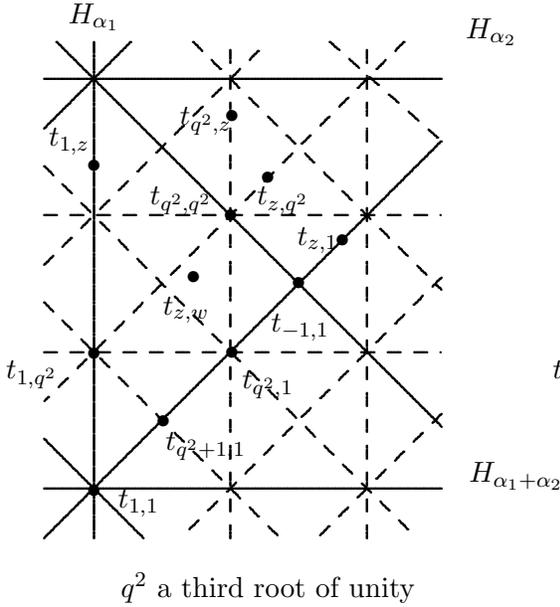
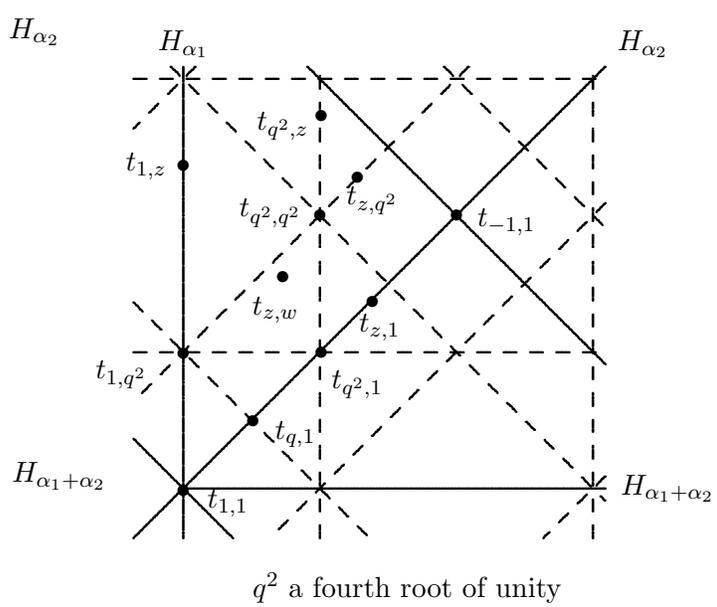
\[\beginpicture
\setcoordinatesystem units <1.3in,1.3in>
\setplotarea x from 0.8 to 5.6, y from 0.8 to 2.8
\setlinear
\setplotsymbol({{\scriptsize .}})
\plot 5 2.7 3.1 0.8 /
\plot 3.3 1 5 1 /
\plot 3.5 0.8 3.1 1.2 /
\plot 3.3 0.8 3.3 2.7 /
\plot 3.8 2.7 5 1.5 /
\setdashes
\plot 3.85 0.8 3.85 2.7 /
\plot 4.95 0.8 4.95 2.7 /
\plot 3.1 2.65 5 2.65 /
\plot 3.1 1.55 5 1.55 /
\plot 3.25 2.7 5 0.95 /
\plot 3.1 1.75 4.05 0.8 /
\plot 4.35 2.7 5 2.05 /
\plot 3.35 2.7 3.1 2.45 /
\plot 4.45 2.7 3.1 1.35 /
\plot 5 2.15 3.65 0.8 /
\plot 5 1.05 4.75 0.8 /
\put{$\bullet$} at 3.3 .99
\put{$\bullet$} at 4.4 2.1
\put{$\bullet$} at 3.3 2.3
\put{$\bullet$} at 3.3 1.54
\put{$\bullet$} at 3.855 1.545
\put{$\bullet$} at 3.58 1.27
\put{$\bullet$} at 4.06 1.75
\put{$\bullet$} at 3.85 2.1
\put{$\bullet$} at 3.855 2.5
\put{$\bullet$} at 4 2.25
\put{$\bullet$} at 3.7 1.85
\put{$t_{1,1}$} at 3.48 .94
\put{$t_{-1,1}$} at 4.6 2.075
\put{$t_{1,z}$} at 3.15 2.3
\put{$t_{1,q^{2}}$} at 3.05 1.46
\put{$t_{q^{2},1}$} at 4 1.42
\put{$t_{q,1}$} at 3.75 1.22
\put{$t_{z,1}$} at 4.09 1.65
\put{$t_{q^{2},q^{2}}$} at 3.65 2.11
\put{$t_{q^{2},z}$} at 3.7 2.46
\put{$t_{z,q^{2}}$} at 4.06 2.16
\put{$t_{z,w}$} at 3.67 1.72
\put{$H_{\alpha_{1}}$} at 3.3 2.79
\put{$H_{\alpha_{1}+\alpha_{2}}$} at 5.25 1
\put{$H_{\alpha_{2}}$} at 5.15 2.79
\put{$q^{2}$ a fourth root of unity} at 4.2 0.6
\setsolid
\plot 2.5 2.4 0.9 0.8 /
\plot 0.9 1 2.5 1 /
\plot 1.3 0.8 0.9 1.2 /
\plot 1.1 0.8 1.1 2.8 /
\plot 0.9 2.65 2.5 2.65 /
\plot 0.95 2.8 2.5 1.25 /
\plot 1.25 2.8 0.9 2.45 /
\setdashes
\plot 0.9 1.9 1.8 2.8 /
\plot 2.2 0.8 2.2 2.8 /
\plot 1.5 2.8 2.5 1.8 /
\plot 1.65 0.8 1.65 2.8 /
\plot 0.9 2.1 2.5 2.1 /
\plot 0.9 1.55 2.5 1.55 /
\plot 0.9 2.3 2.4 0.8 /
\plot 0.9 1.75 1.85 0.8 /
\plot 2.05 2.8 2.5 2.4 /
\plot 2.35 2.8 0.9 1.35 /
\plot 2.5 1.85 1.45 0.8 /
\plot 2.5 1.3 2 0.8 /
\put{$\bullet$} at 1.1 .99
\put{$\bullet$} at 2.1 2
\put{$\bullet$} at 1.1 2.3
\put{$\bullet$} at 1.1 1.54
\put{$\bullet$} at 1.655 1.545
\put{$\bullet$} at 1.38 1.27
\put{$\bullet$} at 1.925 1.825
\put{$\bullet$} at 1.65 2.1
\put{$\bullet$} at 1.655 2.5
\put{$\bullet$} at 1.8 2.25
\put{$\bullet$} at 1.5 1.85
\put{$t_{1,1}$} at 1.28 .94
\put{$t_{z,1}$} at 2 2
\put{$t_{1,z}$} at 1 2.4
\put{$t_{1,q^{2}}$} at 0.85 1.46
\put{$t_{q^{2},1}$} at 1.8 1.42
\put{$t_{q^{2}+1,1}$} at 1.55 1.17
\put{$t_{-1,1}$} at 1.93 1.65
\put{$t_{q^{2},q^{2}}$} at 1.45 2.16
\put{$t_{q^{2},z}$} at 1.55 2.48
\put{$t_{z,q^{2}}$} at 1.86 2.16
\put{$t_{z,w}$} at 1.47 1.72
\put{$H_{\alpha_{1}}$} at 1.1 2.9
\put{$H_{\alpha_{1}+\alpha_{2}}$} at 2.8 1.05
\put{$H_{\alpha_{2}}$} at 2.7 2.84
\put{$q^{2}$ a third root of unity} at 1.8 0.6
\endpicture \]

\clearpage

If $q^{2} = -1$, then $t_{-1,1} = t_{q^{2},1}$, and $t_{1,q^{2}}$ is in the same orbit as $t_{q^{2},q^{2}} = t_{-1,q^{2}}$. Also in this case, $t_{z,q^{2}} = t_{z,-1} = s_{2}t_{-z,-1}$. Finally, if $q = -1$, we have $t_{1,1} = t_{q^{2},1} = t_{1,q^{2}} = t_{q^{2},q^{2}} = t_{-q,1}$. Also, $t_{-1,1} = t_{-1,q^{2}}$, while $t_{q^{2},z} = t_{1,z}$ and $t_{z,q^{2}} = t_{z,1}$.

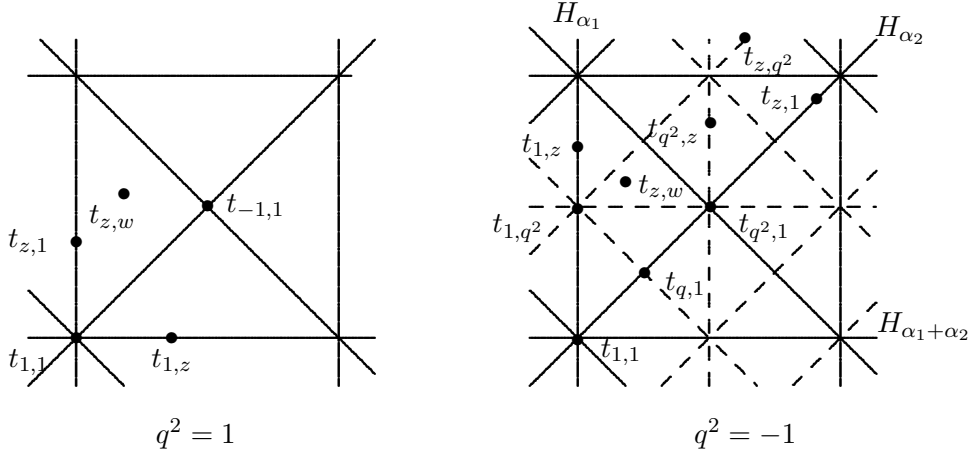
\[\beginpicture
\setcoordinatesystem units <1.25in,1.25in>
\setplotarea x from 1 to 5, y from 0.8 to 2.4
\setlinear
\setplotsymbol({{\scriptsize .}})
\thicklines
\plot 4.75 2.25 3.3 0.8 /
\plot 3.3 1 4.75 1 /
\plot 3.5 0.8 3.5 2.25 /
\plot 4.6 0.8 4.6 2.25 /
\plot 3.3 1.9 3.65 2.25 /
\plot 3.3 1.2 3.7 0.8 /
\plot 3.3 2.1 4.75 2.1 /
\plot 3.3 2.3 4.75 0.85 /
\plot 4.45 2.25 4.75 1.95 /
\setdashes
\plot 3.9 2.25 4.75 1.4 /
\plot 4.05 0.8 4.05 2.25 /
\plot 3.3 1.55 4.75 1.55 /
\plot 3.3 1.75 4.25 0.8 /
\plot 4.2 2.25 3.3 1.35 /
\plot 4.75 1.7 3.85 0.8 /
\plot 4.75 1.15 4.4 0.8 /
\put{$\bullet$} at 3.5 .99
\put{$\bullet$} at 4.5 2
\put{$\bullet$} at 3.5 1.8
\put{$\bullet$} at 3.5 1.54
\put{$\bullet$} at 4.055 1.545
\put{$\bullet$} at 3.78 1.27
\put{$\bullet$} at 4.055 1.9
\put{$\bullet$} at 4.2 2.255
\put{$\bullet$} at 3.7 1.65
\put{$t_{1,1}$} at 3.68 .94
\put{$t_{z,1}$} at 4.35 2
\put{$t_{1,z}$} at 3.35 1.8
\put{$t_{1,q^{2}}$} at 3.25 1.46
\put{$t_{q^{2},1}$} at 4.28 1.46
\put{$t_{q,1}$} at 3.95 1.22
\put{$t_{q^{2},z}$} at 3.9 1.86
\put{$t_{z,q^{2}}$} at 4.3 2.16
\put{$t_{z,w}$} at 3.85 1.63
\put{$H_{\alpha_{1}}$} at 3.5 2.35
\put{$H_{\alpha_{1}+\alpha_{2}}$} at 4.95 1.05
\put{$H_{\alpha_{2}}$} at 4.85 2.29
\put{$q^{2}=-1$} at 4.2 0.6
\setlinear
\setsolid
\setplotsymbol({{\scriptsize .}})
\thicklines
\plot 1.2 1 2.65 1 /
\plot 1.4 0.8 1.4 2.25 /
\plot 1.2 0.8 2.65 2.25 /
\plot 2.5 0.8 2.5 2.25 /
\plot 1.2 2.1 2.55 2.1 /
\plot 1.2 1.2 1.6 0.8 /
\plot 1.25 2.25 2.65 0.85 /
\put{$\bullet$} at 1.4 1
\put{$\bullet$} at 1.95 1.55
\put{$\bullet$} at 1.4 1.4
\put{$\bullet$} at 1.6 1.6
\put{$\bullet$} at 1.8 1
\put{$t_{1,1}$} at 1.2 0.9
\put{$t_{-1,1}$} at 2.15 1.55
\put{$t_{1,z}$} at 1.8 0.9
\put{$t_{z,1}$} at 1.2 1.4
\put{$t_{z,w}$} at 1.55 1.5
\put{$q^{2}=1$} at 1.9 0.6
\endpicture \]

\vspace{.1in}

\noindent \textbf{Analysis of the Characters.} \label{C2analysis} \begin{theorem} \label{eq:B2ones}
The 1-dimensional representations of $\widetilde{H}$ are

\begin{center}

\vspace{.3in}
\renewcommand{\baselinestretch}{1} \normalsize
$\begin{array}{ccc}  L_{q,q,1}: \widetilde{H} \rightarrow \mathbb{C} \quad \,\,\,\,\,\,\, \qquad & L_{q,q,-1}: \widetilde{H} \rightarrow \mathbb{C} \qquad & L_{q,-q^{-1},1}: \widetilde{H} \rightarrow \mathbb{C} \qquad \\ T_{1} \mapsto q & T_{1} \mapsto q & T_{1} \mapsto q \\ T_{2} \mapsto q & T_{2} \mapsto q & T_{2} \mapsto -q^{-1} \\ X^{\omega_{1}} \mapsto q^{3} & X^{\omega_{1}} \mapsto -q^{3} & X^{\omega_{1}} \mapsto q \\ X^{\omega_{2}} \mapsto q^{4} & X^{\omega_{2}} \mapsto q^{4} & X^{\omega_{2}} \mapsto 1 \end{array}$

\vspace{.2in}

$\begin{array}{ccc} L_{q,-q^{-1},-1}: \widetilde{H} \rightarrow \mathbb{C} \qquad & L_{-q^{-1},q,1}: \widetilde{H} \rightarrow \mathbb{C} \quad & L_{-q^{-1},q,-1}: \widetilde{H} \rightarrow \mathbb{C} \quad \\  T_{1} \mapsto q & T_{1} \mapsto -q^{-1} & T_{1} \mapsto -q^{-1} \\ T_{2} \mapsto -q^{-1} & T_{2} \mapsto q & T_{2} \mapsto q \\ X^{\omega_{1}} \mapsto -q & X^{\omega_{1}} \mapsto q^{-1} & X^{\omega_{1}} \mapsto -q^{-1} \\ X^{\omega_{2}} \mapsto 1 & X^{\omega_{2}} \mapsto 1 & X^{\omega_{2}} \mapsto 1 \end{array}$

\vspace{.2in}

$\begin{array}{cc}  L_{-q^{-1},-q^{-1},1}: \widetilde{H} \rightarrow \mathbb{C} \quad & L_{-q^{-1},-q^{-1},-1}: \widetilde{H} \rightarrow \mathbb{C} \\ T_{1} \mapsto -q^{-1} & T_{1} \mapsto -q^{-1}   \\ T_{2} \mapsto -q^{-1} &T_{2} \mapsto -q^{-1} \\ X^{\omega_{1}} \mapsto q^{-3} & X^{\omega_{1}} \mapsto -q^{-3} \\ X^{\omega_{2}} \mapsto q^{-4} & X^{\omega_{2}} \mapsto q^{-4} \end{array}$ \renewcommand{\baselinestretch}{1} \normalsize

\end{center}

\end{theorem}

\noindent \textit{Proof.} A straightforward check shows that the maps above respect the defining relations for $\widetilde{H}$ (\eqref{eq:quadrel} - \eqref{eq:relation}), so that the maps are homomorphisms. Now let $\mathbb{C}v$ be any 1-dimensional $\widetilde{H}$-module. By $\eqref{eq:quadrel}$, \[T_{i}v = qv \quad \textrm{or} \quad T_{i}v = -q^{-1}v,\] but these choices can be made independently for $T_{1}$ and $T_{2}$.

Once the actions of the $T_{i}$ are known, the relation $\eqref{eq:relation}$ forces the actions of the $X^{\alpha_{i}}$. This then defines the action of the $X^{\omega_{i}}$ up to a choice of square root. $\square$

\noindent
\underline{Remark:} If $q$ is a primitive fourth root of unity, then $L_{q,q,\pm 1} \cong L_{q,-q^{-1},\mp 1} \cong L_{-q^{-1},q, \pm 1} \cong L_{-q^{-1},-q^{-1},\mp 1}$.

\vspace{.2in}

Let $t \in T$. The principal series module is \[ M(t) = \textrm{Ind}_{\mathbb{C}[X]}^{\widetilde{H}} \mathbb{C}_{t} = \widetilde{H} \otimes_{\mathbb{C}[X]} \mathbb{C}_{t},\] where $\mathbb{C}_{t}$ is the one-dimensional $\mathbb{C}[X]$-module given by
\[ \mathbb{C}_{t} = \textrm{span}\{v_{t}\} \quad \textrm{ and } \quad X^{\lambda}v_{t} = t(X^{\lambda})v_{t}.\]
By Theorem \ref{PrinQuot}c, every irreducible $\widetilde{H}$ module with central character $t$ is a composition factor of $M(t)$, and thus it is sufficient to find all composition factors of $M(t)$ for all possible weights $t$.

A weight $t|_{Q}$ corresponds to a point in the root lattice $Q$ as described above. The composition structure of $M(t)$ is largely determined by the structure of the operators $\tau_{i}$, which can be encoded in the picture of a small neighborhood of the point $t$ in $Q$. The following pictures show a small neighborhood of $t$, including the hyperplanes $H_{\alpha_{i}}$ that include the point $t$. Each picture shows 8 chambers around $t$, (bounded by the affine hyperplanes through $t$), representing the 8 (not necessarily distinct) elements of $W_{0}t$. The solid lines are hyperplanes $H_{\alpha_{i}}$ with $t(X^{\alpha_{i}}) = 1$, so that $s_{\alpha_{i}} \in W_{t}$. The chambers between any two successive solid lines are in bijection with the distinct elements of the orbit $W_{0}t$, which are in bijection with the cosets $W_{0}/W_{t}$. These weight spaces are $|W_{t}|$-dimensional, so each chamber has $|W_{t}|$ points drawn in it. The dotted lines represent hyperplanes $H_{\alpha_{i}}$ such that $t(X^{\alpha_{i}}) = q^{\pm 2}$.

The operators $\tau_{i}:M_{t} \rightarrow M_{s_{i}t}$ and $\tau_{i}: M_{s_{i}t} \rightarrow M_{t}$ exist and are invertible exactly when $t(X^{\alpha_{i}}) \neq 1$ or $q^{2}$, i.e. exactly when the corresponding chambers in our picture do not have a hyperplane between them. Thus the hyperplanes in the picture of the neighborhood of $t$ show which $\tau$ operators are invertible and which are not. In most cases, this is enough to determine the exact composition factors of $M(t)$.

\underline{Case 1:} $P(t) = \emptyset$

If $P(t) = \emptyset$, then by Kato's criterion (Theorem $\ref{eq:Kato}$), $M(t)$ is irreducible and is the only irreducible module with central character $t$. The weights of $M(t)$ are in bijection with $W_{0}/W_{t}$, the cosets of the centralizer of $t$ in $W$, and dim$(M(t)_{wt}) = |W_{t}|$. If $w$ and $s_{i}w$ are distinct weights in $M(t)$, then $\tau_{i}:M(t)_{t} \rightarrow M(t)_{s_{i}t}$ is a bijection.

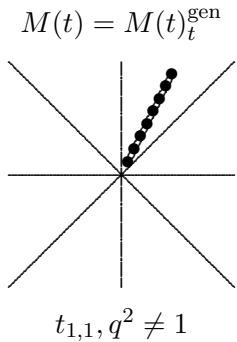
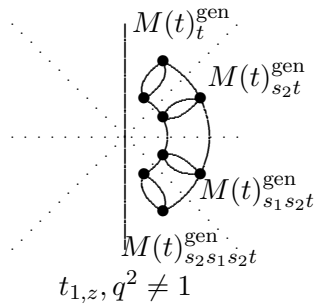
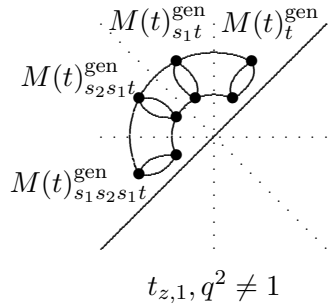
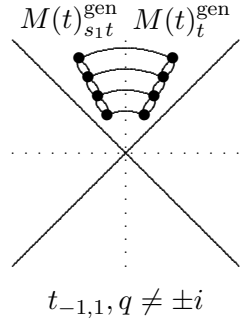
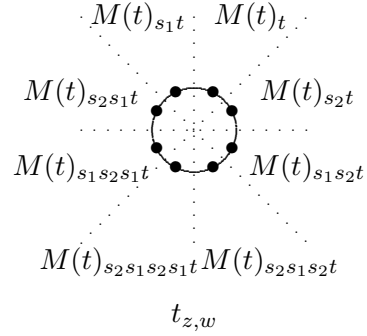
\vspace{.1in}

\[ \beginpicture
\setcoordinatesystem units <0.25cm,0.25cm>
\setplotarea x from -25 to 25, y from -6 to 4
\linethickness=3pt
\put{$\bullet$} at -18.666 0.6666
\put{$\bullet$} at -18 2
\put{$\bullet$} at -17.666 2.6666
\put{$\bullet$} at -17.333 3.3333
\put{$\bullet$} at -17 4
\put{$\bullet$} at -18.333 1.3333
\put{$\bullet$} at -16.666 4.6666
\put{$\bullet$} at -16.333 5.3333
\put{$\bullet$} at 4 2
\put{$\bullet$} at 5 4
\put{$\bullet$} at 5 1
\put{$\bullet$} at 7 2
\put{$\bullet$} at 5 -1
\put{$\bullet$} at 7 -2
\put{$\bullet$} at 4 -2
\put{$\bullet$} at 5 -4
\put{$\bullet$} at 26 2
\put{$\bullet$} at 27 4
\put{$\bullet$} at 23 1
\put{$\bullet$} at 21 2
\put{$\bullet$} at 23 -1
\put{$\bullet$} at 21 -2
\put{$\bullet$} at 24 2
\put{$\bullet$} at 23 4
\put{$t_{1,1}, q^{2} \neq 1$} at -19 -8
\put{$t_{1,z},q^{2} \neq 1$} at 3 -8
\put{$t_{z,1},q^{2} \neq 1$} at 25 -8
\put{$M(t) = M(t)_{t}^{\textrm{gen}}$} at -19 8
\put{$M(t)_{t}^{\textrm{gen}}$} at 6 6
\put{$M(t)_{s_{2}t}^{\textrm{gen}}$} at 10 3
\put{$M(t)_{s_{1}s_{2}t}^{\textrm{gen}}$} at 10 -3
\put{$M(t)_{s_{2}s_{1}s_{2}t}^{\textrm{gen}}$} at 6.4 -6
\put{$M(t)_{t}^{\textrm{gen}}$} at 28 6
\put{$M(t)_{s_{1}t}^{\textrm{gen}}$} at 22 6
\put{$M(t)_{s_{2}s_{1}t}^{\textrm{gen}}$} at 17.8 3
\put{$M(t)_{s_{1}s_{2}s_{1}t}^{\textrm{gen}}$} at 17.8 -2.6
\setlinear
\plot -19 -6 -19 6 /
\plot -25 0 -13 0 /
\plot -25 6 -13 -6 /
\plot -25 -6 -13 6 /
\plot 3 6 3 -6 /
\plot 19 -6 31 6 /
\setdots
\plot -3 6 9 -6 /
\plot -3 -6 9 6 /
\plot -3 0 9 0 /
\plot 25 6 25 -6 /
\plot 19 0 31 0 /
\plot 19 6 31 -6 /
\setdashes
\setsolid
\circulararc 130 degrees from 4 -2 center at 3 0
\circulararc 130 degrees from 5 -4 center at 3 0
\circulararc 140 degrees from 26 2 center at 25 0
\circulararc 140 degrees from 27 4 center at 25 0
\setquadratic
\plot -18.666 0.6666 -18.666 1 -18.333 1.3333 /
\plot -18.666 0.6666 -18.333 1 -18.333 1.3333 /
\plot -18.333 1.3333 -18.333 1.666 -18 2 /
\plot -18.333 1.3333 -18 1.666 -18 2 /
\plot -18 2 -18 2.3333 -17.666 2.6666 /
\plot -18 2 -17.666 2.3333 -17.666 2.6666 /
\plot -17.666 2.6666 -17.666 3 -17.333 3.3333 /
\plot -17.666 2.6666 -17.333 3 -17.333 3.3333 /
\plot -17.333 3.3333 -17.333 3.6666 -17 4 /
\plot -17.333 3.3333 -17 3.6666 -17 4 /
\plot -17 4 -17 4.3333 -16.666 4.6666 /
\plot -17 4 -16.666 4.3333 -16.666 4.6666 /
\plot -16.666 4.6666 -16.6666 5 -16.3333 5.3333 /
\plot -16.666 4.6666 -16.3333 5 -16.3333 5.3333 /
\plot 4 2 4 3 5 4 /
\plot 4 2 5 3 5 4 /
\plot 5 1 5.7 1.8 7 2 /
\plot 5 1 5.9 1.1 7 2 /
\plot 5 -1 5.7 -1.8 7 -2 /
\plot 5 -1 5.9 -1.1 7 -2 /
\plot 4 -2 4 -3 5 -4 /
\plot 4 -2 5 -3 5 -4 /
\plot 26 2 26 3 27 4 /
\plot 26 2 27 3 27 4 /
\plot 23 1 21.7 1.2 21 2 /
\plot 23 1 22.1 1.9 21 2 /
\plot 23 -1 21.7 -1.2 21 -2 /
\plot 23 -1 22.1 -1.9 21 -2 /
\plot 24 2 24 3 23 4 /
\plot 24 2 23 3 23 4 /
\endpicture \]

\[\beginpicture
\setcoordinatesystem units <0.25cm,0.25cm>
\setplotarea x from -25 to 25, y from -6 to 6
\linethickness=3pt
\put{$\bullet$} at -10 2
\put{$\bullet$} at -9.5 3
\put{$\bullet$} at -9 4
\put{$\bullet$} at -8.5 5
\put{$\bullet$} at -12 2
\put{$\bullet$} at -12.5 3
\put{$\bullet$} at -13 4
\put{$\bullet$} at -13.5 5
\put{$\bullet$} at 18 2
\put{$\bullet$} at 19 1
\put{$\bullet$} at 19 -1
\put{$\bullet$} at 18 -2
\put{$\bullet$} at 16 2
\put{$\bullet$} at 15 1
\put{$\bullet$} at 15 -1
\put{$\bullet$} at 16 -2
\put{$t_{-1,1}, q \neq \pm i$} at -11 -8
\put{$t_{z,w}$} at 17 -10
\put{$M(t)_{t}^{\textrm{gen}}$} at -8 7
\put{$M(t)_{s_{1}t}^{\textrm{gen}}$} at -14 7
\put{$M(t)_{t} $} at 20 6
\put{$M(t)_{s_{2}t} $} at 23 2
\put{$M(t)_{s_{1}s_{2}t} $} at 23 -2
\put{$M(t)_{s_{2}s_{1}s_{2}t} $} at 21 -7
\put{$M(t)_{s_{2}s_{1}s_{2}s_{1}t} $} at 13 -7
\put{$M(t)_{s_{1}s_{2}s_{1}t} $} at 11 -2
\put{$M(t)_{s_{2}s_{1}t} $} at 11 2
\put{$M(t)_{s_{1}t} $} at 14 6
\setlinear
\plot -17 6 -5 -6 /
\plot -17 -6 -5 6 /
\setdots
\plot -11 -6 -11 6 /
\plot -17 0 -5 0 /
\plot 17 6 17 -6 /
\plot 11 0 23 0 /
\plot 11 6 23 -6 /
\plot 11 -6 23 6 /
\setdashes
\setsolid
\setquadratic
\plot -10 2 -10 2.5 -9.5 3 /
\plot -10 2 -9.5 2.5 -9.5 3 /
\plot -9.5 3 -9.5 3.5 -9 4 /
\plot -9.5 3 -9 3.5 -9 4 /
\plot -9 4 -9 4.5 -8.5 5 /
\plot -9 4 -8.5 4.5 -8.5 5 /
\plot -12 2 -12 2.5 -12.5 3 /
\plot -12 2 -12.5 2.5 -12.5 3 /
\plot -12.5 3 -12.5 3.5 -13 4 /
\plot -12.5 3 -13 3.5 -13 4 /
\plot -13 4 -13 4.5 -13.5 5 /
\plot -13 4 -13.5 4.5 -13.5 5 /
\circulararc 360 degrees from 16 -2 center at 17 0
\circulararc 55 degrees from -10 2 center at -11 0
\circulararc 55 degrees from -9.5 3 center at -11 0
\circulararc 55 degrees from -9 4 center at -11 0
\circulararc 55 degrees from -8.5 5 center at -11 0
\endpicture \]

\vspace{.2in}

\underline{Case 2:} $Z(t) = \emptyset$, but $P(t) \neq \emptyset$.

If $Z(t) = \emptyset$ then $t$ is a regular central character. Then the irreducibles with central character $t$ are in bijection with the connected components of the calibration graph for $t$, and can be constructed using Theorem $\ref{eq:locreg}$. The components of the calibration graph are the same as the sets of chambers between any two successive dotted hyperplanes in the picture of $t$. The segments in the picture represent the fact that the $\tau$ operators between the two corresponding weight vectors are invertible, and thus those vectors must lie in the same composition factor. (See proposition $\ref{calibsame}$b.)

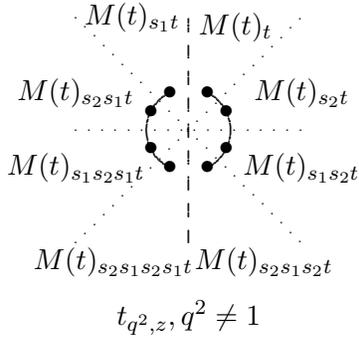
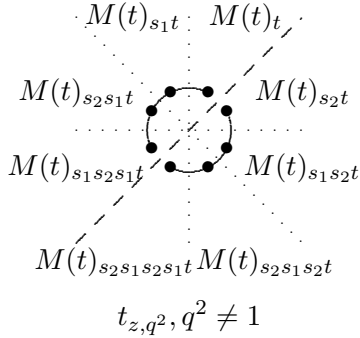
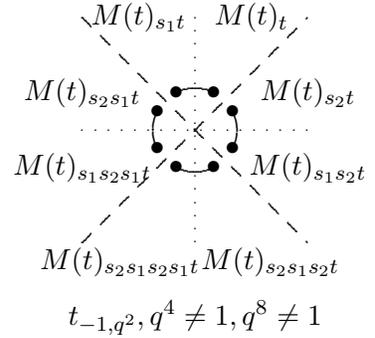
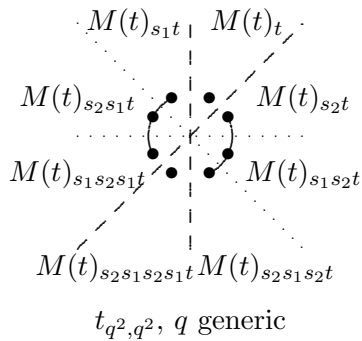
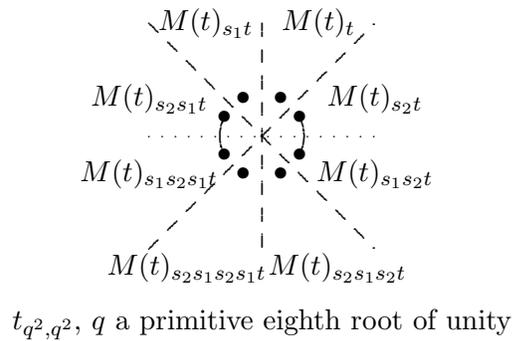
\[ \beginpicture
\setcoordinatesystem units <0.25cm,0.25cm>
\setplotarea x from -25 to 25, y from -6 to 6
\linethickness=3pt
\put{$\bullet$} at -17 -1
\put{$\bullet$} at -18 2
\put{$\bullet$} at -18 -2
\put{$\bullet$} at -20 2
\put{$\bullet$} at -17 1
\put{$\bullet$} at -21 1
\put{$\bullet$} at -21 -1
\put{$\bullet$} at -20 -2
\put{$\bullet$} at 4 2
\put{$\bullet$} at 2 2
\put{$\bullet$} at 5 1
\put{$\bullet$} at 1 1
\put{$\bullet$} at 5 -1
\put{$\bullet$} at 1 -1
\put{$\bullet$} at 4 -2
\put{$\bullet$} at 2 -2
\put{$\bullet$} at 26 2
\put{$\bullet$} at 24 2
\put{$\bullet$} at 23 1
\put{$\bullet$} at 24 -2
\put{$\bullet$} at 23 -1
\put{$\bullet$} at 27 1
\put{$\bullet$} at 26 -2
\put{$\bullet$} at 27 -1
\put{$t_{q^{2},z}, q^{2} \neq 1$} at -19 -10
\put{$t_{z,q^{2}},q^{2} \neq 1$} at 3 -10
\put{$t_{-1,q^{2}},q^{4} \neq 1, q^{8} \neq 1$} at 25 -10
\put{$M(t)_{t}$} at -16.5 5.5
\put{$M(t)_{s_{2}t} $} at -13 2
\put{$M(t)_{s_{1}s_{2}t} $} at -13 -2
\put{$M(t)_{s_{2}s_{1}s_{2}t} $} at -15 -7
\put{$M(t)_{s_{2}s_{1}s_{2}s_{1}t} $} at -23 -7
\put{$M(t)_{s_{1}s_{2}s_{1}t} $} at -25 -2
\put{$M(t)_{s_{2}s_{1}t} $} at -25 2
\put{$M(t)_{s_{1}t} $} at -22 6
\put{$M(t)_{t} $} at 6 6
\put{$M(t)_{s_{2}t} $} at 9 2
\put{$M(t)_{s_{1}s_{2}t} $} at 9 -2
\put{$M(t)_{s_{2}s_{1}s_{2}t} $} at 7 -7
\put{$M(t)_{s_{2}s_{1}s_{2}s_{1}t} $} at -1 -7
\put{$M(t)_{s_{1}s_{2}s_{1}t} $} at -3 -2
\put{$M(t)_{s_{2}s_{1}t} $} at -3 2
\put{$M(t)_{s_{1}t} $} at 0 6
\put{$M(t)_{t} $} at 28 6
\put{$M(t)_{s_{2}t} $} at 31 2
\put{$M(t)_{s_{1}s_{2}t} $} at 31 -2
\put{$M(t)_{s_{2}s_{1}s_{2}t} $} at 29 -7
\put{$M(t)_{s_{2}s_{1}s_{2}s_{1}t} $} at 21 -7
\put{$M(t)_{s_{1}s_{2}s_{1}t} $} at 19 -2
\put{$M(t)_{s_{2}s_{1}t} $} at 19 2
\put{$M(t)_{s_{1}t} $} at 22 6
\setlinear
\setdots
\plot -25 0 -13 0 /
\plot -25 6 -13 -6 /
\plot -25 -6 -13 6 /
\plot -3 6 9 -6 /
\plot -3 0 9 0 /
\plot 3 6 3 -6 /
\plot 19 0 31 0 /
\plot 25 6 25 -6 /
\setdashes
\plot 19 6 31 -6 /
\plot 19 -6 31 6 /
\plot -3 -6 9 6 /
\plot -19 -6 -19 6 /
\setsolid
\circulararc 130 degrees from -18 -2 center at -19 0
\circulararc 130 degrees from -20 2 center at -19 0
\circulararc 135 degrees from 2 -2 center at 3 0
\circulararc 135 degrees from 4 2 center at 3 0
\circulararc 48 degrees from 27 -1 center at 25 0
\circulararc 48 degrees from 23 1 center at 25 0
\circulararc 48 degrees from 26 2 center at 25 0
\circulararc 48 degrees from 24 -2 center at 25 0
\endpicture \]

\vspace{.2in}

\[\beginpicture
\setcoordinatesystem units <0.25cm,0.25cm>
\setplotarea x from -25 to 25, y from -6 to 6
\linethickness=3pt
\put{$\bullet$} at -10 2
\put{$\bullet$} at -10 -2
\put{$\bullet$} at -9 1
\put{$\bullet$} at -9 -1
\put{$\bullet$} at -12 2
\put{$\bullet$} at -12 -2
\put{$\bullet$} at -13 1
\put{$\bullet$} at -13 -1
\put{$\bullet$} at 18 2
\put{$\bullet$} at 19 1
\put{$\bullet$} at 19 -1
\put{$\bullet$} at 18 -2
\put{$\bullet$} at 16 2
\put{$\bullet$} at 15 1
\put{$\bullet$} at 15 -1
\put{$\bullet$} at 16 -2
\put{$t_{q^{2},q^{2}}$, $q$ generic} at -11 -10
\put{$t_{q^{2},q^{2}}$, $q$ a primitive eighth root of unity} at 17 -10
\put{$M(t)_{t} $} at 20 6
\put{$M(t)_{s_{2}t} $} at 23 2
\put{$M(t)_{s_{1}s_{2}t} $} at 23 -2
\put{$M(t)_{s_{2}s_{1}s_{2}t} $} at 21 -7
\put{$M(t)_{s_{2}s_{1}s_{2}s_{1}t} $} at 13 -7
\put{$M(t)_{s_{1}s_{2}s_{1}t} $} at 11 -2
\put{$M(t)_{s_{2}s_{1}t} $} at 11 2
\put{$M(t)_{s_{1}t} $} at 14 6
\put{$M(t)_{t} $} at -8 6
\put{$M(t)_{s_{2}t} $} at -5 2
\put{$M(t)_{s_{1}s_{2}t} $} at -5 -2
\put{$M(t)_{s_{2}s_{1}s_{2}t} $} at -7 -7
\put{$M(t)_{s_{2}s_{1}s_{2}s_{1}t} $} at -15 -7
\put{$M(t)_{s_{1}s_{2}s_{1}t} $} at -17 -2
\put{$M(t)_{s_{2}s_{1}t} $} at -17 2
\put{$M(t)_{s_{1}t} $} at -14 6
\setlinear
\setdots
\plot -17 6 -5 -6 /
\plot -17 0 -5 0 /
\plot 11 0 23 0 /
\setdashes
\plot 17 6 17 -6 /
\plot 11 6 23 -6 /
\plot 11 -6 23 6 /
\plot -17 -6 -5 6 /
\plot -11 -6 -11 6 /
\setsolid
\setquadratic
\circulararc 55 degrees from 19 -1 center at 17 0
\circulararc 55 degrees from 15 1 center at 17 0
\circulararc 95 degrees from -10 -2 center at -11 0
\circulararc 95 degrees from -12 2 center at -11 0
\endpicture \]

\underline{Case 3:} $Z(t) \neq \emptyset, P(t) \neq \emptyset$.

The only central characters not covered in cases 1 and 2 are those in the orbits of $t_{1,q^{2}}$, $t_{q^{2},1}$, and $t_{\pm q,1}$.

\noindent
\underline{$t|_{Q} = t_{q^{2},1}$.}

If $q^{2} = 1$, then the results of section $\ref{clifford}$ show that $\widetilde{H}$ has five irreducible representations - four of them 1-dimensional, and one 2-dimensional. Specifically, an irreducible $\widetilde{H}$ module is an irreducible $W_{0}$-module (via the identification $\mathbb{C}[W_{0}] = H $) on which $X^{\lambda}$ acts by the constant $t(X^{\lambda})$, and $M(t)$ is isomorphic to the regular representation of $\mathbb{C}[W_{0}]$ as a $W_{0}$-module.

Assume $q^{2} \neq 1$ and let
\renewcommand{\baselinestretch}{1} \normalsize \[ w_{1} = \begin{cases} s_{1} & \textrm{ if } q^{2} = -1 \\ s_{1}s_{2}s_{1} & \textrm{ if } q^{2} \neq -1 \end{cases}.\] \renewcommand{\baselinestretch}{1} \normalsize
Then let $\mathbb{C}_{q^{2},1}$ and $\mathbb{C}_{q^{-2},1}$ be the 1-dimensional $\widetilde{H}_{\{1\}}$-modules spanned by $v_{t}$ and $v_{w_{1}t}$, respectively, given by
\renewcommand{\baselinestretch}{1} \normalsize \[ X^{\lambda}v_{t}  = t(X^{\lambda})v_{t}, \quad \textrm{ and } \quad T_{1}v_{t} = qv_{t},  \quad \textrm{ and }\] \[ X^{\lambda}v_{w_{1}t} = (w_{1}t)(X^{\lambda})v_{w_{1}t}  \quad  \textrm{ and } \quad T_{1}v_{w_{1}t} = -q^{-1}v_{w_{1}t}.\] \renewcommand{\baselinestretch}{1} \normalsize
Then \[M = \widetilde{H} \otimes_{\widetilde{H}_{\{1\}}} \mathbb{C}_{q^{2},1} \quad \textrm{ and } \quad N = \widetilde{H} \otimes_{\widetilde{H}_{\{1\}}} \mathbb{C}_{q^{-2},1}\] are 4-dimensional $\widetilde{H}$-modules.

Each dot in the chamber $w^{-1}C$ in the following picture represents a basis element of the $wt$ weight space of $M$ or $N$. The dots that are connected by arcs represent basis vectors in the same module, $M$ or $N$.

\[\beginpicture
\setcoordinatesystem units <0.35cm,0.35cm>
\setplotarea x from -20 to 20, y from -6 to 6
\linethickness=3pt
\put{$\bullet$} at 15 2
\put{$\bullet$} at 16 4
\put{$\bullet$} at 16.5 5
\put{$\bullet$} at 15.5 3
\put{$\bullet$} at 13 2
\put{$\bullet$} at 12 4
\put{$\bullet$} at 12.5 3
\put{$\bullet$} at 11.5 5
\put{$\bullet$} at -10 2
\put{$\bullet$} at -9 4
\put{$\bullet$} at -12 2
\put{$\bullet$} at -13 4
\put{$\bullet$} at -13 1
\put{$\bullet$} at -15 2
\put{$\bullet$} at -13 -1
\put{$\bullet$} at -15 -2
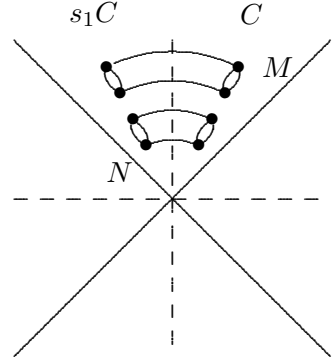
\put{$t_{q^{2},1}$, $q$ a primitive fourth root of unity} at 14 -8
\put{$t_{q^{2},1}$, $q^{2} \neq \pm 1$} at -11 -8
\put{$C$} at 17 7
\put{$s_{1}C$} at 11 7
\put{$C$} at -8 7
\put{$s_{1}C$} at -14 7
\put{$s_{1}s_{2}C$} at -19 3
\put{$s_{1}s_{2}s_{1}C$} at -19 -3
\put{$M$} at -9 5
\put{$N$} at -16 -2
\put{$M$} at 18 5
\put{$N$} at 12 1
\setlinear
\plot 8 6 20 -6 /
\plot 8 -6 20 6 /
\plot -17 -6 -5 6 /
\setdots
\plot -17 6 -5 -6 /
\setdashes
\plot -17 0 -5 0 /
\plot -11 6 -11 -6 /
\plot 8 0 20 0 /
\plot 14 6 14 -6 /
\setsolid
\setquadratic
\plot 15 2 15 2.5 15.5 3 /
\plot 15 2 15.5 2.5 15.5 3 /
\plot 16 4 16 4.5 16.5 5 /
\plot 16 4 16.5 4.5 16.5 5 /
\plot 13 2 13 2.5 12.5 3 /
\plot 13 2 12.5 2.5 12.5 3 /
\plot 12 4 12 4.5 11.5 5 /
\plot 12 4 11.5 4.5 11.5 5 /
\plot -10 2 -10 3 -9 4 /
\plot -10 2 -9 3 -9 4 /
\plot -13 -1 -13.9 -1.9 -15 -2 /
\plot -13 -1 -14.1 -1.1 -15 -2 /
\plot -15 -2 -15 2 -13 4 /
\plot -10 2 -11.5 2.2 -13 1 /
\plot -13 -1 -14 0 -15 2 /
\plot -9 4 -11 3 -12 2 /
\circulararc 50 degrees from 15 2 center at 14 0
\circulararc 50 degrees from 15.5 3 center at 14 0
\circulararc 50 degrees from 16 4 center at 14 0
\circulararc 50 degrees from 16.5 5 center at 14 0
\endpicture \]

\begin{proposition}
If $q^{2} = -1$ and $M = \widetilde{H} \otimes_{\widetilde{H}_{\{1\}}} \mathbb{C}_{q^{2},1} \quad \textrm{ and } \quad N = \widetilde{H} \otimes_{\widetilde{H}_{\{1\}}} \mathbb{C}_{q^{-2},1}$ then

(a) $M$ is irreducible, and

(b) The map \[ \begin{array}{rcl} \phi: N & \,\, \rightarrow \,\, & M \\ hv_{wt} & \mapsto & hv, \quad \textrm{ for } h \in \widetilde{H}\end{array}\] is a $\widetilde{H}$-module isomorphism, where $v = T_{1}T_{2}v_{t} - qT_{2}v_{t} - v_{t} \in M$, and

(c) Any irreducible $\widetilde{H}$-module $L$ with central character $t$ is isomorphic to $M$.

\end{proposition}

\noindent \textit{Proof.} (a) If $q$ is a primitive fourth root of unity, then $M$ has weight spaces $M_{t}^{\textrm{gen}}$, and $M_{s_{1}t}^{\textrm{gen}}$, each of which is 2-dimensional. By Lemma $\ref{eq:lemma}$ and Theorem $\ref{2dims}$, $M$ is irreducible.

(b) Let $v = T_{1}T_{2}v_{t} - qT_{2}v_{t} - v_{t}$. Then a straightforward computation using equation $\ref{eq:relation}$ shows that $v$ spans a 1-dimensional $\widetilde{H}_{\{1\}}$-submodule of $M$, given by \[ T_{1}v = qv, \quad \textrm{ and } \quad X^{\lambda}v = s_{1}t(X^{\lambda})v.\]

Then the $\widetilde{H}_{\{1\}}$-module map given by $v_{wt} \mapsto v$ corresponds to $\phi$ under the adjunction \[ \textrm{Hom}_{\widetilde{H}}(\widetilde{H} \otimes_{\widetilde{H}_{\{1\}}} \mathbb{C}_{q^{-2},1}, M) = \textrm{Hom}_{\widetilde{H}_{\{1\}}}(\mathbb{C}_{q^{-2},1}, M|_{\widetilde{H}_{\{1\}}}).\] Thus $\phi$ is a $\widetilde{H}$-module map and since $M$ is irreducible, the map is surjective. Then since $M$ and $N$ have the same dimension, $\phi$ is an isomorphism.

(c) Let $L$ be an irreducible $\widetilde{H}$-module with central character $t$, which must have weights $t$ and $s_{1}t$. Then, viewing $L$ as a $\widetilde{H}_{\{1\}}$-module, it must have all 1-dimensional composition factors, and it must have a 1-dimensional $\widetilde{H}_{\{1\}}$-submodule, with weight $t$ or $s_{1}t$. Then the same argument as in part b gives an isomorphism from $M$ to $L$ or from $N$ to $L$. $\square$

\vspace{.2in}

\[\beginpicture
\setcoordinatesystem units <0.23cm,0.23cm>
\setplotarea x from -25 to 25, y from -6 to 6
\linethickness=3pt
\put{$\bullet$} at -10 2
\put{$\bullet$} at -15 -2
\put{$\bullet$} at -9 4
\put{$\bullet$} at -13 4
\put{$\bullet$} at -12 2
\put{$\bullet$} at -13 -1
\put{$\bullet$} at -13 1
\put{$\bullet$} at -15 2
\put{$\bullet$} at 18 2
\put{$\bullet$} at 19 4
\put{$\bullet$} at 16 2
\put{$\bullet$} at 15 4
\put{$\bullet$} at 15 1
\put{$\bullet$} at 13 2
\put{$\bullet$} at 15 -1
\put{$\bullet$} at 13 -2
\put{$t_{q^{2},1}$, $q^{2}$ a primitive third root of unity} at -11 -8
\put{$t_{q^{2},1}$, $q$ generic} at 17 -8
\put{$C$} at -8 7
\put{$s_{1}C$} at -14 7
\put{$s_{1}s_{2}C$} at -19 3
\put{$s_{1}s_{2}s_{1}C$} at -19 -3
\put{$C$} at 20 7
\put{$s_{1}C$} at 14 7
\put{$s_{1}s_{2}C$} at 9 3
\put{$s_{1}s_{2}s_{1}C$} at 9 -3
\setlinear
\plot -17 -6 -5 6 /
\plot 11 -6 23 6 /
\setdots
\plot 11 6 23 -6 /
\setdashes
\plot -11 -6 -11 6 /
\plot -17 0 -5 0 /
\plot -17 6 -5 -6 /
\plot 11 0 23 0 /
\plot 17 6 17 -6 /
\setsolid
\setquadratic
\plot -10 2 -10 3 -9 4 /
\plot -10 2 -9 3 -9 4 /
\plot -13 -1 -13.9 -1.9 -15 -2 /
\plot -13 -1 -14.1 -1.1 -15 -2 /
\plot -15 -2 -15.6 0 -15 2 /
\plot -13 -1 -14 0 -15 2 /
\plot -10 2 -11 3.3333 -13 4 /
\plot -9 4 -11 4.6 -13 4 /
\plot 18 2 18 3 19 4 /
\plot 18 2 19 3 19 4 /
\plot 15 -1 14.1 -1.9 13 -2 /
\plot 15 -1 13.9 -1.1 13 -2 /
\plot 13 -2 13 2 15 4 /
\plot 18 2 16.5 2.2 15 1 /
\plot 15 -1 14 0 13 2 /
\plot 19 4 17 3 16 2 /
\endpicture \]

\begin{proposition}

(a) If $q^{2}$ is a primitive third root of unity then $M_{s_{2}s_{1}t}$ is a submodule of $M$ isomorphic to $L_{-q^{-1},q, \pm 1}$ and $M/M_{s_{2}s_{1}t}$ is irreducible. In addition, $N_{s_{1}t}$ is a submodule of $N$ isomorphic to $L_{q,-q^{-1}}$, and $N/N_{s_{1}t}$ is irreducible.

(b) If $q^{2}$ is not $\pm 1$ or a primitive third root of unity then $M$ and $N$ are irreducible and nonisomorphic.

\end{proposition}

\noindent \textit{Proof.} (a) Assume $q^{2}$ is a primitive third root of unity. Then Theorem \ref{weightbasis} shows that $\tau_{2} : M_{s_{1}t} \rightarrow M_{s_{2}s_{1}t}$ is non-zero. But $s_{2}s_{1}t(X^{\alpha_{2}}) = q^{2}$ so that $\tau_{2} : M_{s_{2}s_{1}t} \rightarrow M_{s_{1}t}$ is the zero map by Theorem \ref{tauthm}, and $M_{s_{1}s_{2}t}$ is a submodule of $M$. By Lemma $\ref{eq:lemma}$, $M/M_{s_{1}s_{2}t}$ is irreducible. A parallel argument shows that $N_{s_{1}t}$ is a submodule of $N$, with $N/N_{s_{1}t}$ irreducible.

(b) If $q^{4} \neq 1$ and $q^{6} \neq 1$, then $P(t) = \{ \alpha_{1}, \alpha_{1} + \alpha_{2}\}$. Then Lemma \ref{eq:lemma} shows that the composition factor $M'$ of $M$ with $(M')_{t} \neq 0$ has dim $(M')_{t}^{\textrm{gen}} \geq 2$ and $(M')_{s_{1}t} \neq 0$. Then by Theorem \ref{calibsame}b, $(M')_{s_{2}s_{1}t} \neq 0$, so that $M' = M$. Similarly, Lemma \ref{eq:lemma} and Theorem \ref{calibsame}b show that $N$ is irreducible. Since they have different weight spaces, they are not isomorphic. $\square$

\vspace{.2in}

\underline{$t|_{Q} = t_{1,q^{2}}$}

Note that if $q^{2} = 1$, then $t_{1,q^{2}} = t_{q^{2},1} = t_{1,1}$, so this case has already been addressed.

Let $\mathbb{C}_{1,q^{2}}$ and $\mathbb{C}_{1,q^{-2}}$ be the 1-dimensional $\widetilde{H}_{\{2\}}$-modules spanned by $v_{t}$ and $v_{w_{0}t}$, respectively, and given by

\renewcommand{\baselinestretch}{1} \normalsize
\[ T_{2} v_{t} = qv_{t} \quad \textrm{ and } \quad X^{\lambda}v_{t} = t(X^{\lambda})v_{t}, \quad \textrm{ and } \]

\[ T_{2} v_{w_{0}t} = -q^{-1}v_{w_{0}t} \quad \textrm { and } X^{\lambda}v_{w_{0}t} = w_{0}t(X^{\lambda})v_{w_{0}t}.\] Then \[M = \widetilde{H} \otimes_{\widetilde{H}_{\{2\}}} \mathbb{C}_{1,q^{2}} \quad \textrm{ and } \quad N = \widetilde{H} \otimes_{\widetilde{H}_{\{2\}}} \mathbb{C}_{1,q^{-2}}\] \renewcommand{\baselinestretch}{1} \normalsize
are 4-dimensional $\widetilde{H}$-modules.

Each dot in the chamber $w^{-1}C$ in the following picture represents a basis element of the $wt$ generalized weight space of $M$ or $N$. The dots that are connected by arcs represent basis vectors in the same module, $M$ or $N$.

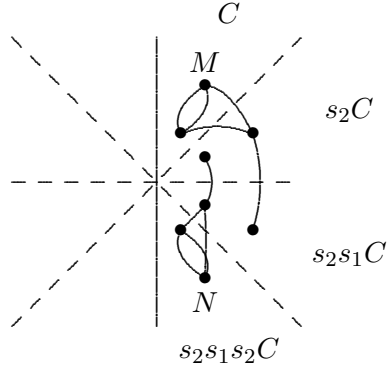
\[\beginpicture
\setcoordinatesystem units <0.32cm,0.32cm>
\setplotarea x from -15 to 15, y from -6 to 6
\linethickness=3pt
\put{$\bullet$} at -2 2
\put{$\bullet$} at -1 4
\put{$\bullet$} at -2 -2
\put{$\bullet$} at -1 -4
\put{$\bullet$} at -1 1
\put{$\bullet$} at 1 2
\put{$\bullet$} at -1 -1
\put{$\bullet$} at 1 -2
\put{$t_{1,q^{2}}$, $q$ generic} at -3 -10
\put{$C$} at 0 7
\put{$s_{2}C$} at 5 3
\put{$s_{2}s_{1}C$} at 5 -3
\put{$s_{2}s_{1}s_{2}C$} at 0 -7
\put{$M$} at -1 5
\put{$N$} at -1 -5
\setlinear
\plot -3 6 -3 -6 /
\plot -2 -2 -1 -1 /
\setdots
\setdashes
\plot -9 6 3 -6 /
\plot -9 -6 3 6 /
\plot -9 0 3 0 /
\setsolid
\setquadratic
\plot -2 2 -2 3 -1 4 /
\plot -2 2 -1 3 -1 4 /
\plot -2 -2 -2 -3 -1 -4 /
\plot -2 -2 -1 -3 -1 -4 /
\plot -1 -4 -.95 -2.5 -1 -1 /
\plot 1 2 1.3 0 1 -2 /
\plot -2 2 -0.5 2.3 1 2 /
\plot -1 -1 -0.7 0 -1 1 /
\plot -1 4 0 3.5 1 2 /
\endpicture \]

\begin{proposition} \label{t1q2B2fourth}
Assume $q^{2} = -1$ and let $M = \widetilde{H} \otimes_{\widetilde{H}_{\{2\}}} \mathbb{C}_{1,q^{2}} \quad \textrm{ and } \quad N = \widetilde{H} \otimes_{\widetilde{H}_{\{2\}}} \mathbb{C}_{1,q^{-2}}$. Then

(a) $M_{s_{1}s_{2}t}$ is a submodule of $M$, and the image of $M_{s_{2}t}$ is a submodule of $M/M_{s_{1}s_{2}t}$. The resulting 2-dimensional quotient of $M$ is irreducible. Also, $N_{s_{2}t}$ is a submodule of $N$ and the image of $N_{s_{1}s_{2}t}$ in $N/N_{s_{2}t}$ is a submodule of $N/N_{s_{2}t}$. The resulting 2-dimensional quotient of $N$ is irreducible, and

(b) Any composition factor of $M(t)$ is a composition factor of either $M$ or $N$.

\end{proposition}

(a) If $q^{2} = -1$, then $M$ has weight spaces $M_{t}^{\textrm{gen}}$, which is two dimensional, and $M_{s_{2}t}$ and $M_{s_{1}s_{2}t}$, both of which are 1-dimensional. Theorems $\ref{weightbasis}$ and $\ref{tauthm}$ show that $\tau_{1}$ is non-zero on $M_{s_{2}t}$, but zero on $M_{s_{1}s_{2}t}$, so that $M_{s_{1}s_{2}t}$ is a submodule of $M$. The resulting quotient must be reducible, but since $v_{t}$ generates all of $M$, the generalized $t$ weight space cannot be a submodule. This the $s_{i}t$ weight space is the submodule, and its quotient must be the 2-dimensional module constructed in Theorem \ref{2dims}, since it accounts for the entire $t$ weight space of $M(t)$. A similar argument shows the result for $N$.

(b) By counting dimensions of weight spaces, the remaining composition factor(s) of $M(t)$ must have weights $s_{2}t$ and $s_{1}s_{2}t$. If there were only one composition factor $L$ left, it would contain both weight spaces which would each have dimension 1, which is impossible by Theorem $\ref{2dims}$. Thus the remaining composition factors are more copies of the 1-dimensional modules. $\square$

\vspace{.2in}

\begin{proposition}
Assume $q^{2} \neq \pm 1$. Then $M_{s_{1}s_{2}t}$ is a submodule of $M$ and $M/M_{s_{1}s_{2}t}$ is irreducible. Similarly, $N_{s_{2}t}$ is a submodule of $N$ and $N/N_{s_{2}t}$ is irreducible.
\end{proposition}

\noindent \textit{Proof.} If $q^{4} \neq 1$, then by the same reasoning as in Proposition \ref{t1q2B2fourth}, $M_{s_{1}s_{2}t}$ must be a submodule of $M$. Similarly, $N_{s_{2}t}$ is a submodule of $N$. Then Lemma \ref{eq:lemma} shows that the resulting 3-dimensional quotients of $M$ and $N$ are irreducible. $\square$

If $q^{2} \neq 1$, the composition factors of $M$ and $N$ account for all 8 dimensions of $M(t)$. The following pictures show the composition factors of $M(t)$. Each dot in $w^{-1}C$ represents one basis element of $M(t)_{wt}^{\textrm{gen}}$, and the basis elements corresponding to connected dots are in the same composition factor.

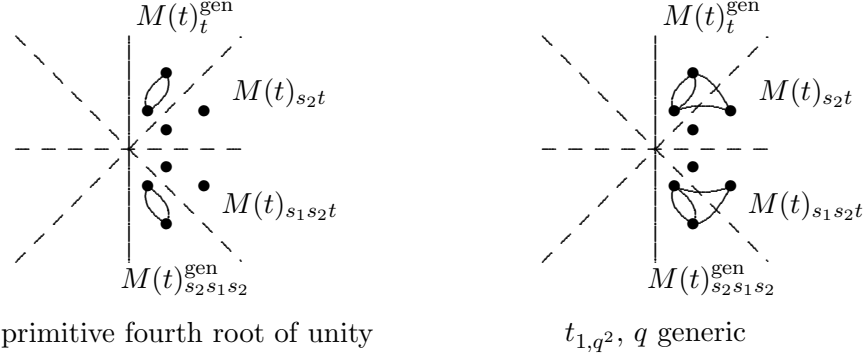
\[\beginpicture
\setcoordinatesystem units <0.25cm,0.25cm>
\setplotarea x from -25 to 25, y from -6 to 6
\linethickness=3pt
\put{$\bullet$} at -10 2
\put{$\bullet$} at -10 -2
\put{$\bullet$} at -9 4
\put{$\bullet$} at -9 -4
\put{$\bullet$} at -9 1
\put{$\bullet$} at -9 -1
\put{$\bullet$} at -7 2
\put{$\bullet$} at -7 -2
\put{$\bullet$} at 18 2
\put{$\bullet$} at 19 4
\put{$\bullet$} at 18 -2
\put{$\bullet$} at 19 -4
\put{$\bullet$} at 19 1
\put{$\bullet$} at 21 2
\put{$\bullet$} at 19 -1
\put{$\bullet$} at 21 -2
\put{$t_{1,q^{2}}$, $q$ a primitive fourth root of unity} at -11 -10
\put{$t_{1,q^{2}}$, $q$ generic} at 17 -10
\put{$M(t)_{t}^{\textrm{gen}}$} at -8 7
\put{$M(t)_{s_{2}t}$} at -3 3
\put{$M(t)_{s_{1}s_{2}t}$} at -3 -3
\put{$M(t)_{s_{2}s_{1}s_{2}}^{\textrm{gen}}$} at -8 -7
\put{$M(t)_{t}^{\textrm{gen}}$} at 20 7
\put{$M(t)_{s_{2}t}$} at 25 3
\put{$M(t)_{s_{1}s_{2}t}$} at 25 -3
\put{$M(t)_{s_{2}s_{1}s_{2}}^{\textrm{gen}}$} at 20 -7
\setlinear
\plot 17 6 17 -6 /
\plot -11 -6 -11 6 /
\setdots
\setdashes
\plot 11 6 23 -6 /
\plot -17 -6 -5 6 /
\plot 11 -6 23 6 /
\plot -17 0 -5 0 /
\plot -17 6 -5 -6 /
\plot 11 0 23 0 /
\setsolid
\setquadratic
\plot -10 2 -10 3 -9 4 /
\plot -10 2 -9 3 -9 4 /
\plot -10 -2 -10 -3 -9 -4 /
\plot -10 -2 -9 -3 -9 -4 /
\plot 18 2 18 3 19 4 /
\plot 18 2 19 3 19 4 /
\plot 18 -2 18 -3 19 -4 /
\plot 18 -2 19 -3 19 -4 /
\plot 18 -2 19.5 -2.3 21 -2 /
\plot 19 -4 20 -3.5 21 -2 /
\plot 18 2 19.5 2.3 21 2 /
\plot 19 4 20 3.5 21 2 /
\endpicture \]

$\underline{t|_{Q} = t_{\pm q,1}}$.

Let $\mathbb{C}_{\pm q,1}$ and $\mathbb{C}_{\pm q^{-1},1}$ be the 1-dimensional $\widetilde{H}_{\{2\}}$-modules spanned by $v_{s_{1}t}$ and $v_{s_{2}s_{1}t}$, respectively, and given by

\[ T_{2} v_{s_{1}t} = qv_{s_{1}t} \quad \textrm{ and } \quad X^{\lambda}v_{s_{1}t} = s_{1}t(X^{\lambda})v_{s_{1}t}, \quad \textrm{ and } \]

\[ T_{2} v_{s_{2}s_{1}t} = -q^{-1}v_{s_{2}s_{1}t} \quad \textrm { and } X^{\lambda}v_{s_{2}s_{1}t} = s_{2}s_{1}t(X^{\lambda})v_{s_{2}s_{1}t}.\] Then \[M = \widetilde{H} \otimes_{\widetilde{H}_{\{2\}}} \mathbb{C}_{\pm q,1} \quad \textrm{ and } \quad N = \widetilde{H} \otimes_{\widetilde{H}_{\{2\}}} \mathbb{C}_{\pm q^{-1}, 1}\] are 4-dimensional $\widetilde{H}$-modules.

If $t|_{Q} = t_{-q,1}$ and $q$ is a primitive sixth root of unity or if $t|_{Q} = t_{q,1}$ and $q$ is a primitive third root of unity, then $t|_{Q} = t_{q^{-2},1}$, which is in the same orbit as $t_{q^{2},1}$, and the irreducibles with central character $t$ have already been analyzed.

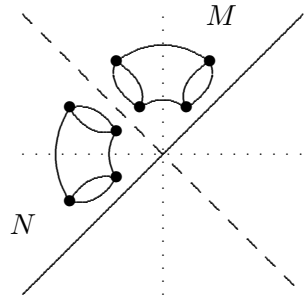
\[\beginpicture
\setcoordinatesystem units <0.31cm,0.31cm>
\setplotarea x from -15 to 15, y from -6 to 6
\linethickness=3pt
\put{$\bullet$} at 1 2
\put{$\bullet$} at 2 4
\put{$\bullet$} at -1 2
\put{$\bullet$} at -2 4
\put{$\bullet$} at -2 1
\put{$\bullet$} at -4 2
\put{$\bullet$} at -2 -1
\put{$\bullet$} at -4 -2
\put{$t_{\pm q ,1}$ $q^{2} \neq 1$,} at 0 -8
\put{(excluding $t_{-q,1}$ when $q$ is a primitive sixth root of unity, and} at 0 -9.8
\put{$t_{q,1}$ when $q$ is a primitive third root of unity.)} at 0 -11.6
\put{$M$} at 2.5 6
\put{$N$} at -6 -3
\setlinear
\plot -6 -6 6 6 /
\setdots
\plot 0 -6 0 6 /
\plot -6 0 6 0 /
\setdashes
\plot -6 6 6 -6 /
\setsolid
\setquadratic
\plot -2 1 -2.3 0 -2 -1 /
\plot -4 2 -4.6 0 -4 -2 /
\plot 1 2 0 2.35 -1 2 /
\plot 2 4 0 4.7 -2 4 /
\plot 1 2 1 3 2 4 /
\plot 1 2 2 3 2 4 /
\plot -1 2 -1 3 -2 4 /
\plot -1 2 -2 3 -2 4 /
\plot -2 1 -3.1 1.1 -4 2 /
\plot -2 1 -2.9 1.9 -4 2 /
\plot -2 -1 -3.1 -1.1 -4 -2 /
\plot -2 -1 -2.9 -1.9 -4 -2 /
\endpicture \]

\begin{proposition}
Let $M = \widetilde{H} \otimes_{\widetilde{H}_{\{2\}}} \mathbb{C}_{\pm q,1} \quad \textrm{ and } \quad N = \widetilde{H} \otimes_{\widetilde{H}_{\{2\}}} \mathbb{C}_{\pm q^{-1}, 1}$. Unless $t|_{Q} = t_{-q,1}$ and $q$ is a primitive sixth root of unity or $t|_{Q} = t_{q,1}$ and $q$ is a primitive third root of unity, $M$ and $N$ are irreducible.

\end{proposition}

The assumption gives that $P(t) = \{2\alpha_{1} + \alpha_{2}\}$. Then Lemma \ref{eq:lemma} shows that both $M$ and $N$ are irreducible. $\square$

Since they have different weight spaces and are thus not isomorphic, $M$ and $N$ are the only two irreducibles with central character $t$.

\vspace{.1in}

\noindent \textbf{Summary.} The following tables summarize the classification. It should be noted that for any value of $q$ with $q^{2}$ not a root of unity of order 4 or less, the representation theory of $\widetilde{H}$ can be described in terms of $q$ only. If $q^{2}$ is a primitive root of unity of order 4 or less, then the representation theory of $\widetilde{H}$ does not fit that same description. This fact can be seen through a number of different lenses. It is a reflection of the fact that the sets $P(t)$ and $Z(t)$ for all possible central characters $t$ can be described solely in terms of $q$, as long as $q^{2}$ is not one of these small roots of unity. In the local region pictures, this is reflected in the fact that the hyperplanes $H_{\alpha}$ and $H_{\alpha \pm \delta}$ are distinct \textit{unless} $q^{2}$ is a root of unity of order 4 or less. When these hyperplanes coincide, the sets $P(t)$ and $Z(t)$ change for characters on those hyperplanes. Note also that in this table, the entries for $t_{\pm q,1}$ assume that $\pm q \neq q^{-2}$, as described above.

\begin{table}[htpb]
\centering
$\begin{array}{| c | c | c | c | c | c |}
\hline \, & \multicolumn{5}{|c|}{\textrm{Dims. of Irreds.}} \\
\hline t & q^{8} \neq 1, q^{6} \neq 1 & q^{8} = 1, q^{4} \neq 1 & q^{6} = 1 & q^{2} = -1 & q = -1 \\
\hline t_{1,1} & 8 & 8 & 8 & 8 & 1,1,1,1,2 \\
\hline t_{-1,1} & 8 & 8 & 8 & 4 & 2,2,2,2 \\
\hline t_{1,z} & 8 & 8 & 8 & 8 & 4,4 \\
\hline t_{1,q^{2}} & 1,1,3,3 & 1,1,3,3 & 1,1,3,3 & 1,1,2,2 & \textrm{N/A} \\
\hline t_{q^{2},1} & 4,4 & 4,4 & 1,1,3,3 & \textrm{N/A} & \textrm{N/A} \\
\hline t_{q,1} & 4,4 & 4,4 & 4,4 & 4,4 & \textrm{N/A} \\
\hline t_{-q,1} & 4,4 & 4,4 & 4,4 & \textrm{N/A} & \textrm{N/A} \\
\hline t_{z,1} & 8 & 8 & 8 & 8 & 4,4 \\
\hline t_{q^{2},q^{2}} & 1,1,3,3 & 1,1,1,1,2,2 & \textrm{N/A} & \textrm{N/A} & \textrm{N/A} \\
\hline t_{q^{2},z} & 4,4 & 4,4 & 4,4 & 4,4 & \textrm{N/A} \\
\hline t_{-1,q^{2}} & 2,2,2,2 & \textrm{N/A} & 2,2,2,2 & \textrm{N/A} & \textrm{N/A} \\
\hline t_{z,q^{2}} & 4,4 & 4,4 & 4,4 & 4,4 & \textrm{N/A} \\
\hline t_{z,w} & 8 & 8 & 8 & 8 & 8 \\
\hline \end{array}$
\caption{Table of possible central characters in Type $C_{2}$, with varying values of $q$.}
\end{table}

\section{Type $G_{2}$}

The type $G_{2}$ root system is \[ R = \{ \pm \alpha_{1}, \pm \alpha_{2}, \pm (\alpha_{1} + \alpha_{2}), \pm (2\alpha_{1} + \alpha_{2}), \pm (3\alpha_{1} + \alpha_{2}), \pm(3\alpha_{1} + 2\alpha_{2})\},\] with $\langle \alpha_{1}, \alpha_{2}^{\vee} \rangle = -1$ and $\langle \alpha_{2}, \alpha_{1}^{\vee} \rangle = -3$. Then the Weyl group is \[ W_{0} = \langle s_{1}, s_{2} \, | \, s_{1}^{2} = s_{2}^{2} = 1, s_{1}s_{2}s_{1}s_{2}s_{1}s_{2} = s_{2}s_{1}s_{2}s_{1}s_{2}s_{1} \rangle,\] isomorphic to the dihedral group of order 12. The simple roots are $\alpha_{1}$ and $\alpha_{2}$, and $\alpha_{1}, \alpha_{1} + \alpha_{2}, 3\alpha_{1} + \alpha_{2}$ will be referred to as \textit{short} roots, while $\alpha_{2}, 2\alpha_{1} + \alpha_{2}$, and $3\alpha_{1} + 2\alpha_{2}$ will be referred to as \textit{long} roots.

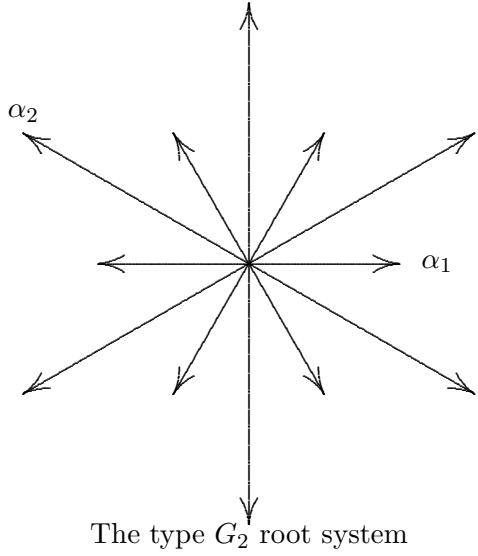
\[ \beginpicture
\setcoordinatesystem units <0.25cm,0.25cm>
\setplotarea x from -25 to 25, y from -13 to 13
\linethickness=5pt
\setlinear
\thicklines
\arrow <10 pt> [.2,.67] from 0 0 to 0 13.856
\arrow <10 pt> [.2,.67] from 0 0 to 0 -13.856
\arrow <10 pt> [.2,.67] from 0 0 to 12 6.928
\arrow <10 pt> [.2,.67] from 0 0 to 12 -6.928
\arrow <10 pt> [.2,.67] from 0 0 to -12 6.928
\arrow <10 pt> [.2,.67] from 0 0 to -12 -6.928
\arrow <10 pt> [.2,.67] from 0 0 to 4 6.928
\arrow <10 pt> [.2,.67] from 0 0 to 4 -6.928
\arrow <10 pt> [.2,.67] from 0 0 to -4 6.928
\arrow <10 pt> [.2,.67] from 0 0 to -4 -6.928
\arrow <10 pt> [.2,.67] from 0 0 to 8 0
\arrow <10 pt> [.2,.67] from 0 0 to -8 0
\put{The type $G_{2}$ root system} at 0 -14.5
\put{$\alpha_{1}$} at 10 0
\put{$\alpha_{2}$} at -12 8
\endpicture\]

The fundamental weights satisfy \[ \begin{matrix} \omega_{1} = 2\alpha_{1} + \alpha_{2} & \quad & \alpha_{1} = 2\omega_{1} - \omega_{2} \\ \omega_{2} = 3\alpha_{1} + 2\alpha_{2} & \quad & \alpha_{2} = 2\omega_{2} - 3\omega_{1}. \end{matrix}\] Let \[P = \mathbb{Z}\textrm{-span}\{\omega_{1},\omega_{2}\}.\] This is the same lattice spanned by $\alpha_{1}$ and $\alpha_{2}$. Then $W_{0}$ acts on $X$ by \begin{eqnarray*} s_{1} \cdot X^{\omega_{1}} & = & X^{\omega_{2}-\omega_{1}}, \\ s_{1} \cdot X^{\omega_{2}} & = & X^{\omega_{2}}, \\ s_{2} \cdot X^{\omega_{1}} & = & X^{\omega_{1}}, \quad \textrm{ and } \\  s_{2} \cdot X^{\omega_{2}} & = & X^{3\omega_{1}-\omega_{2}}.\end{eqnarray*}
The affine Hecke algebra $\widetilde{H}$ of type $G_{2}$ is defined as in section $\ref{subsectHecke}$.

Let \[\mathbb{C}[X] = \{ X^{\lambda} \, | \, \lambda \in P\}, \] a subalgebra of $\widetilde{H}$, and let \[ T = \textrm{Hom}_{\mathbb{C}\textrm{-alg}}(\mathbb{C}[X],\mathbb{C}).\] Define \[ t_{z,w} : T \rightarrow \mathbb{C} \textrm{ by } t_{z,w}(X^{\alpha_{1}}) = z \textrm{ and } t_{z,w}(X^{\alpha_{2}}) = w.\]

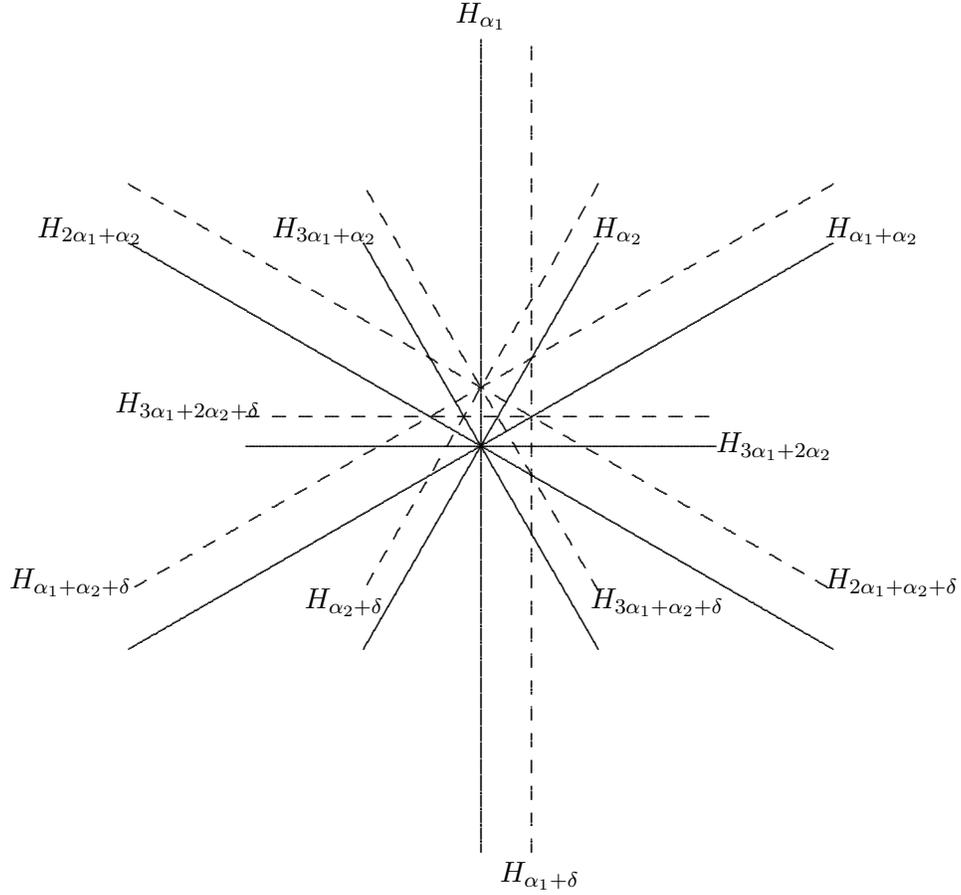
\[ \beginpicture
\setcoordinatesystem units <0.26cm,0.26cm>
\setplotarea x from -25 to 25, y from -15 to 15
\linethickness=5pt
\plot 0 -20.784 0 20.784 /
\plot 6 10.392 -6 -10.392 /
\plot 6 -10.392 -6 10.392 /
\plot -18 10.392 18 -10.392 /
\plot 18 10.392 -18 -10.392 /
\plot -12 0 12 0 /
\put{$H_{\alpha_{1}}$} at 0 22
\put{$H_{\alpha_{2}}$} at 7 11
\put{$H_{\alpha_{1} + \alpha_{2}}$} at 20 11
\put{$H_{2\alpha_{1} + \alpha_{2}}$} at -20 11
\put{$H_{3\alpha_{1} + \alpha_{2}}$} at -8 11
\put{$H_{3\alpha_{1} + 2\alpha_{2}}$} at 15 0
\setdashes
\plot 2.598 -20.784 2.598 20.784 /
\plot 6 13.392 -6 -7.392 /
\plot 6 -7.392 -6 13.392 /
\plot -18 13.392 18 -7.392 /
\plot 18 13.392 -18 -7.392 /
\plot -12 1.5 12 1.5 /
\put{$H_{\alpha_{1}+\delta}$} at 3 -22
\put{$H_{\alpha_{2}+\delta}$} at -7 -8
\put{$H_{\alpha_{1} + \alpha_{2}+\delta}$} at -21 -7
\put{$H_{2\alpha_{1} + \alpha_{2}+\delta}$} at 21 -7
\put{$H_{3\alpha_{1} + \alpha_{2}+\delta}$} at 9 -8
\put{$H_{3\alpha_{1} + 2\alpha_{2}+\delta}$} at -15 2
\endpicture \]

The structure of the modules with weight $t$ depends virtually exclusively on $P(t) = \{ \alpha \in R^{+} \, | \, t(X^{\alpha})=q^{\pm 2} \}$ and $Z(t) = \{ \alpha \in R^{+} \, | \, t(X^{\alpha}) = 1 \}$. For a generic weight $t$, $P(t)$ and $Z(t)$ are empty, so we examine only the non-generic orbits.

\begin{theorem} \label{G2charsgen} If $q^{2}$ is not a primitive $\ell$th root of unity for $\ell \leq 6$ and $Z(t) \cup P(t) \neq \emptyset$, then $t$ is in the same $W_{0}$-orbits as one of the following weights.

\[ t_{1,1}, t_{1,-1}, t_{1^{1/3},1}, t_{1,q^{2}}, t_{1,\pm q}, t_{q^{2},1}, t_{\pm q, 1}, t_{q^{2/3},1}, t_{q^{2},-q^{-2}}, t_{1^{1/3},q^{2}}, t_{q^{2},q^{2}}, \] \[ \{ t_{1,z} \, | \, z \in \mathbb{C}^{\times}, z \neq \pm 1, q^{\pm 2}, \pm q^{\pm 1} \}, \{ t_{z,1} \, | \, z \in \mathbb{C}^{\times}, z \neq \pm 1, 1^{1/3}, q^{\pm 2}, \pm q^{\pm 1}, q^{\pm 2/3} \},\] \[ \left\{ t_{q^{2},z} \, | \, z \in \mathbb{C}^{\times}, \{ 1,q^{2},q^{-2}\} \cap \{z,q^{2}z,q^{4}z,q^{6}z,q^{6}z^{2}\} = \emptyset \right\}, \textrm{ or } \] \[ \left\{ t_{z,q^{2}} \, | \, z \in \mathbb{C}^{\times}, \{ 1,q^{2},q^{-2}\} \cap \{z,q^{2}z,q^{2}z^{2},q^{2}z^{3},q^{4}z^{3}\} = \emptyset \right\}.\]\end{theorem}

\noindent \textit{Proof.} In general, the third roots of unity in this theorem are assumed to be primitive, so that there are two different weights that we call $t_{1^{1/3},1}$ and $t_{1^{1/3},q^{2}}$. Similarly, $t_{q^{2/3},1}$ typically refers to one of three different characters, corresponding to the three third roots of $q^{2}$. Exceptions to this principle are noted as they arise. We refer to $\alpha_{1}$, $\alpha_{1}+\alpha_{2}$, and $2\alpha_{1}+\alpha_{2}$ as ``short'' roots. The other roots are referred to as ``long'' roots.

\vspace{.1in}

\underline{Case 1: $|Z(t)| \geq 2$}

If $Z(t)$ contains at least two roots, and one of them is short, then $Z(wt)$ contains $\alpha_{1}$ for some $w \in W_{0}$. If $Z(wt)$ also contains any of $\alpha_{2}, \alpha_{1}+\alpha_{2}, 2\alpha_{1}+\alpha_{2},$ or $3\alpha_{1}+\alpha_{2}$, then it contains both simple roots and thus $wt = t_{1,1}$. It is also possible that $Z(wt) = \{ \alpha_{1}, 3\alpha_{1}+2\alpha_{2}\}$, in which case $wt(X^{\alpha_{1}+\alpha_{2}}) = -1$, and $wt(X^{\alpha_{2}}) = -1$, so that $wt = t_{1,-1}$.

If $Z(t)$ contains no short roots, it contains two of $\alpha_{2}$, $3\alpha_{1}+\alpha_{2}$, and $3\alpha_{1}+2\alpha_{2}$. But then it must also contain the third, and $Z(t) = \{ \alpha_{2}, 3\alpha_{1}+\alpha_{2}, 3\alpha_{1}+2\alpha_{2}\}$. In this case, $wt(X^{\alpha_{2}}) = 1$, but $wt(X^{\alpha_{1}})$ is a third root of unity, so that $wt = t_{1^{1/3},1}$.

\vspace{.1in}

\underline{Case 2: $|Z(t)| = 1$}

If $Z(t)$ has exactly one root, then there is some $w \in W_{0}$ with $Z(wt) = \{\alpha_{1}\}$ or $Z(wt) = \{ \alpha_{2}\}$.

If $Z(wt) = \{\alpha_{1}\}$, then $P(t)$ either contains all of $\alpha_{2}, \alpha_{1}+\alpha_{2}, 2\alpha_{1}+\alpha_{2},$ and $3\alpha_{1}+\alpha_{2}$, or it contains none of them. If it contains all of them, $wt(X^{\alpha_{2}}) = q^{\pm 2}$, and $t$ is in the same $W_{0}$-orbit as $t_{1,q^{2}}$. If $3\alpha_{1}+2\alpha_{2} \in P(wt)$, then $wt(X^{\alpha_{2}}) = \pm q^{\pm 1}$, and $t$ is in the same orbit as $t_{1,\pm q}$. Otherwise, $wt = t_{1,z}$ for some $z$ besides $\pm q^{\pm 1}$ and $q^{\pm 2}$. Also then, $z \neq \pm 1$ by assumption on $Z(t)$.

If $Z(wt) = \{ \alpha_{2}\}$, then any two roots that differ by a multiple of $\alpha_{2}$ are either both in $P(wt)$ or both not in $P(wt)$. By applying $w_{0}$ if necessary, we can assume that $wt(X^{\alpha}) = q^{2}$ for the $\alpha$ that are in $P(wt)$. If $\alpha_{1} \in P(wt)$, then $wt(X^{\alpha_{1}}) = q^{2}$, and $wt = t_{q^{2},1}$. If $2\alpha_{1}+\alpha_{2} \in P(wt)$, then $wt = t_{\pm q,1}$. If $3\alpha_{1}+\alpha_{2} \in P(wt)$, then $wt(X^{\alpha_{1}})$ is a third root of $q^{2}$ and $wt = t_{1,q^{2/3}}$. Otherwise, $wt = t_{z,1}$ for some $z$ so that none of $z,z^{2},z^{3}$ is equal to $q^{\pm 2}$ or $1$. That is, $z \neq \pm 1, 1^{1/3}, q^{\pm 2}, \pm q^{\pm 1}, q^{\pm 2/3}$.

\vspace{.1in}

\underline{Case 3: $|Z(t)| = \emptyset$}

If $Z(t)$ is empty but $P(t)$ contains a short root, then $\alpha_{1} \in P(wt)$ for some $w \in W_{0}$. If $P(wt)$ contains another short root, then we can apply $s_{1}$ if necessary so that $P(t)$ contains $\alpha_{1}$ and $\alpha_{1}+\alpha_{2}$. Then either $wt(X^{\alpha_{1}}) = wt(X^{\alpha_{1}+\alpha_{2}})$ so that $wt(X^{\alpha_{2}}) =1$, or $wt(X^{\alpha_{1}})$ and $wt(X^{\alpha_{1}+\alpha_{2}})$ are $q^{2}$ and $q^{-2}$ in some order, so that $wt(X^{2\alpha_{1}+\alpha_{2}}) = 1$. Thus $P(wt)$ contains at most one short root. If $P(wt)$ also contains a long root, then applying $s_{1}$ if necessary, we can assume $P(wt)$ contains either $\alpha_{2}$ or $3\alpha_{1}+2\alpha_{2}$. If $P(wt)$ contains $\alpha_{1}$ and $\alpha_{2}$, then we can apply $w_{0}$ to assume that $wt(X{\alpha_{1}}) = q^{2}$. If $wt(X^{\alpha_{2}}) = q^{-2}$, then $\alpha_{1}+\alpha_{2} \in Z(wt)$. Then $wt(X^{\alpha_{2}}) = q^{2}$ and $wt = t_{q^{2},q^{2}}$. If $P(wt)$ contains $ \alpha_{1}$ and $3\alpha_{1}+2\alpha_{2}$, then since $\alpha_{1}$ is perpendicular to $3\alpha_{1}+2\alpha_{2}$, we can apply $s_{1}$ and/or $s_{3\alpha_{1}+2\alpha_{2}}$ to assume $wt(X^{\alpha_{1}}) = q^{2} = wt(X^{3\alpha_{1}+2\alpha_{2}})$. Hence $wt(X^{2\alpha_{2}}) = q^{-4}$ and by assumption, $wt = t_{q^{2},-q^{-2}}$. If $P(wt) = \{ \alpha_{1}\}$, then $wt = t_{q^{2},z}$ does not take the value $1$ or $q^{\pm 2}$ on any other positive root. Then $\{ 1,q^{2},q^{-2}\} \cap \{z,q^{2}z,q^{4}z,q^{6}z,q^{6}z^{2}\} = \emptyset$

If $P(t)$ contains no short roots, but at least two long roots, then $wt(X^{\alpha_{2}})  = q^{2} = wt(X^{3\alpha_{1}+\alpha_{2}})$ for some $w \in W_{0}$. (If $wt(X^{\alpha_{2}}) = q^{-2}$, then $wt(X^{3\alpha_{1}+2\alpha_{2}}) = 1$, a contradiction.) Hence $wt(X^{\alpha_{1}})$ is a primitive third root of unity and $wt = t_{1^{1/3},q^{2}}$. If $P(t)$ contains exactly one long root, then $wt = t_{z,q^{2}}$ for some $z \in \mathbb{C}^{\times}$ so that $wt$ does not take the value $1$ or $q^{\pm 2}$ on any other positive root. Thus $\{ 1,q^{2},q^{-2}\} \cap \{z,q^{2}z,q^{2}z^{2},q^{2}z^{3},q^{4}z^{3}\} = \emptyset$. $\square$

\noindent \textbf{Remark.} There are some redundancies in this list for specific values of $q$. If $q^{2}$ is a primitive fifth root of unity, then $q$ and $-q$ are equal to $q^{-4}$ and $-q^{-4}$ in some order depending on whether $q^{5} = 1$ or $-1$. Then $t_{q^{2},q^{2}}$ is in the same orbit as $t_{q^{-4},1}$, which is equal to either $t_{q,1}$ or $t_{-q,1}$.

If $q^{2}$ is a primitive fourth root of unity, then one note is necessary on the weight $t_{q^{2/3},1}$. Since $q^{-2}$ is a third root of unity, we take $q^{2/3}$ to mean a different third root of $q^{2}$ so that $t_{q^{2/3},1}$ and $t_{q^{2},1}$ are in different orbits. In addition, $t_{q^{2},-q^{-2}} = t_{q^{2},q^{2}}$, which is in the same orbit as $t_{q^{2},1}$.

If $q^{2}$ is a primitive third root of unity, then $1^{1/3} = q^{2}, q^{-2}$, or 1. Then $t_{1^{1/3},1}$ is in the same orbit as $t_{q^{2},1}$ or $t_{1,1}$. Also, $t_{1^{1/3},q^{2}}$ is in the same orbit as $t_{q^{2},q^{2}}$, which is in turn in the same orbit as $t_{1,q^{2}}$. In addition $q$ and $-q$ are equal to $q^{-2}$ and $-q^{-2}$ in some order depending on whether $q^{3}$ is $1$ or $-1$. Then $t_{1,q^{2}}$ is in the same orbit as either $t_{1,q}$ or $t_{1,-q}$, and $t_{q^{2},1}$ is in the same orbit as either $t_{q,1}$ or $t_{-q,1}$.

If $q^{2} = -1$, then $t_{q^{2},-q^{-2}} = t_{q^{2},1} = t_{-1,1}$, which is in the same orbit as $t_{q^{2},q^{2}}$, while $t_{1,q^{2}} = t_{1,-1}$. In fact, $t_{-1,1} = s_{1}s_{2}t_{1,-1}$. Also, since $q = -q^{-1}$, the weights $t_{1, \pm q}$ are in the same orbit as each other, as are the weights $t_{\pm q, 1}$. Finally, $t_{1^{1/3},q^{2}}$ is in the same orbit as $t_{q^{2/3},1}$.

Finally, if $q = -1$, then $t_{1,1} = t_{q^{2},1} = t_{1,q^{2}} = t_{q^{2},q^{2}} = t_{1,-q} = t_{-q,1}$. Also, $t_{q,1} = t_{-1,1},$ which is in the same orbit as $t_{1,-1} = t_{1,q} = t_{q^{2},-q^{-2}}$. Finally, $t_{1^{1/3},1} = t_{q^{2/3},1} = t_{1^{1/3},q^{2}}$, while $t_{1,z} = t_{q^{2},z}$ and $t_{z,1} = t_{z,q^{2}}$.

\vspace{.1in}

\noindent \textbf{Analysis of the characters.}

\begin{theorem}\label{G2ones} The 1-dimensional representations of $\widetilde{H}$ are:

\renewcommand{\baselinestretch}{1} \normalsize $\begin{array}{cccc}  L_{q,q}: \widetilde{H} \rightarrow \mathbb{C} \quad & L_{q,-q^{-1}}: \widetilde{H} \rightarrow \mathbb{C} \quad & L_{q,-q^{-1}}: \widetilde{H} \rightarrow \mathbb{C} \quad & L_{-q^{-1},-q^{-1}}: \widetilde{H} \rightarrow \mathbb{C} \quad \\ T_{1} \mapsto q & T_{1} \mapsto q & T_{1} \mapsto -q^{-1} & T_{1} \mapsto -q^{-1} \\ T_{2} \mapsto q & T_{2} \mapsto -q^{-1} & T_{2} \mapsto q & T_{2} \mapsto -q^{-1} \\ X^{\omega_{1}} \mapsto q^{6} & X^{\omega_{1}} \mapsto q^{2} & X^{\omega_{1}} \mapsto q^{-2} & X^{\omega_{1}} \mapsto q^{-6} \\ X^{\omega_{2}} \mapsto q^{10} & X^{\omega_{2}} \mapsto q^{2} & X^{\omega_{2}} \mapsto q^{-2} & X^{\omega_{2}} \mapsto q^{-10} \end{array}$   \normalsize

\end{theorem}

\noindent \textit{Proof.} This is achieved in the same way as Theorem \ref{eq:B2ones}.\\ $\square$

Remark: If $q^{2} = -1$, then $q = -q^{-1}$ and all four 1-dimensional representations are isomorphic. Let $t \in T$. The principal series module is \[ M(t) = \textrm{Ind}_{\mathbb{C}[X]}^{\widetilde{H}} \mathbb{C}_{t} = \widetilde{H} \otimes_{\mathbb{C}[X]} \mathbb{C}_{t},\] where $\mathbb{C}_{t}$ is the one-dimensional $\mathbb{C}[X]$-module given by
\[ \mathbb{C}_{t} = \textrm{span}\{v_{t}\} \quad \textrm{ and } \quad X^{\lambda}v_{t} = t(X^{\lambda})v_{t}.\]
By \eqref{PrinQuot}c, every irreducible $\widetilde{H}$ module is a quotient of some principal series module $M(t)$. Thus, finding all the composition factors of $M(t)$ for all central characters $t$ suffices to find all the irreducible $\widetilde{H}$-modules.

\underline{Case 1:} $P(t) = \emptyset$.

If $P(t) = \emptyset$ then by Kato's criterion, Theorem \ref{eq:Kato}, $M(t)$ is irreducible and thus is the only irreducible module with central character $t$.

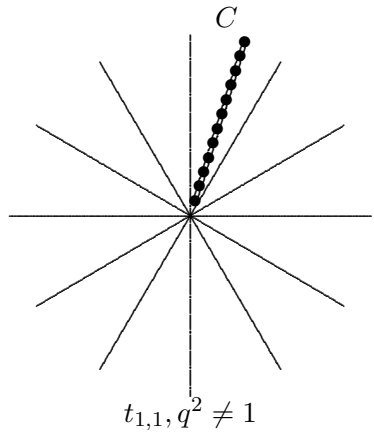
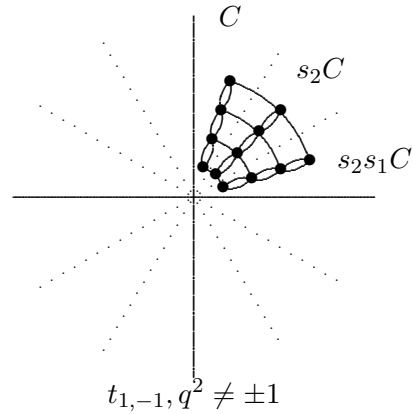
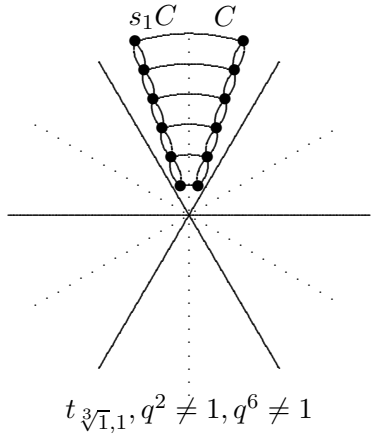
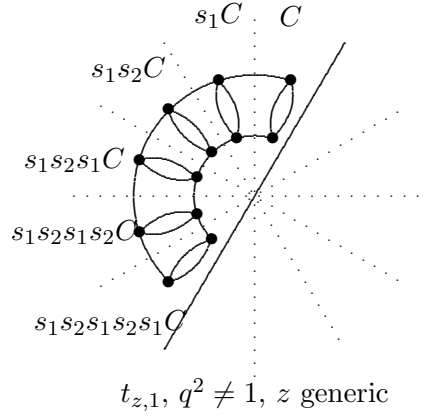
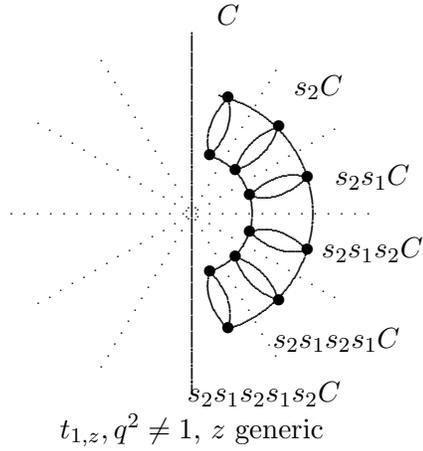
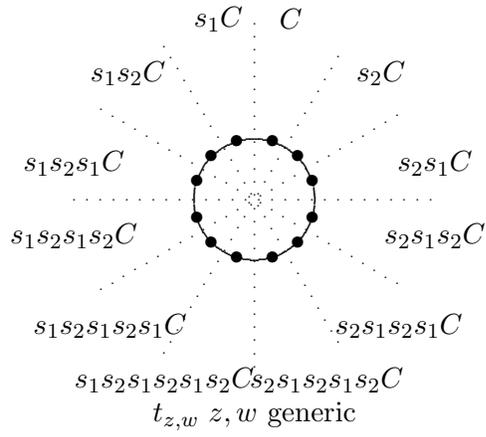
\[\beginpicture
\setcoordinatesystem units <0.24cm,0.24cm>
\setplotarea x from -25 to 25, y from -15 to 12
\linethickness=3pt
\plot -20 -8.5 -10 8.5 / 
\plot -20 8.5 -10 -8.5 / 
\plot -25 0 -5 0 / 
\plot -23.5 -5 -6.5 5 / 
\plot -23.5 5 -6.5 -5 / 
\plot -15 -10 -15 10 / 
\plot 5 0 25 0 / 
\plot 15 -10 15 10 / 
\setdots
\plot 10 -8.5 20 8.5 / 
\plot 10 8.5 20 -8.5 / 
\plot 6.5 -5 23.5 5 / 
\plot 6.5 5 23.5 -5 / 
\put{$\bullet$} at -14.75 .8
\put{$\bullet$} at -12 9.6
\put{$\bullet$} at -14.5 1.6
\put{$\bullet$} at -14.25 2.4
\put{$\bullet$} at -14 3.2
\put{$\bullet$} at -13.75 4
\put{$\bullet$} at -13.5 4.8
\put{$\bullet$} at -13.25 5.6
\put{$\bullet$} at -13 6.4
\put{$\bullet$} at -12.75 7.2
\put{$\bullet$} at -12.5 8
\put{$\bullet$} at -12.25 8.8
\put{$\bullet$} at 15.5 1.6
\put{$\bullet$} at 16 3.2
\put{$\bullet$} at 16.5 4.8
\put{$\bullet$} at 17 6.4
\put{$\bullet$} at 16.6 0.5
\put{$\bullet$} at 18.2 1
\put{$\bullet$} at 19.8 1.5
\put{$\bullet$} at 21.4 2
\put{$\bullet$} at 16.2 1.2
\put{$\bullet$} at 17.4 2.4
\put{$\bullet$} at 18.6 3.6
\put{$\bullet$} at 19.8 4.8
\put{$t_{1,1}, q^{2} \neq 1$} at -15 -11
\put{$t_{1,-1}, q^{2} \neq \pm 1$} at 15 -11
\setsolid
\circulararc 55 degrees from 16.6 0.5 center at 15 0
\circulararc 55 degrees from 18.2 1 center at 15 0
\circulararc 55 degrees from 19.8 1.5 center at 15 0
\circulararc 55 degrees from 21.4 2 center at 15 0
\setquadratic
\plot 16.6 0.5 17.4 0.5 18.2 1 /
\plot 16.6 0.5 17.4 1 18.2 1 /
\plot 19.8 1.5 19 1 18.2 1 /
\plot 18.2 1 19 1.5 19.8 1.5 /
\plot 19.8 1.5 20.6 2 21.4 2 /
\plot 19.8 1.5 20.6 1.5 21.4 2 /
\plot 16.2 1.2 16.9 1.6 17.4 2.4 /
\plot 16.2 1.2 16.7 2 17.4 2.4 /
\plot 18.6 3.6 18.1 2.8 17.4 2.4 /
\plot 17.4 2.4 17.9 3.2 18.6 3.6 /
\plot 18.6 3.6 19.3 4 19.8 4.8 /
\plot 18.6 3.6 19.1 4.4 19.8 4.8 /
\plot 15.5 1.6 15.5 2.4 16 3.2 /
\plot 15.5 1.6 16 2.4 16 3.2 /
\plot 16 3.2 16 4 16.5 4.8 /
\plot 16 3.2 16.5 4 16.5 4.8 /
\plot 16.5 4.8 17 5.6 17 6.4 /
\plot 16.5 4.8 16.5 5.6 17 6.4 /
\plot -14.75 .8 -14.75 1.2 -14.5 1.6 /
\plot -14.5 1.6 -14.5 2 -14.25 2.4 /
\plot -14.25 2.4 -14.25 2.8 -14 3.2 /
\plot -14 3.2 -14 3.6 -13.75 4 /
\plot -13.75 4 -13.75 4.4 -13.5 4.8 /
\plot -13.5 4.8 -13.5 5.2 -13.25 5.6 /
\plot -13.25 5.6 -13.25 6 -13 6.4 /
\plot -13 6.4 -13 6.8 -12.75 7.2 /
\plot -12.75 7.2 -12.75 7.6 -12.5 8 /
\plot -12.5 8 -12.5 8.4 -12.25 8.8 /
\plot -12.25 8.8 -12.25 9.2 -12 9.6 /
\plot -14.75 .8 -14.5 1.2 -14.5 1.6 /
\plot -14.5 1.6 -14.25 2 -14.25 2.4 /
\plot -14.25 2.4 -14 2.8 -14 3.2 /
\plot -14 3.2 -13.75 3.6 -13.75 4 /
\plot -13.75 4 -13.5 4.4 -13.5 4.8 /
\plot -13.5 4.8 -13.25 5.2 -13.25 5.6 /
\plot -13.25 5.6 -13 6 -13 6.4 /
\plot -13 6.4 -12.75 6.8 -12.75 7.2 /
\plot -12.75 7.2 -12.5 7.6 -12.5 8 /
\plot -12.5 8 -12.25 8.4 -12.25 8.8 /
\plot -12.25 8.8 -12 9.2 -12 9.6 /
\put{$C$} at -13 11
\put{$C$} at 17 10
\put{$s_{2}C$} at 22 7
\put{$s_{2}s_{1}C$} at 25 2
\endpicture\]

\[\beginpicture
\setcoordinatesystem units <0.24cm,0.24cm>
\setplotarea x from -25 to 25, y from -15 to 12
\linethickness=3pt
\plot -20 -8.5 -10 8.5 / 
\plot -20 8.5 -10 -8.5 / 
\plot -25 0 -5 0 / 
\plot 10 -8.5 20 8.5 / 
\setdots
\plot -23.5 -5 -6.5 5 / 
\plot -23.5 5 -6.5 -5 / 
\plot -15 -10 -15 10 / 
\plot 5 0 25 0 / 
\plot 15 -10 15 10 / 
\plot 10 8.5 20 -8.5 / 
\plot 6.5 -5 23.5 5 / 
\plot 6.5 5 23.5 -5 / 
\put{$\bullet$} at -15.5 1.6
\put{$\bullet$} at -12 9.6
\put{$\bullet$} at -14.5 1.6
\put{$\bullet$} at -16 3.2
\put{$\bullet$} at -14 3.2
\put{$\bullet$} at -16.5 4.8
\put{$\bullet$} at -13.5 4.8
\put{$\bullet$} at -17 6.4
\put{$\bullet$} at -13 6.4
\put{$\bullet$} at -17.5 8
\put{$\bullet$} at -12.5 8
\put{$\bullet$} at -18 9.6
\put{$\bullet$} at 14 3.2
\put{$\bullet$} at 16 3.2
\put{$\bullet$} at 13 6.4
\put{$\bullet$} at 17 6.4
\put{$\bullet$} at 8.6 -2
\put{$\bullet$} at 11.8 1
\put{$\bullet$} at 11.8 -1
\put{$\bullet$} at 8.6 2
\put{$\bullet$} at 12.6 -2.4
\put{$\bullet$} at 12.6 2.4
\put{$\bullet$} at 10.2 -4.8
\put{$\bullet$} at 10.2 4.8
\setdashes
\put{$t_{\sqrt[3]{1},1}, q^{2} \neq 1, q^{6} \neq 1$} at -15 -11
\put{$t_{z,1}$, $q^{2} \neq 1$, $z$ generic} at 15 -11
\setquadratic
\setsolid
\plot -15.5 1.6 -15.5 2.4 -16 3.2 /
\plot -15.5 1.6 -16 2.4 -16 3.2 /
\plot -16 3.2 -16 4 -16.5 4.8 /
\plot -16 3.2 -16.5 4 -16.5 4.8 /
\plot -16.5 4.8 -17 5.6 -17 6.4 /
\plot -16.5 4.8 -16.5 5.6 -17 6.4 /
\plot -17 6.4 -17.5 7.2 -17.5 8 /
\plot -17 6.4 -17 7.2 -17.5 8 /
\plot -17.5 8 -18 8.8 -18 9.6 /
\plot -17.5 8 -17.5 8.8 -18 9.6 /
\plot -14.5 1.6 -14.5 2.4 -14 3.2 /
\plot -14.5 1.6 -14 2.4 -14 3.2 /
\plot -14 3.2 -14 4 -13.5 4.8 /
\plot -14 3.2 -13.5 4 -13.5 4.8 /
\plot -13.5 4.8 -13 5.6 -13 6.4 /
\plot -13.5 4.8 -13.5 5.6 -13 6.4 /
\plot -13 6.4 -12.5 7.2 -12.5 8 /
\plot -13 6.4 -13 7.2 -12.5 8 /
\plot -12.5 8 -12 8.8 -12 9.6 /
\plot -12.5 8 -12.5 8.8 -12 9.6 /
\circulararc 35 degrees from -14.5 1.6 center at -15 0
\circulararc 35 degrees from -14 3.2 center at -15 0
\circulararc 35 degrees from -13.5 4.8 center at -15 0
\circulararc 35 degrees from -13 6.4 center at -15 0
\circulararc 35 degrees from -12.5 8 center at -15 0
\circulararc 35 degrees from -12 9.6 center at -15 0
\circulararc 150 degrees from 16 3.2 center at 15 0
\circulararc 150 degrees from 17 6.4 center at 15 0
\plot 16 3.2 16 4.8 17 6.4 /
\plot 16 3.2 17 4.8 17 6.4 /
\plot 14 3.2 14 4.8 13 6.4 /
\plot 14 3.2 13 4.8 13 6.4 /
\plot 11.8 1 10.2 1 8.6 2 /
\plot 11.8 1 10.2 2 8.6 2 /
\plot 11.8 -1 10.2 -2 8.6 -2 /
\plot 11.8 -1 10.2 -1 8.6 -2 /
\plot 12.6 2.4 11.1 3.3 10.2 4.8 /
\plot 12.6 2.4 11.7 3.9 10.2 4.8 /
\plot 12.6 -2.4 11.1 -3.3 10.2 -4.8 /
\plot 12.6 -2.4 11.7 -3.9 10.2 -4.8 /
\put{$C$} at -13 11
\put{$s_{1}C$} at -17 11
\put{$C$} at 17 10
\put{$s_{1}C$} at 13 10
\put{$s_{1}s_{2}C$} at 8 7
\put{$s_{1}s_{2}s_{1}C$} at 5 2
\put{$s_{1}s_{2}s_{1}s_{2}C$} at 5 -2
\put{$s_{1}s_{2}s_{1}s_{2}s_{1}C$} at 7 -7
\endpicture\]

\[\beginpicture
\setcoordinatesystem units <0.24cm,0.24cm>
\setplotarea x from -25 to 25, y from -15 to 12
\linethickness=3pt
\plot -15 -10 -15 10 / 
\setdots
\plot -20 -8.5 -10 8.5 / 
\plot -20 8.5 -10 -8.5 / 
\plot -25 0 -5 0 / 
\plot 10 -8.5 20 8.5 / 
\plot -23.5 -5 -6.5 5 / 
\plot -23.5 5 -6.5 -5 / 
\plot 5 0 25 0 / 
\plot 15 -10 15 10 / 
\plot 10 8.5 20 -8.5 / 
\plot 6.5 -5 23.5 5 / 
\plot 6.5 5 23.5 -5 / 
\put{$\bullet$} at 14 3.2
\put{$\bullet$} at 16 3.2
\put{$\bullet$} at 14 -3.2
\put{$\bullet$} at 16 -3.2
\put{$\bullet$} at 18.2 1
\put{$\bullet$} at 11.8 1
\put{$\bullet$} at 11.8 -1
\put{$\bullet$} at 18.2 -1
\put{$\bullet$} at 12.6 -2.4
\put{$\bullet$} at 12.6 2.4
\put{$\bullet$} at 17.4 2.4
\put{$\bullet$} at 17.4 -2.4
\put{$\bullet$} at -14 3.2
\put{$\bullet$} at -14 -3.2
\put{$\bullet$} at -13 6.4
\put{$\bullet$} at -13 -6.4
\put{$\bullet$} at -8.6 -2
\put{$\bullet$} at -11.8 1
\put{$\bullet$} at -11.8 -1
\put{$\bullet$} at -8.6 2
\put{$\bullet$} at -12.6 -2.4
\put{$\bullet$} at -12.6 2.4
\put{$\bullet$} at -10.2 -4.8
\put{$\bullet$} at -10.2 4.8
\setdashes
\put{$t_{1,z}, q^{2} \neq 1$, $z$ generic} at -15 -12
\put{$t_{z,w}$ $z,w$ generic} at 15 -12
\setquadratic
\setsolid
\circulararc 150 degrees from -14 -3.2 center at -15 0
\circulararc 150 degrees from -13 -6.4 center at -15 0
\plot -14 -3.2 -14 -4.8 -13 -6.4 /
\plot -14 -3.2 -13 -4.8 -13 -6.4 /
\plot -14 3.2 -14 4.8 -13 6.4 /
\plot -14 3.2 -13 4.8 -13 6.4 /
\plot -11.8 1 -10.2 1 -8.6 2 /
\plot -11.8 1 -10.2 2 -8.6 2 /
\plot -11.8 -1 -10.2 -2 -8.6 -2 /
\plot -11.8 -1 -10.2 -1 -8.6 -2 /
\plot -12.6 2.4 -11.1 3.3 -10.2 4.8 /
\plot -12.6 2.4 -11.7 3.9 -10.2 4.8 /
\plot -12.6 -2.4 -11.1 -3.3 -10.2 -4.8 /
\plot -12.6 -2.4 -11.7 -3.9 -10.2 -4.8 /
\circulararc 360 degrees from 14 3.2 center at 15 0
\put{$C$} at -13 11
\put{$s_{2}C$} at -8 7
\put{$s_{2}s_{1}C$} at -5 2
\put{$s_{2}s_{1}s_{2}C$} at -5 -2
\put{$s_{2}s_{1}s_{2}s_{1}C$} at -7 -7
\put{$s_{2}s_{1}s_{2}s_{1}s_{2}C$} at -11 -10
\put{$C$} at 17 10
\put{$s_{1}C$} at 13 10
\put{$s_{1}s_{2}C$} at 8 7
\put{$s_{1}s_{2}s_{1}C$} at 5 2
\put{$s_{1}s_{2}s_{1}s_{2}C$} at 5 -2
\put{$s_{1}s_{2}s_{1}s_{2}s_{1}C$} at 7 -7
\put{$s_{2}C$} at 22 7
\put{$s_{2}s_{1}C$} at 25 2
\put{$s_{2}s_{1}s_{2}C$} at 25 -2
\put{$s_{2}s_{1}s_{2}s_{1}C$} at 23 -7
\put{$s_{2}s_{1}s_{2}s_{1}s_{2}C$} at 19 -10
\put{$s_{1}s_{2}s_{1}s_{2}s_{1}s_{2}C$} at 10 -10
\endpicture\]

\underline{Case 2:} $Z(t) = \emptyset$.

If $Z(t) = \emptyset$ then $t$ is a regular central character. Then the irreducibles with central character $t$ are in bijection with the connected components of the calibration graph for $t$, and can be constructed using Theorem $\ref{eq:locreg}$. In particular, the module $H^{(t,C)}$ corresponding to a component $C$ of the calibration graph has all 1-dimensional weight spaces, and the weights of $H^{(t,C)}$ are precisely the vertices in $C$.

The following pictures show the modules $H^{(t,C)}$. There is one dot in the chamber $w^{-1}C$ for each $w \in C$, and the dots corresponding to $w$ and $s_{i}w$ are connected exactly when $wt$ and $s_{i}wt$ are weights in the same composition factor of $M(t)$.

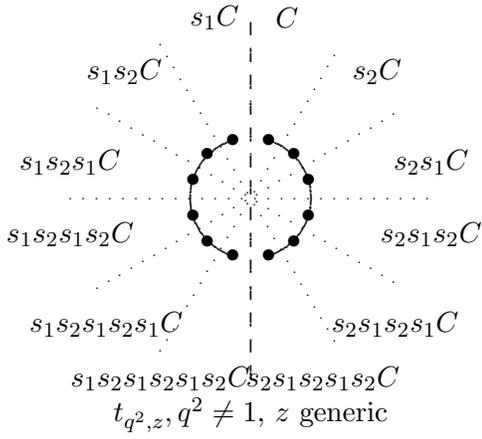
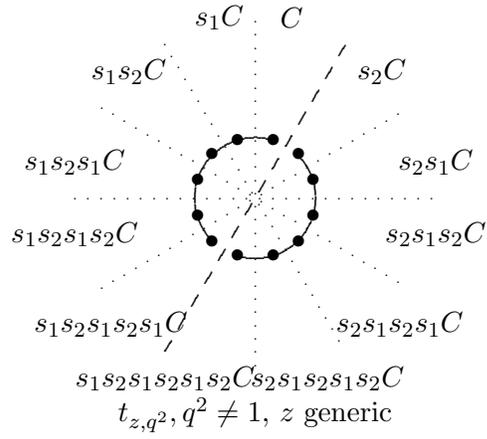
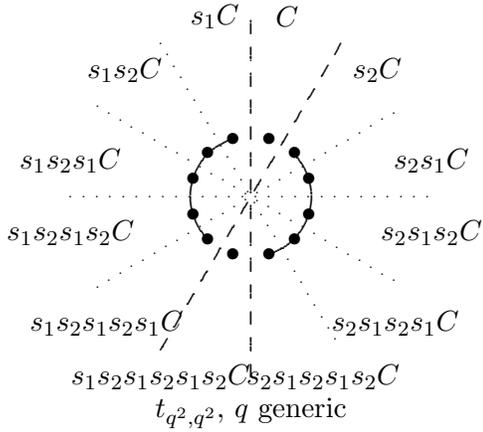
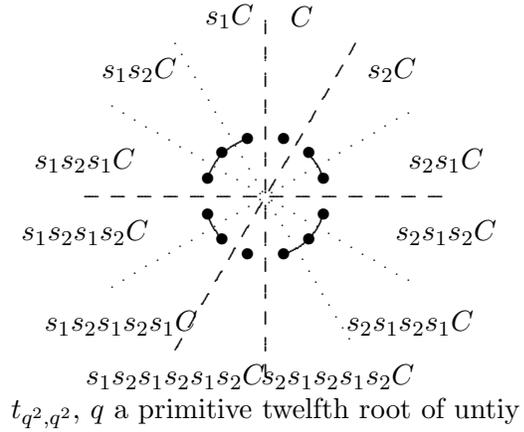
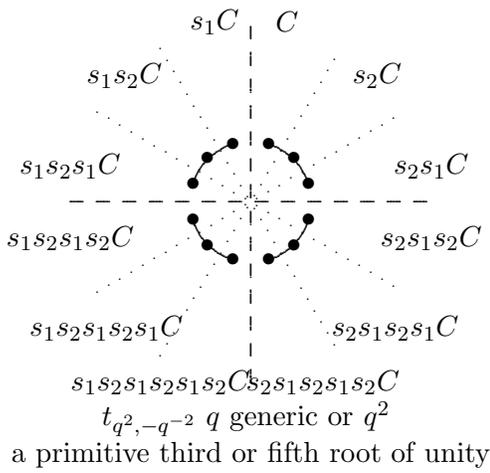
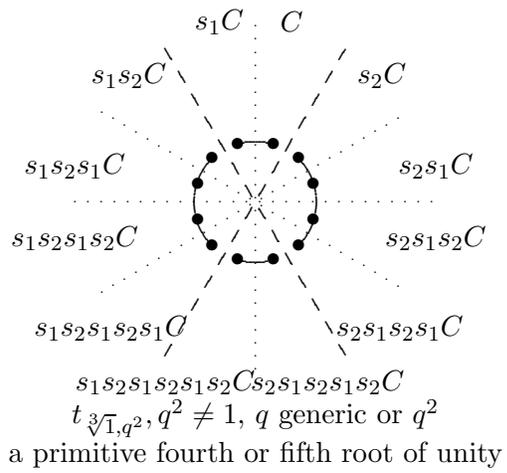
\[\beginpicture
\setcoordinatesystem units <0.24cm,0.24cm>
\setplotarea x from -25 to 25, y from -15 to 12
\linethickness=3pt
\setdashes
\plot -15 -10 -15 10 / 
\plot 10 -8.5 20 8.5 / 
\setdots
\plot -20 -8.5 -10 8.5 / 
\plot -20 8.5 -10 -8.5 / 
\plot -25 0 -5 0 / 
\plot -23.5 -5 -6.5 5 / 
\plot -23.5 5 -6.5 -5 / 
\plot 5 0 25 0 / 
\plot 15 -10 15 10 / 
\plot 10 8.5 20 -8.5 / 
\plot 6.5 -5 23.5 5 / 
\plot 6.5 5 23.5 -5 / 
\put{$\bullet$} at 14 3.2
\put{$\bullet$} at 16 3.2
\put{$\bullet$} at 14 -3.2
\put{$\bullet$} at 16 -3.2
\put{$\bullet$} at 18.2 1
\put{$\bullet$} at 11.8 1
\put{$\bullet$} at 11.8 -1
\put{$\bullet$} at 18.2 -1
\put{$\bullet$} at 12.6 -2.4
\put{$\bullet$} at 12.6 2.4
\put{$\bullet$} at 17.4 2.4
\put{$\bullet$} at 17.4 -2.4
\put{$\bullet$} at -14 3.2
\put{$\bullet$} at -16 3.2
\put{$\bullet$} at -14 -3.2
\put{$\bullet$} at -16 -3.2
\put{$\bullet$} at -18.2 1
\put{$\bullet$} at -11.8 1
\put{$\bullet$} at -11.8 -1
\put{$\bullet$} at -18.2 -1
\put{$\bullet$} at -12.6 -2.4
\put{$\bullet$} at -12.6 2.4
\put{$\bullet$} at -17.4 2.4
\put{$\bullet$} at -17.4 -2.4
\put{$t_{q^{2},z}, q^{2} \neq 1$, $z$ generic} at -15 -12
\put{$t_{z,q^{2}}, q^{2} \neq 1$, $z$ generic} at 15 -12
\setquadratic
\setsolid
\circulararc 150 degrees from -14 -3.2 center at -15 0
\circulararc 150 degrees from -16 3.2 center at -15 0
\circulararc 150 degrees from 16 3.2 center at 15 0
\circulararc 150 degrees from 14 -3.2 center at 15 0
\put{$C$} at -13 10
\put{$s_{2}C$} at -8 7
\put{$s_{2}s_{1}C$} at -5 2
\put{$s_{2}s_{1}s_{2}C$} at -5 -2
\put{$s_{2}s_{1}s_{2}s_{1}C$} at -7 -7
\put{$s_{2}s_{1}s_{2}s_{1}s_{2}C$} at -11 -10
\put{$s_{1}C$} at -17 10
\put{$s_{1}s_{2}C$} at -22 7
\put{$s_{1}s_{2}s_{1}C$} at -25 2
\put{$s_{1}s_{2}s_{1}s_{2}C$} at -25 -2
\put{$s_{1}s_{2}s_{1}s_{2}s_{1}C$} at -23 -7
\put{$s_{1}s_{2}s_{1}s_{2}s_{1}s_{2}C$} at -20 -10
\put{$C$} at 17 10
\put{$s_{1}C$} at 13 10
\put{$s_{1}s_{2}C$} at 8 7
\put{$s_{1}s_{2}s_{1}C$} at 5 2
\put{$s_{1}s_{2}s_{1}s_{2}C$} at 5 -2
\put{$s_{1}s_{2}s_{1}s_{2}s_{1}C$} at 7 -7
\put{$s_{2}C$} at 22 7
\put{$s_{2}s_{1}C$} at 25 2
\put{$s_{2}s_{1}s_{2}C$} at 25 -2
\put{$s_{2}s_{1}s_{2}s_{1}C$} at 23 -7
\put{$s_{2}s_{1}s_{2}s_{1}s_{2}C$} at 19 -10
\put{$s_{1}s_{2}s_{1}s_{2}s_{1}s_{2}C$} at 10 -10
\endpicture\]

\[\beginpicture
\setcoordinatesystem units <0.24cm,0.24cm>
\setplotarea x from -25 to 25, y from -15 to 12
\linethickness=3pt
\setdots
\plot -20 8.5 -10 -8.5 / 
\plot -25 0 -5 0 / 
\plot -23.5 -5 -6.5 5 / 
\plot -23.5 5 -6.5 -5 / 
\plot 10 8.5 20 -8.5 / 
\plot 6.5 -5 23.5 5 / 
\plot 6.5 5 23.5 -5 / 
\put{$\bullet$} at 14 3.2
\put{$\bullet$} at 16 3.2
\put{$\bullet$} at 14 -3.2
\put{$\bullet$} at 16 -3.2
\put{$\bullet$} at 18.2 1
\put{$\bullet$} at 11.8 1
\put{$\bullet$} at 11.8 -1
\put{$\bullet$} at 18.2 -1
\put{$\bullet$} at 12.6 -2.4
\put{$\bullet$} at 12.6 2.4
\put{$\bullet$} at 17.4 2.4
\put{$\bullet$} at 17.4 -2.4
\put{$\bullet$} at -14 3.2
\put{$\bullet$} at -16 3.2
\put{$\bullet$} at -14 -3.2
\put{$\bullet$} at -16 -3.2
\put{$\bullet$} at -18.2 1
\put{$\bullet$} at -11.8 1
\put{$\bullet$} at -11.8 -1
\put{$\bullet$} at -18.2 -1
\put{$\bullet$} at -12.6 -2.4
\put{$\bullet$} at -12.6 2.4
\put{$\bullet$} at -17.4 2.4
\put{$\bullet$} at -17.4 -2.4
\setdashes
\plot 5 0 25 0 / 
\plot 15 -10 15 10 / 
\plot 10 -8.5 20 8.5 / 
\plot -15 -10 -15 10 / 
\plot -20 -8.5 -10 8.5 / 
\put{$t_{q^{2},q^{2}}$, $q$ generic} at -15 -12
\put{$t_{q^{2},q^{2}}$, $q$ a primitive twelfth root of untiy} at 15 -12
\setquadratic
\setsolid
\circulararc 120 degrees from -14 -3.2 center at -15 0
\circulararc 120 degrees from -16 3.2 center at -15 0
\circulararc 60 degrees from 14 3.2 center at 15 0
\circulararc 60 degrees from 16 -3.2 center at 15 0
\circulararc 30 degrees from 18.2 1 center at 15 0
\circulararc 30 degrees from 11.8 -1 center at 15 0
\put{$C$} at -13 10
\put{$s_{2}C$} at -8 7
\put{$s_{2}s_{1}C$} at -5 2
\put{$s_{2}s_{1}s_{2}C$} at -5 -2
\put{$s_{2}s_{1}s_{2}s_{1}C$} at -7 -7
\put{$s_{2}s_{1}s_{2}s_{1}s_{2}C$} at -11 -10
\put{$s_{1}C$} at -17 10
\put{$s_{1}s_{2}C$} at -22 7
\put{$s_{1}s_{2}s_{1}C$} at -25 2
\put{$s_{1}s_{2}s_{1}s_{2}C$} at -25 -2
\put{$s_{1}s_{2}s_{1}s_{2}s_{1}C$} at -23 -7
\put{$s_{1}s_{2}s_{1}s_{2}s_{1}s_{2}C$} at -20 -10
\put{$C$} at 17 10
\put{$s_{1}C$} at 13 10
\put{$s_{1}s_{2}C$} at 8 7
\put{$s_{1}s_{2}s_{1}C$} at 5 2
\put{$s_{1}s_{2}s_{1}s_{2}C$} at 5 -2
\put{$s_{1}s_{2}s_{1}s_{2}s_{1}C$} at 7 -7
\put{$s_{2}C$} at 22 7
\put{$s_{2}s_{1}C$} at 25 2
\put{$s_{2}s_{1}s_{2}C$} at 25 -2
\put{$s_{2}s_{1}s_{2}s_{1}C$} at 23 -7
\put{$s_{2}s_{1}s_{2}s_{1}s_{2}C$} at 19 -10
\put{$s_{1}s_{2}s_{1}s_{2}s_{1}s_{2}C$} at 10 -10
\endpicture\]

\[\beginpicture
\setcoordinatesystem units <0.24cm,0.24cm>
\setplotarea x from -25 to 25, y from -15 to 12
\linethickness=3pt
\setdots
\plot -20 8.5 -10 -8.5 / 
\plot -20 -8.5 -10 8.5 / 
\plot -23.5 -5 -6.5 5 / 
\plot -23.5 5 -6.5 -5 / 
\plot 5 0 25 0 / 
\plot 15 -10 15 10 / 
\plot 6.5 -5 23.5 5 / 
\plot 6.5 5 23.5 -5 / 
\put{$\bullet$} at 14 3.2
\put{$\bullet$} at 16 3.2
\put{$\bullet$} at 14 -3.2
\put{$\bullet$} at 16 -3.2
\put{$\bullet$} at 18.2 1
\put{$\bullet$} at 11.8 1
\put{$\bullet$} at 11.8 -1
\put{$\bullet$} at 18.2 -1
\put{$\bullet$} at 12.6 -2.4
\put{$\bullet$} at 12.6 2.4
\put{$\bullet$} at 17.4 2.4
\put{$\bullet$} at 17.4 -2.4
\put{$\bullet$} at -14 3.2
\put{$\bullet$} at -16 3.2
\put{$\bullet$} at -14 -3.2
\put{$\bullet$} at -16 -3.2
\put{$\bullet$} at -18.2 1
\put{$\bullet$} at -11.8 1
\put{$\bullet$} at -11.8 -1
\put{$\bullet$} at -18.2 -1
\put{$\bullet$} at -12.6 -2.4
\put{$\bullet$} at -12.6 2.4
\put{$\bullet$} at -17.4 2.4
\put{$\bullet$} at -17.4 -2.4
\setdashes
\plot 10 8.5 20 -8.5 / 
\plot 10 -8.5 20 8.5 / 
\plot -25 0 -5 0 / 
\plot -15 -10 -15 10 / 
\put{$t_{\sqrt[3]{1},q^{2}}, q^{2} \neq 1$, $q$ generic or $q^{2}$} at 15 -12
\put{a primitive third or fifth root of unity} at -15 -14
\put{$t_{q^{2},-q^{-2}}$ $q$ generic or $q^{2}$ } at -15 -12
\put{a primitive fourth or fifth root of unity} at 15 -14
\setquadratic
\setsolid
\circulararc 60 degrees from -14 -3.2 center at -15 0
\circulararc 60 degrees from -16 3.2 center at -15 0
\circulararc 60 degrees from -18.2 -1 center at -15 0
\circulararc 60 degrees from -11.8 1 center at -15 0
\circulararc 30 degrees from 16 3.2 center at 15 0
\circulararc 30 degrees from 14 -3.2 center at 15 0
\circulararc 90 degrees from 17.4 -2.4 center at 15 0
\circulararc 90 degrees from 12.6 2.4 center at 15 0
\put{$C$} at -13 10
\put{$s_{2}C$} at -8 7
\put{$s_{2}s_{1}C$} at -5 2
\put{$s_{2}s_{1}s_{2}C$} at -5 -2
\put{$s_{2}s_{1}s_{2}s_{1}C$} at -7 -7
\put{$s_{2}s_{1}s_{2}s_{1}s_{2}C$} at -11 -10
\put{$s_{1}C$} at -17 10
\put{$s_{1}s_{2}C$} at -22 7
\put{$s_{1}s_{2}s_{1}C$} at -25 2
\put{$s_{1}s_{2}s_{1}s_{2}C$} at -25 -2
\put{$s_{1}s_{2}s_{1}s_{2}s_{1}C$} at -23 -7
\put{$s_{1}s_{2}s_{1}s_{2}s_{1}s_{2}C$} at -20 -10
\put{$C$} at 17 10
\put{$s_{1}C$} at 13 10
\put{$s_{1}s_{2}C$} at 8 7
\put{$s_{1}s_{2}s_{1}C$} at 5 2
\put{$s_{1}s_{2}s_{1}s_{2}C$} at 5 -2
\put{$s_{1}s_{2}s_{1}s_{2}s_{1}C$} at 7 -7
\put{$s_{2}C$} at 22 7
\put{$s_{2}s_{1}C$} at 25 2
\put{$s_{2}s_{1}s_{2}C$} at 25 -2
\put{$s_{2}s_{1}s_{2}s_{1}C$} at 23 -7
\put{$s_{2}s_{1}s_{2}s_{1}s_{2}C$} at 19 -10
\put{$s_{1}s_{2}s_{1}s_{2}s_{1}s_{2}C$} at 10 -10
\endpicture\]

\underline{Case 3:} $Z(t), P(t) \neq \emptyset$.

For these central characters, rather than analyzing $M(t)$ directly, it is easier to construct several irreducible $\widetilde{H}$-modules and show that they include all the composition factors of $M(t)$.

\noindent \textbf{Case 3a: $t_{1,q^{2}}$}

Assume $\alpha_{1} \in Z(t)$ and $\alpha_{2} \in P(t)$. Then $t = t_{1,q^{\pm 2}}$, but $s_{2}s_{1}s_{2}s_{1}s_{2}t_{1,q^{-2}} = t_{1,q^{2}}$, so that analyzing $M(t_{1,q^{2}})$ is sufficient. Then let $t = t_{1,q^{2}}$. We have $Z(t) = \{ \alpha_{1} \}$ and $P(t) = \{ \alpha_{2}\}$ unless $q^{2} = \pm 1$. Hence the cases $q^{2} = 1$ and $q^{2} = -1$ will be treated separately.

If $q^{2} = 1$, then $Z(t) = P(t) = R^{+}$, and the irreducibles with this central character can be constructed using the results of $\eqref{clifford}$. In this case, since the action of $\mathbb{C}[X]$ is trivial, we can consider the irreducible $\widetilde{H}$-modules as $\mathbb{C}[W_{0}]$-modules, and the principal series module is the regular representation of $W_{0}$. Specifically, there are four 1-dimensional modules and two 2-dimensional modules with central character $t$, and the 2-dimensional modules appear in $M(t)$ with multiplicity 2.

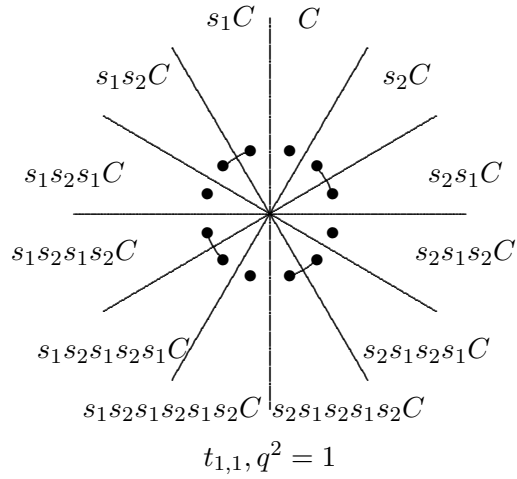
\[\beginpicture
\setcoordinatesystem units <0.26cm,0.26cm>
\setplotarea x from -25 to 25, y from -15 to 11
\linethickness=3pt
\plot -5 8.5 5 -8.5 / 
\plot -8.5 -5 8.5 5 / 
\plot -8.5 5 8.5 -5 / 
\plot -5 -8.5 5 8.5 / 
\plot -10 0 10 0 / 
\plot 0 -10 0 10 / 
\put{$\bullet$} at 1 3.2
\put{$\bullet$} at -1 3.2
\put{$\bullet$} at 1 -3.2
\put{$\bullet$} at -1 -3.2
\put{$\bullet$} at -3.2 1
\put{$\bullet$} at 3.2 1
\put{$\bullet$} at 3.2 -1
\put{$\bullet$} at -3.2 -1
\put{$\bullet$} at 2.4 -2.4
\put{$\bullet$} at 2.4 2.4
\put{$\bullet$} at -2.4 2.4
\put{$\bullet$} at -2.4 -2.4
\setdashes
\put{$t_{1,1}, q^{2} = 1 $} at 0 -12.5
\setquadratic
\setsolid
\circulararc 30 degrees from 1 -3.2 center at 0 0
\circulararc 30 degrees from -1 3.2 center at 0 0
\circulararc 30 degrees from -3.2 -1 center at 0 0
\circulararc 30 degrees from 3.2 1 center at 0 0
\put{$C$} at 2 10
\put{$s_{2}C$} at 7 7
\put{$s_{2}s_{1}C$} at 10 2
\put{$s_{2}s_{1}s_{2}C$} at 10 -2
\put{$s_{2}s_{1}s_{2}s_{1}C$} at 8 -7
\put{$s_{2}s_{1}s_{2}s_{1}s_{2}C$} at 4 -10
\put{$s_{1}C$} at -2 10
\put{$s_{1}s_{2}C$} at -7 7
\put{$s_{1}s_{2}s_{1}C$} at -10 2
\put{$s_{1}s_{2}s_{1}s_{2}C$} at -10 -2
\put{$s_{1}s_{2}s_{1}s_{2}s_{1}C$} at -8 -7
\put{$s_{1}s_{2}s_{1}s_{2}s_{1}s_{2}C$} at -5 -10
\endpicture\]

If $q^{2} \neq 1$, let $w_{0} =s_{1}s_{2}s_{1}s_{2}s_{1}s_{2}$ and define \[\widetilde{H}_{\{2\}} = \mathbb{C}\textrm{-span}\{ T_{2}X^{\lambda}, X^{\lambda} \, | \, \lambda \in P \},\] the subalgebra of $\widetilde{H}$ generated by $T_{2}$ and $\mathbb{C}[X]$. Let $\mathbb{C}v_{t}$ and $\mathbb{C}v_{w_{0}t}$ be the 1-dimensional $\widetilde{H}_{\{2\}}$-modules spanned by $v_{t}$ and $v_{w_{0}t}$, respectively, and given by \[ \begin{matrix} T_{2}v_{t} = qv_{t}, & \quad & X^{\lambda}v_{t} = t(X^{\lambda})v_{t}, \\ T_{2}v_{w_{0}t} = -q^{-1}v_{w_{0}t}, & \textrm{ and } & X^{\lambda}v_{w_{0}t} = w_{0}t(X^{\lambda})v_{w_{0}t}. \end{matrix}\] Then define \[ M = \mathbb{C}v_{t} \otimes_{\widetilde{H}_{\{2\}}} \widetilde{H} \quad \textrm{ and } \] \[ N = \mathbb{C}v_{w_{0}t} \otimes_{\widetilde{H}_{\{2\}}} \widetilde{H}.\]

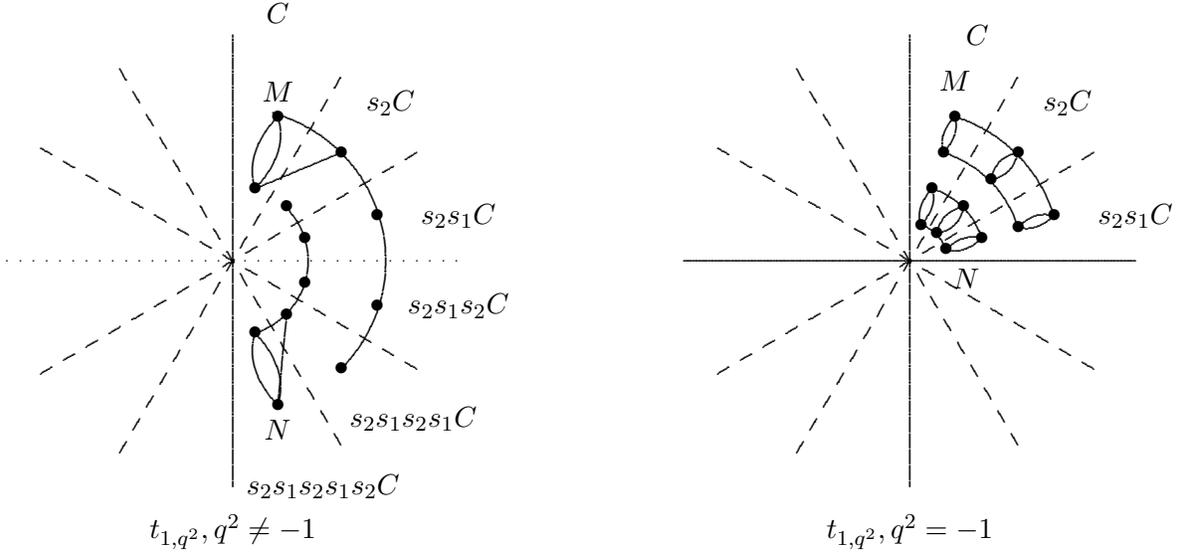
\[\beginpicture
\setcoordinatesystem units <0.3cm,0.3cm>
\setplotarea x from -25 to 25, y from -15 to 15
\linethickness=3pt
\plot -15 -10 -15 10 / 
\plot 15 -10 15 10 / 
\plot 5 0 25 0 / 
\setdots
\plot -25 0 -5 0 / 
\put{$\bullet$} at 15.5 1.6
\put{$\bullet$} at 16 3.2
\put{$\bullet$} at 16.5 4.8
\put{$\bullet$} at 17 6.4
\put{$\bullet$} at 21.4 2
\put{$\bullet$} at 16.6 0.5
\put{$\bullet$} at 18.2 1
\put{$\bullet$} at 19.8 1.5
\put{$\bullet$} at 16.2 1.2
\put{$\bullet$} at 17.4 2.4
\put{$\bullet$} at 18.6 3.6
\put{$\bullet$} at 19.8 4.8
\put{$\bullet$} at -14 3.2
\put{$\bullet$} at -14 -3.2
\put{$\bullet$} at -13 6.4
\put{$\bullet$} at -13 -6.4
\put{$\bullet$} at -8.6 -2
\put{$\bullet$} at -11.8 1
\put{$\bullet$} at -11.8 -1
\put{$\bullet$} at -8.6 2
\put{$\bullet$} at -12.6 -2.4
\put{$\bullet$} at -12.6 2.4
\put{$\bullet$} at -10.2 -4.8
\put{$\bullet$} at -10.2 4.8
\setdashes
\plot 10 -8.5 20 8.5 / 
\plot -23.5 -5 -6.5 5 / 
\plot -23.5 5 -6.5 -5 / 
\plot 10 8.5 20 -8.5 / 
\plot 6.5 -5 23.5 5 / 
\plot 6.5 5 23.5 -5 / 
\plot -20 -8.5 -10 8.5 / 
\plot -20 8.5 -10 -8.5 / 
\put{$t_{1,q^{2}},q^{2} \neq -1$} at -15 -12
\put{$t_{1,q^{2}}, q^{2} = -1$} at 15 -12
\put{$M$} at -13 7.5
\put{$N$} at -13 -7.5
\setquadratic
\setsolid
\circulararc 120 degrees from -14 -3.2 center at -15 0
\circulararc 120 degrees from -10.2 -4.8 center at -15 0
\circulararc 57 degrees from 16.6 0.5 center at 15 0
\circulararc 57 degrees from 18.2 1 center at 15 0
\circulararc 57 degrees from 19.8 1.5 center at 15 0
\circulararc 57 degrees from 21.4 2 center at 15 0
\plot -14 -3.2 -14 -4.8 -13 -6.4 /
\plot -14 -3.2 -13 -4.8 -13 -6.4 /
\plot -14 3.2 -14 4.8 -13 6.4 /
\plot -14 3.2 -13 4.8 -13 6.4 /
\plot -14 3.2 -12.1 4 -10.2 4.8 /
\plot -13 -6.4 -12.8 -4.4 -12.6 -2.4 /
\plot 15.5 1.6 16 2.4 16 3.2 /
\plot 15.5 1.6 15.5 2.4 16 3.2 /
\plot 16.5 4.8 16.5 5.6 17 6.4 /
\plot 16.5 4.8 17 5.6 17 6.4 /
\plot 16.2 1.2 16.6 2 17.4 2.4 /
\plot 16.2 1.2 17 1.6 17.4 2.4 /
\plot 18.6 3.6 19 4.4 19.8 4.8 /
\plot 18.6 3.6 19.4 4 19.8 4.8 /
\plot 16.6 0.5 17.4 1 18.2 1 /
\plot 16.6 0.5 17.4 0.5 18.2 1 /
\plot 19.8 1.5 20.6 2 21.4 2 /
\plot 19.8 1.5 20.6 1.5 21.4 2 /
\put{$M$} at 17 8
\put{$N$} at 17.5 -0.8
\put{$C$} at -13 11
\put{$s_{2}C$} at -8 7
\put{$s_{2}s_{1}C$} at -5 2
\put{$s_{2}s_{1}s_{2}C$} at -5 -2
\put{$s_{2}s_{1}s_{2}s_{1}C$} at -7 -7
\put{$s_{2}s_{1}s_{2}s_{1}s_{2}C$} at -11 -10
\put{$C$} at 18 10
\put{$s_{2}C$} at 22 7
\put{$s_{2}s_{1}C$} at 25 2
\endpicture\]

\begin{proposition}
Assume $q^{2} \neq \pm 1$. Let $M = \mathbb{C}v_{t} \otimes_{\widetilde{H}_{\{2\}}} \widetilde{H} \quad \textrm{ and } N = \mathbb{C}v_{w_{0}t} \otimes_{\widetilde{H}_{\{2\}}} \widetilde{H},$ where $t = t_{1,q^{2}}$.

(a) $M_{s_{1}s_{2}s_{1}s_{2}t}$ is a 1-dimensional submodule of $M$. $M'$, the image of the weight spaces $M_{s_{2}s_{1}s_{2}t}$ and $M_{s_{1}s_{2}t}$ in $M/M_{s_{1}s_{2}s_{1}s_{2}t}$, is a submodule of $M/M_{s_{1}s_{2}s_{1}s_{2}t}$. The resulting quotient of $M$ is irreducible.

(b) If $q^{2}$ is not a primitive third root of unity, then $M'$ is irreducible.

(c) If $q^{2}$ is a primitive third root of unity, then $(M')_{s_{2}s_{1}s_{2}t}$ is a submodule of $M'$.

(d) $N_{s_{2}t}$ is a 1-dimensional submodule of $N$. $N'$, the image of the weight spaces $N_{s_{1}s_{2}t}$ and $N_{s_{2}s_{1}s_{2}t}$ in $N/N_{s_{2}t}$, is a submodule of $N/N_{s_{2}t}$. The resulting quotient of $N$ is irreducible.

(e) If $q^{2}$ is not a primitive third root of unity, then $N'$ is irreducible.

(f) If $q^{2}$ is a primitive third root of unity, then $(N')_{s_{1}s_{2}t}$ is a submodule of $N'$.

\end{proposition}

\noindent \textit{Proof.} Assume $q^{2} \neq -1$. Then $Z(t) = \{ \alpha_{1} \}$ and $P(t)$ contains $\alpha_{2}, \alpha_{1} + \alpha_{2}, 2\alpha_{1} + \alpha_{2},$ and $3\alpha_{1} + \alpha_{2}$. If $q^{2}$ is a primitive third root of unity, then $P(t)$ also contains $3\alpha_{1} + 2\alpha_{2}$. Then $M$ has one 2-dimensional weight space $M_{t}^{\textrm{gen}}$ and four 1-dimensional weight spaces $M_{s_{2}t}$, $M_{s_{1}s_{2}t}$, $M_{s_{2}s_{1}s_{2}t}$, and $M_{s_{1}s_{2}s_{1}s_{2}t}$. For $w \in \{ s_{2}, s_{1}s_{2}, s_{2}s_{1}s_{2}, s_{1}s_{2}s_{1}s_{2}\}$, let $m_{wt}$ be a non-zero vector in $M_{wt}$. By a calculation as in \eqref{weightbasis}, \[ m_{wt} = T_{w}T_{1}v_{t} + \sum_{w' < w} a_{w,w'}T_{w'}T_{1}v_{t},\] for $w \in \{ s_{2}, s_{1}s_{2}, s_{2}s_{1}s_{2}, s_{1}s_{2}s_{1}s_{2}\}$, where $a_{w,w'} \in \mathbb{C}$. Then if $s_{i}w > w$, $\tau_{i}m_{wt} \neq 0$ for $w \in \{ s_{2}, s_{1}s_{2}, s_{2}s_{1}s_{2} \}$, since the term $T_{i}T_{w}T_{1}$ cannot be canceled by any other term in $\tau_{i}m_{wt}$.

Thus $\tau_{1}: M_{s_{1}s_{2}s_{1}s_{2}t} \rightarrow M_{s_{2}s_{1}s_{2}t}$ is the zero map since, by Theorem \ref{tauthm}, $\tau_{1}^{2}: M_{s_{2}s_{1}s_{2}}t \rightarrow M_{s_{2}s_{1}s_{2}t}$ is the zero map. Hence $M_{s_{1}s_{2}s_{1}s_{2}t}$ is a submodule of $M$. Similarly, $\tau_{1}:M_{s_{1}s_{2}t} \rightarrow M_{s_{2}t}$ must be the zero map since, by Theorem \ref{tauthm}, $\tau_{1}^{2}:M_{s_{2}t} \rightarrow M_{s_{2}t}$ is the zero map. Let $M_{1} = M/M_{s_{1}s_{2}s_{1}s_{2}t}$. Then $M'$, the subspace spanned by $\overline{m_{s_{1}s_{2}t}}$ and $\overline{m_{s_{2}s_{1}s_{2}t}}$ in $M_{1}$, is a submodule of $M_{1}$. Lemma $\ref{eq:lemma}$ shows that $M_{2} = M_{1}/M'$ is irreducible.

(b) If $q^{2}$ is not a primitive third root of unity, $\tau_{2}^{2}: (M')_{s_{1}s_{2}t} \rightarrow (M')_{s_{1}s_{2}t}$ is invertible, so that $M'$ is irreducible.

(c) If $q^{2}$ is a primitive third root of unity, then $\tau_{2}:(M'_{1})_{s_{2}s_{1}s_{2}t} \rightarrow (M'_{1})_{s_{1}s_{2}t}$ is the zero map and $(M'_{1})_{s_{2}s_{1}s_{2}t}$ is a 1-dimensional submodule of $M'_{1}$, and $M'_{1}/(M'_{1})_{s_{2}s_{1}s_{2}t}$ is 1-dimensional as well.

(d)-(f) The same argument used in (a)-(c) applies, with each weight space $M_{wt}$ replaced by $N_{ww_{0}t}$. $\square$

\vspace{.2in}

However, the composition factors of $M$ and $N$ are not distinct. If $q^{2}$ is not a primitive third root of unity, then $M'$ and $N'$ are irreducible 2-dimensional modules with the same weight space structure. Then Theorem \ref{2dims} shows that $M' \cong N'$. If $q^{2}$ is a primitive third root of unity, then note that two 1-dimensional modules $\mathbb{C}v_{t}$ and $\mathbb{C}v_{t'}$ are isomorphic if and only if they have the same weight. Then $M'_{1}$ and $N'_{1}$ have the same composition factors. In any case, the 3-dimensional modules are different since their weight space structures are different.

\[\beginpicture
\setcoordinatesystem units <0.274cm,0.274cm>
\setplotarea x from -25 to 25, y from -15 to 15
\linethickness=3pt
\plot -15 -10 -15 10 / 
\plot 15 -10 15 10 / 
\setdots
\plot 5 0 25 0 / 
\put{$\bullet$} at 16 3.2
\put{$\bullet$} at 16 -3.2
\put{$\bullet$} at 17 6.4
\put{$\bullet$} at 17 -6.4
\put{$\bullet$} at 21.4 -2
\put{$\bullet$} at 18.2 1
\put{$\bullet$} at 18.2 -1
\put{$\bullet$} at 21.4 2
\put{$\bullet$} at 17.4 -2.4
\put{$\bullet$} at 17.4 2.4
\put{$\bullet$} at 19.8 -4.8
\put{$\bullet$} at 19.8 4.8
\put{$\bullet$} at -14 3.2
\put{$\bullet$} at -14 -3.2
\put{$\bullet$} at -13 6.4
\put{$\bullet$} at -13 -6.4
\put{$\bullet$} at -8.6 -2
\put{$\bullet$} at -11.8 1
\put{$\bullet$} at -11.8 -1
\put{$\bullet$} at -8.6 2
\put{$\bullet$} at -12.6 -2.4
\put{$\bullet$} at -12.6 2.4
\put{$\bullet$} at -10.2 -4.8
\put{$\bullet$} at -10.2 4.8
\setdashes
\plot -25 0 -5 0 / 
\plot 10 -8.5 20 8.5 / 
\plot -23.5 -5 -6.5 5 / 
\plot -23.5 5 -6.5 -5 / 
\plot 10 8.5 20 -8.5 / 
\plot 6.5 -5 23.5 5 / 
\plot 6.5 5 23.5 -5 / 
\plot -20 -8.5 -10 8.5 / 
\plot -20 8.5 -10 -8.5 / 
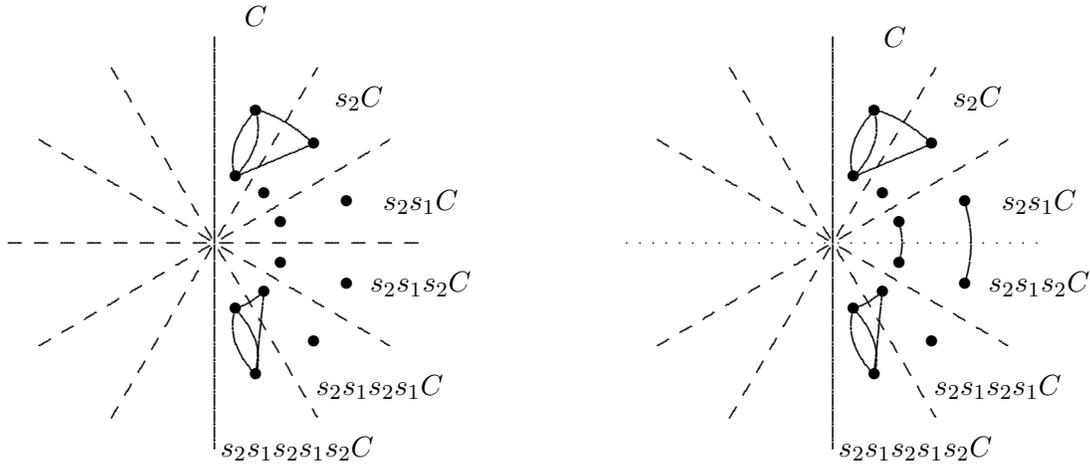
\put{$t_{1,q^{2}}, q^{2}$ a primitive third root of unity} at -15 -12
\put{$t_{1,q^{2}}, q^{2}$ not a primitive third root of unity, $q^{2} \neq -1$} at 15 -12
\setquadratic
\setsolid
\circulararc 30 degrees from -10.2 4.8 center at -15 0
\circulararc 30 degrees from -14 -3.2 center at -15 0
\circulararc 30 degrees from 16 -3.2 center at 15 0
\circulararc 30 degrees from 19.8 4.8 center at 15 0
\circulararc 35 degrees from 21.4 -2 center at 15 0
\circulararc 35 degrees from 18.2 -1 center at 15 0
\plot -14 -3.2 -14 -4.8 -13 -6.4 /
\plot -14 -3.2 -13 -4.8 -13 -6.4 /
\plot -14 3.2 -14 4.8 -13 6.4 /
\plot -14 3.2 -13 4.8 -13 6.4 /
\plot -14 3.2 -12.1 4 -10.2 4.8 /
\plot -13 -6.4 -12.8 -4.4 -12.6 -2.4 /
\plot 16 -3.2 16 -4.8 17 -6.4 /
\plot 16 -3.2 17 -4.8 17 -6.4 /
\plot 16 3.2 16 4.8 17 6.4 /
\plot 16 3.2 17 4.8 17 6.4 /
\plot 16 3.2 17.9 4 19.8 4.8 /
\plot 17 -6.4 17.2 -4.4 17.4 -2.4 /
\put{$C$} at -13 11
\put{$s_{2}C$} at -8 7
\put{$s_{2}s_{1}C$} at -5 2
\put{$s_{2}s_{1}s_{2}C$} at -5 -2
\put{$s_{2}s_{1}s_{2}s_{1}C$} at -7 -7
\put{$s_{2}s_{1}s_{2}s_{1}s_{2}C$} at -11 -10
\put{$C$} at 18 10
\put{$s_{2}C$} at 22 7
\put{$s_{2}s_{1}C$} at 25 2
\put{$s_{2}s_{1}s_{2}C$} at 25 -2
\put{$s_{2}s_{1}s_{2}s_{1}C$} at 23 -7
\put{$s_{2}s_{1}s_{2}s_{1}s_{2}C$} at 19 -10
\endpicture\]

\begin{proposition} If $q^{2} \neq  \pm 1$ then the composition factors of $M$ and $N$ are the only irreducibles with central character $t_{1,q^{2}}$.
\end{proposition}

\noindent \textit{Proof.} Counting multiplicities of weight spaces in $M(t)$ and the distinct composition factors of $M$ and $N$ shows that the remaining composition factor(s) of $M(t)$ must contain an $s_{1}s_{2}t$ weight space and an $s_{2}s_{1}s_{2}t$ weight space, each of dimension 1.

If $q^{2}$ is not a primitive third root of unity then Theorem \eqref{calibsame}b shows that there must be one remaining composition factor with an $s_{1}s_{2}t$ weight space and an $s_{2}s_{1}s_{2}t$ weight space, and Theorem \ref{2dims} shows that it is isomorphic to $M_{1}$.

If $q^{2}$ is a primitive third root of unity then $\tau_{2}^{2} M_{s_{1}s_{2}t} \rightarrow M_{s_{1}s_{2}t}$ is not invertible. Hence there cannot be an irreducible module consisting of an $s_{1}s_{2}t$ weight space and an $s_{2}s_{1}s_{2}$ weight space, and the remaining composition factors of $M(t)$ are 1-dimensional. $\square$

\vspace{.2in}

If $q^{2} = -1$, then dim $M_{t}^{\textrm{gen}} = $ dim $M_{s_{2}t}^{\textrm{gen}} =$ dim $M_{s_{1}s_{2}t}^{\textrm{gen}} = 2$ and dim $N_{t}^{\textrm{gen}} = $ dim $N_{s_{2}t}^{\textrm{gen}} =$ dim $N_{s_{1}s_{2}t}^{\textrm{gen}} = 2$.

\begin{proposition} Assume $q^{2} = -1$ and $t = t_{1,q^{2}}$. Let $M = \widetilde{H} \otimes_{\widetilde{H}_{\{2\}}} \mathbb{C}v_{t} \quad \textrm{ and } N = \widetilde{H} \otimes_{\widetilde{H}_{\{1\}}} \mathbb{C}v_{s_{1}s_{2}t}.$

(a) $M$ and $N$ each have two 1-dimensional modules and two 2-dimensional modules as composition factors.

(b) The composition factors of $M$ and $N$ are the only irreducible modules with central character $t$.

\end{proposition}

\noindent \textit{Proof.} By Proposition \ref{2dims}, there is a 2-dimensional module $P$ with $P = P_{t}^{\textrm{gen}}$. Let $v \in P_{t}$ be non-zero. The map
\renewcommand{\baselinestretch}{1} \normalsize\[ \begin{array}{lll} \mathbb{C}v_{t} &  \rightarrow & P \\ v_{t} &  \mapsto & v \end{array}\] \renewcommand{\baselinestretch}{1} \normalsize
is a $\widetilde{H}_{\{2\}}$-module homomorphism. Since \[ \textrm{Hom}_{\widetilde{H}} (M, P) = \textrm{Hom}_{\widetilde{H}_{\{2\}}}(\mathbb{C}v_{t}, P),\] there is a non-zero map from $M$ to $P$. Since $P$ is irreducible, this map is surjective and $P$ is a quotient of $M$. The kernel of any map from $M$ to $P$ must be \[M_{1} = M_{s_{2}t}^{\textrm{gen}} \oplus M_{s_{1}s_{2}t}^{\textrm{gen}},\] which is then a submodule of $M$.

Then we note that $m = T_{1}T_{2}T_{1}T_{2}T_{1}v - qT_{2}T_{1}T_{2}T_{1}v - T_{1}T_{2}T_{1}v + qT_{2}T_{1}v + T_{1}v - qv$ spans a 1-dimensional submodule of $M_{1}$. Then let $M_{2} = M_{1}/ m$, so that $T_{1}T_{2}T_{1}T_{2}T_{1}v = qT_{2}T_{1}T_{2}T_{1}v + T_{1}T_{2}T_{1}v - qT_{2}T_{1}v - T_{1}v + qv$ in $M_{2}$.

Then by Theorem $\ref{weightbasis}$, $M_{2}$ contains an element $m' = T_{2}T_{1}v - qT_{1}v -3v \in M_{s_{2}t}$. Then $m'$, $\tau_{1}(m')$ and $T_{2} \cdot \tau_{1}(m')$ are linearly independent (since their leading terms cannot be canceled) and span $M_{2}$. However, $M_{3} = \langle \tau_{1}(m'), T_{2} \cdot \tau_{1}(m') \rangle$ is clearly closed under the action $T_{2}$. Also, $\tau_{1}(m') \in M_{s_{1}s_{2}t}$, so that \[ X^{\lambda} \cdot T_{2} \tau_{1}(m') = T_{2}X^{s_{2}\lambda} \tau_{1}(m') + (q-q^{-1})\frac{X^{\lambda}-X^{s_{2}\lambda}}{1-X^{-\alpha_{2}}} \tau_{1}(m'),\] which again lies in $M_{3}$. Finally, one can compute that \renewcommand{\baselinestretch}{1} \normalsize \begin{eqnarray*} T_{1} \cdot \tau_{1}(m') & = & -q^{-1} \tau_{1}(m'), \end{eqnarray*}   \normalsize and \renewcommand{\baselinestretch}{1} \normalsize \begin{eqnarray*} T_{1} \cdot T_{2}\tau_{1}(m') & = & q(\tau_{1}(m')) + T_{2}\tau_{1}(m').\end{eqnarray*}

Thus $M_{3}$ is a submodule of $M_{2}$. By Theorem $\ref{eq:lemma}$, $M_{3}$ is irreducible, and $M_{2}/M_{3}$ is a 1-dimensional module which is isomorphic to the 1-dimensional module spanned by $m$.

An analoguous argument proves the same result for $N$. Let $Q$ be the 2-dimensional module with $Q = Q_{s_{1}s_{2}t}^{\textrm{gen}}$. Then there is a surjection from $N$ to $Q$, and the kernel of this map, $N_{1}$,  consists of the $t$ and $s_{2}t$ weight spaces of $N$. Then $n = T_{2}T_{1}T_{2}T_{1}T_{2}v - qT_{1}T_{2}T_{1}T_{2}v - T_{2}T_{1}T_{2}v + qT_{1}T_{2}v + T_{2}v - qv$ spans a 1-dimensional submodule of $N_{1}$. Let $N_{2} = N_{1}/n$.

Then $N_{s_{2}t}$ contains a non-zero element $n'$, and $n', \tau_{2}(n')$, and $T_{1}\tau_{2}(n')$ are linearly independent and span $N_{2}$. But $\tau_{2}(n')$ and $T_{1}\tau_{2}(n')$ span a submodule of $N_{2}$, which is irreducible by Theorem $\ref{eq:lemma}$.

(b) Let $\mathbb{C}v_{s_{2}t}$ be the one-dimensional $\widetilde{H}_{\{1\}}$-module with weight $s_{2}t$, and define $L = \widetilde{H} \otimes_{\widetilde{H}_{\{1\}}} \mathbb{C}v_{s_{2}t}$. We claim that the composition factors of $L$ are the same as those of $M$. First, note that the one-dimensional $\widetilde{H}$-module $L_{q,q}$ restricted to $\widetilde{H}_{\{1\}}$ is $\mathbb{C}v_{s_{2}t}$. Then there is a $\widetilde{H}_{\{1\}}$-module map from $\mathbb{C}v_{s_{2}t}$ to $L_{q,q}$, and thus there is a map from $L$ to $L_{q,q}$. Let $L_{1}$ be the kernel of this map. Then $L_{1}$ has a 1-dimensional $s_{2}t$ weight space, and 2-dimensional generalized $t$ and $s_{1}s_{2}t$ weight spaces. Also, $L_{1}$ contains $l = \tau_{2}(v_{s_{2}t}) = T_{2}v_{s_{2}t} - qv_{s_{2}t}$, an element of the $t$ weight space of $L_{1}$. Then we note that $T_{2} \cdot (T_{2}v_{s_{2}t}-qv_{s_{2}t}) = (q-q^{-1})T_{2}v_{s_{2}t} + v_{s_{2}t} - qT_{2}v_{s_{2}t} = q(T_{2}v_{s_{2}t} - qv_{s_{2}t})$, so that $l$ spans a 1-dimensional $\widetilde{H}_{\{2\}}$-submodule of $L_{1}$, with weight $t$.

Thus, there is a $\widetilde{H}_{2}$ map from $\mathbb{C}v_{t}$ to $L_{1}$, and thus a $\widetilde{H}$ map from $M$ to $L$. This map is surjective since $l, T_{2}l, T_{1}T_{2}l, T_{2}T_{1}T_{2}l$, and $T_{1}T_{2}T_{1}T_{2}l$ are linearly independent and span $L_{1}$. Then $L_{1}$ is a quotient of $M$ and its composition factors are composition factors of $M$.

Now, Let $P$ be any irreducible $\widetilde{H}$-module with central character $t_{1,q^{2}}$. If $P$ is not a composition factor of $M$ or $N$, then $P$ must be in the kernel of the (surjective) map from $M(t)$ to $M$. Hence $P$ is at most 6-dimensional, and each of its generalized weight spaces is at most 2-dimensional. If $P = P_{t}^{\textrm{gen}}$, then $P$ is 2-dimensional and must be the module described in Theorem $\ref{2dims}$. Otherwise, we note that $P_{s_{2}t}^{\textrm{gen}} \oplus P_{s_{1}s_{2}t}^{\textrm{gen}}$ is a $\widetilde{H}_{\{1\}}$-submodule of $P$, since the action of $\tau_{1}$ fixes this subspace of $P$. Thus $P_{s_{2}t}^{\textrm{gen}} \oplus P_{s_{1}s_{2}t}^{\textrm{gen}}$ contains an irreducible $\widetilde{H}_{\{1\}}$-submodule. This subspace must be either $P_{s_{1}s_{2}t}^{\textrm{gen}}$ or a 1-dimensional module with weight $s_{2}t$. Hence $P$ is a quotient of either $L$ or $M$ and is isomorphic to a composition factor of $M$.

$\square$

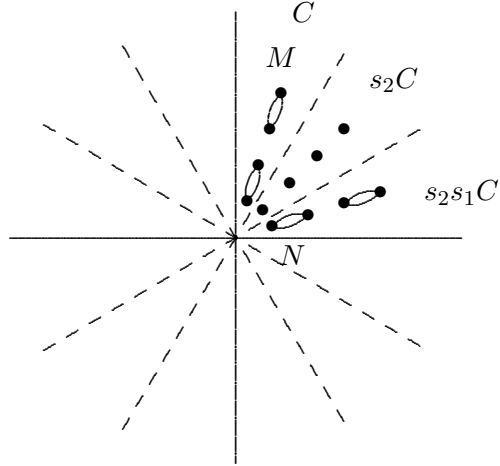
\[\beginpicture
\setcoordinatesystem units <0.3cm,0.3cm>
\setplotarea x from 10 to 25, y from -10 to 10
\linethickness=3pt
\plot 15 -10 15 10 / 
\plot 5 0 25 0 / 
\setdots
\put{$\bullet$} at 15.5 1.6
\put{$\bullet$} at 16 3.2
\put{$\bullet$} at 16.5 4.8
\put{$\bullet$} at 17 6.4
\put{$\bullet$} at 21.4 2
\put{$\bullet$} at 16.6 0.5
\put{$\bullet$} at 18.2 1
\put{$\bullet$} at 19.8 1.5
\put{$\bullet$} at 16.2 1.2
\put{$\bullet$} at 17.4 2.4
\put{$\bullet$} at 18.6 3.6
\put{$\bullet$} at 19.8 4.8
\setdashes
\plot 10 -8.5 20 8.5 / 
\plot 10 8.5 20 -8.5 / 
\plot 6.5 -5 23.5 5 / 
\plot 6.5 5 23.5 -5 / 
\put{$t_{1,q^{2}}, q^{2} = -1$} at 15 -12
\setquadratic
\setsolid
\plot 15.5 1.6 16 2.4 16 3.2 /
\plot 15.5 1.6 15.5 2.4 16 3.2 /
\plot 16.5 4.8 16.5 5.6 17 6.4 /
\plot 16.5 4.8 17 5.6 17 6.4 /
\plot 16.6 0.5 17.4 1 18.2 1 /
\plot 16.6 0.5 17.4 0.5 18.2 1 /
\plot 19.8 1.5 20.6 2 21.4 2 /
\plot 19.8 1.5 20.6 1.5 21.4 2 /
\put{$M$} at 17 8
\put{$N$} at 17.5 -0.8
\put{$C$} at 18 10
\put{$s_{2}C$} at 22 7
\put{$s_{2}s_{1}C$} at 25 2
\endpicture\]

\noindent \textbf{Case 3b: $t_{1, \pm q}$ }

Let $t' \in T$ and assume $\alpha_{1} \in Z(t')$ but $\alpha_{2} \notin P(t')$, so that none of $\alpha_{1} + \alpha_{2}, 2\alpha_{1} + \alpha_{2},$ or $3\alpha_{1} + \alpha_{2}$ are in $P(t')$. Since $P(t') \neq \emptyset$, $3\alpha_{1} + 2\alpha_{2} \in P(t')$, so that $t'(X^{2\alpha_{2}}) = q^{\pm 2}$ and $t'(X^{\alpha_{2}}) = \pm q^{\pm 1}$. By applying $w_{0}$ if necessary, we may assume $t'(X^{\alpha_{2}}) = \pm q$. Thus we will analyze the weights $t_{1,\pm q}$ together, except in one case. If $q$ is a primitive third root of unity then $q = q^{-2}$ and $q^{-1} = q^{2}$, so that $t_{1,q} = t_{1,q^{-2}}$ was analyzed in case 3a. If $q$ is a primitive sixth root of unity then $-q = q^{-2}$ and $-q^{-1} = q^{2}$ so that $t_{1,-q} = t_{1,q^{-2}}$ was analyzed in Case 3a. Thus these cases are excluded from the following analysis by simply assuming that $t'(X^{\alpha_{2}}) \neq q^{-2}$.

If $q^{2} = 1$, then $Z(t') = P(t') = \{\alpha_{1}, 3\alpha_{1} + 2\alpha_{2} \}$, and the irreducibles with central character $t$ can be constructed using the results of $\eqref{clifford}$. Specifically, there are four 3-dimensional modules with central character $t'$.

If $q^{2} \neq 1$, then $s_{1}s_{2}t'(X^{\alpha_{1}}) = t'(X^{-\alpha_{1} -\alpha_{2}}) = \pm q^{\mp 1}$ and $s_{1}s_{2}t'(X^{\alpha_{2}}) = t'(X^{3\alpha_{1} + 2\alpha_{2}}) = q^{\pm 2}$. Then by theorem \ref{PrinOrbit}, $M(t')$ and $M(t)$ have the same composition factors, where $t = s_{1}s_{2}t'$. Also by assuming that $t'(X^{\alpha_{2}}) \neq q^{-2}$, we have $Z(t) = \{ 2\alpha_{1} + \alpha_{2} \}$ and $P(t) = \{\alpha_{2} \}$. Let $w_{0} =s_{1}s_{2}s_{1}s_{2}s_{1}s_{2}$ and define \[\widetilde{H}_{\{2\}} = \mathbb{C}\textrm{-span}\{ T_{2}X^{\lambda}, X^{\lambda} \, | \, \lambda \in P \},\] the subalgebra of $\widetilde{H}$ generated by $T_{2}$ and $\mathbb{C}[X]$. Let $\mathbb{C}v_{t}$ and $\mathbb{C}v_{w_{0}t}$ be the 1-dimensional $\widetilde{H}_{\{2\}}$-modules spanned by $v_{t}$ and $v_{w_{0}t}$, respectively, and given by \[ \begin{matrix} T_{2}v_{t} = qv_{t}, & \quad & X^{\lambda}v_{t} = t(X^{\lambda})v_{t}, \\ T_{2}v_{w_{0}t} = -q^{-1}v_{w_{0}t}, & \textrm{ and } & X^{\lambda}v_{w_{0}t} = w_{0}t(X^{\lambda})v_{w_{0}t}. \end{matrix}\] Then define \[ M = \mathbb{C}v_{t} \otimes_{\widetilde{H}_{\{2\}}} \widetilde{H} \quad \textrm{ and } \] \[ N = \mathbb{C}v_{w_{0}t} \otimes_{\widetilde{H}_{\{2\}}} \widetilde{H}.\]

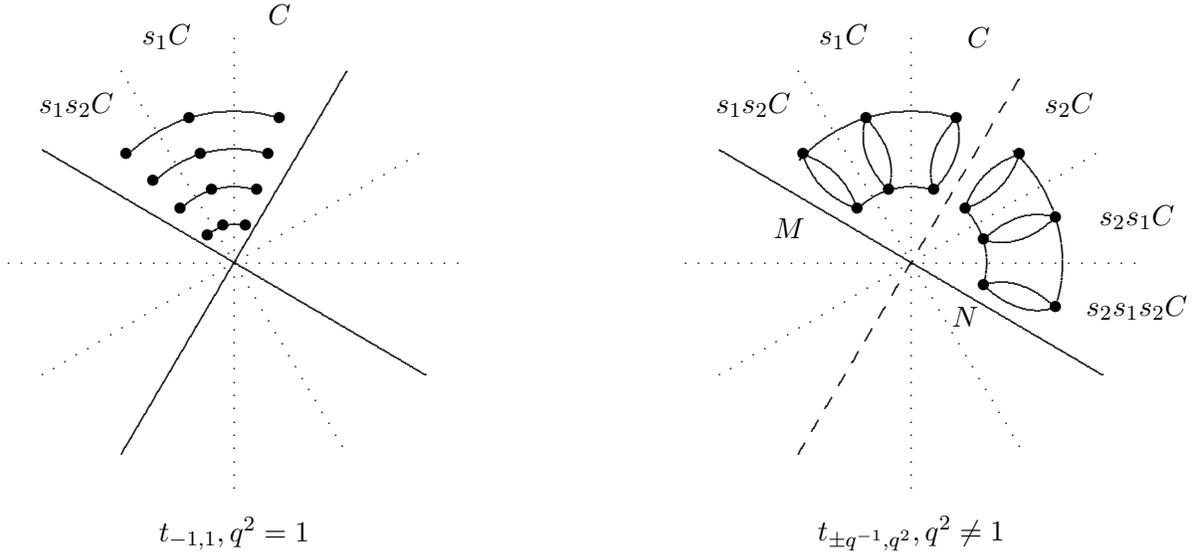
\[\beginpicture
\setcoordinatesystem units <0.3cm,0.3cm>
\setplotarea x from -25 to 25, y from -15 to 15
\linethickness=3pt
\plot 6.5 5 23.5 -5 / 
\plot -23.5 5 -6.5 -5 / 
\plot -20 -8.5 -10 8.5 / 
\setdots
\plot -25 0 -5 0 / 
\plot -23.5 -5 -6.5 5 / 
\plot 10 8.5 20 -8.5 / 
\plot 6.5 -5 23.5 5 / 
\plot -20 8.5 -10 -8.5 / 
\plot -15 -10 -15 10 / 
\plot 15 -10 15 10 / 
\plot 5 0 25 0 / 
\put{$\bullet$} at 16 3.2
\put{$\bullet$} at 14 3.2
\put{$\bullet$} at 17 6.4
\put{$\bullet$} at 13 6.4
\put{$\bullet$} at 21.4 -2
\put{$\bullet$} at 18.2 1
\put{$\bullet$} at 18.2 -1
\put{$\bullet$} at 21.4 2
\put{$\bullet$} at 12.6 2.4
\put{$\bullet$} at 17.4 2.4
\put{$\bullet$} at 10.2 4.8
\put{$\bullet$} at 19.8 4.8
\put{$\bullet$} at -14 3.2
\put{$\bullet$} at -16 3.2
\put{$\bullet$} at -13 6.4
\put{$\bullet$} at -17 6.4
\put{$\bullet$} at -14.5 1.6
\put{$\bullet$} at -13.5 4.8
\put{$\bullet$} at -15.5 1.6
\put{$\bullet$} at -16.5 4.8
\put{$\bullet$} at -17.4 2.4
\put{$\bullet$} at -16.2 1.2
\put{$\bullet$} at -19.8 4.8
\put{$\bullet$} at -18.6 3.6
\setdashes
\plot 10 -8.5 20 8.5 / 
\put{$t_{-1,1}, q^{2} = 1$} at -15 -12
\put{$t_{\pm q^{-1},q^{2}}, q^{2} \neq 1$} at 15 -12
\setquadratic
\setsolid
\circulararc 60 degrees from -14.5 1.6 center at -15 0
\circulararc 60 degrees from -13.5 4.8  center at -15 0
\circulararc 60 degrees from -14 3.2 center at -15 0
\circulararc 60 degrees from -13 6.4 center at -15 0
\circulararc 65 degrees from 21.4 -2 center at 15 0
\circulararc 65 degrees from 18.2 -1 center at 15 0
\circulararc 65 degrees from 16 3.2 center at 15 0
\circulararc 65 degrees from 17 6.4 center at 15 0
\plot 14 3.2 14 4.8 13 6.4 /
\plot 14 3.2 13 4.8 13 6.4 /
\plot 16 3.2 16 4.8 17 6.4 /
\plot 16 3.2 17 4.8 17 6.4 /
\plot 10.2 4.8 11.7 3.9 12.6 2.4 /
\plot 10.2 4.8 11.1 3.3 12.6 2.4 /
\plot 19.8 4.8 18.3 3.9 17.4 2.4 /
\plot 19.8 4.8 18.9 3.3 17.4 2.4 /
\plot 18.2 1 19.8 1 21.4 2 /
\plot 18.2 1 19.8 2 21.4 2 /
\plot 18.2 -1 19.8 -1 21.4 -2 /
\plot 18.2 -1 19.8 -2 21.4 -2 /
\put{$C$} at -13 11
\put{$s_{1}C$} at -18 10
\put{$s_{1}s_{2}C$} at -22 7
\put{$C$} at 18 10
\put{$s_{2}C$} at 22 7
\put{$s_{2}s_{1}C$} at 25 2
\put{$s_{2}s_{1}s_{2}C$} at 25 -2
\put{$s_{1}C$} at 12 10
\put{$s_{1}s_{2}C$} at 8 7
\put{$M$} at 9.6 1.5
\put{$N$} at 17.4 -2.4
\endpicture\]

\begin{proposition} Assume $q^{2} \neq 1$. Let $t' = t_{1,\pm q}$, and define $M$ and $N$ as above. Assume that it is not true that $t'(X^{\alpha_{2}}) = q^{-2}$. Then $M$ and $N$ are irreducible.
\end{proposition}

Let $t = s_{1}s_{2}t$. Under the assumptions, $Z(t) = \{ 2\alpha_{1} + \alpha_{2} \}$ and $P(t) = \{ \alpha_{2}\}$. Then dim $M_{t}^{\textrm{gen}} = $ dim $M_{s_{1}t}^{\textrm{gen}} = $ dim $M_{s_{2}s_{1}t}^{\textrm{gen}} = 2$. By Lemma \ref{eq:lemma}, $M$ has some composition factor $M'$ with dim $M_{s_{2}s_{1}t}^{\textrm{gen}} = 2$, and by Theorem \ref{calibsame}b, $M' = M$ and $M$ is irreducible. Similarly, Lemma \ref{eq:lemma} and Theorem \ref{calibsame}b show that $N$ is irreducible. $\square$

\vspace{.2in}

Under the assumptions of this theorem, since $M$ and $N$ are each 6-dimensional, they must be the only composition factors of $M(t)$. If $t'(X^{\alpha_{2}}) = q^{-2}$, then $w_{0}t' = t_{1,q^{2}}$, which was discussed in the case above.

These two cases are the only weights $t$ with $P(t)$ non-empty and $Z(t)$ containing a short root. Specifically, if $t$ is any weight such that $Z(t)$ contains $\alpha_{1}, \alpha_{1} + \alpha_{2}$, or $2\alpha_{1} + \alpha_{2}$, there exists $w \in W_{0}$ so that $\alpha_{1} \in Z(wt)$. Then $t$ is in the orbit of one of the weights in the previous cases. Then for the following cases, assume $\alpha_{1}, \alpha_{1} + \alpha_{2}, 2\alpha_{1} + \alpha_{2} \notin Z(t)$.

\noindent \textbf{Case 3c: $t_{q^{2},1}$}

If $\alpha_{2} \in Z(t)$ and $\alpha_{1} \in P(t)$, then $\alpha_{1} + \alpha_{2} \in P(t)$ as well, and $t = t_{q^{\pm 2},1}$. These weights are in the same orbit, so we examine $M(t_{q^{2},1})$. If $q^{2} = -1$ then $t(X^{2\alpha_{1} + \alpha_{2}}) = 1$, so that $t$ is in the orbit of one of the weights considered in cases 3a and 3b. If $q^{2} = 1$, then $t = t_{1,1}$ which has also already been considered. Then we assume $q^{2} \neq \pm 1$.

Let $w_{0} =s_{1}s_{2}s_{1}s_{2}s_{1}s_{2}$ and define \[\widetilde{H}_{\{1\}} = \mathbb{C}\textrm{-span}\{ T_{1}X^{\lambda}, X^{\lambda} \, | \, \lambda \in P \},\] the subalgebra of $\widetilde{H}$ generated by $T_{1}$ and $\mathbb{C}[X]$. Let $\mathbb{C}v_{t}$ and $\mathbb{C}v_{w_{0}t}$ be the 1-dimensional $\widetilde{H}_{\{1\}}$-modules spanned by $v_{t}$ and $v_{wt}$, respectively, and given by \[ \begin{matrix} T_{1}v_{t} = qv_{t}, & \quad & X^{\lambda}v_{t} = t(X^{\lambda})v_{t}, \\ T_{1}v_{w_{0}t} = -q^{-1}v_{w_{0}t}, & \textrm{ and } & X^{\lambda}v_{w_{0}t} = w_{0}t(X^{\lambda})v_{w_{0}t}. \end{matrix}\] Then define \[ M = \mathbb{C}v_{t} \otimes_{\widetilde{H}_{\{1\}}} \widetilde{H} \quad \textrm{ and } \] \[ N = \mathbb{C}v_{w_{0}t} \otimes_{\widetilde{H}_{\{1\}}} \widetilde{H}.\]

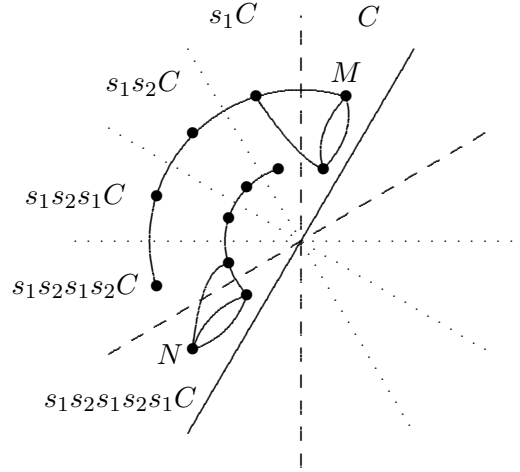
\[\beginpicture
\setcoordinatesystem units <0.3cm,0.3cm>
\setplotarea x from -25 to 25, y from -10 to 10
\linethickness=3pt
\plot -20 8.5 -10 -8.5 / 
\plot -25 0 -5 0 / 
\plot -20 -8.5 -10 8.5 / 
\plot 10 -8.5 20 8.5 / 
\setdots
\plot 6.5 5 23.5 -5 / 
\put{$\bullet$} at 16 3.2
\put{$\bullet$} at 14 3.2
\put{$\bullet$} at 17 6.4
\put{$\bullet$} at 13 6.4
\put{$\bullet$} at 8.6 -2
\put{$\bullet$} at 11.8 1
\put{$\bullet$} at 11.8 -1
\put{$\bullet$} at 8.6 2
\put{$\bullet$} at 12.6 2.4
\put{$\bullet$} at 12.6 -2.4
\put{$\bullet$} at 10.2 4.8
\put{$\bullet$} at 10.2 -4.8
\put{$\bullet$} at -14 3.2
\put{$\bullet$} at -16 3.2
\put{$\bullet$} at -13 6.4
\put{$\bullet$} at -17 6.4
\put{$\bullet$} at -14.5 1.6
\put{$\bullet$} at -13.5 4.8
\put{$\bullet$} at -15.5 1.6
\put{$\bullet$} at -16.5 4.8
\put{$\bullet$} at -17.5 8
\put{$\bullet$} at -18 9.6
\put{$\bullet$} at -12.5 8
\put{$\bullet$} at -12 9.6
\setdashes
\plot 10 8.5 20 -8.5 / 
\plot 6.5 -5 23.5 5 / 
\plot 15 -10 15 10 / 
\plot 5 0 25 0 / 
\plot -15 -10 -15 10 / 
\plot -23.5 5 -6.5 -5 / 
\plot -23.5 -5 -6.5 5 / 
\put{$t_{q^{2},1}, q^{2}$ a primitive third root of unity} at -15 -12
\put{$t_{q^{2},1}, q^{2}$ a primitive fourth root of unity} at 15 -12
\setquadratic
\setsolid
\circulararc 35 degrees from -14 3.2 center at -15 0
\circulararc 35 degrees from -14.5 1.6 center at -15 0
\circulararc 35 degrees from -12.5 8 center at -15 0
\circulararc 35 degrees from -12 9.6 center at -15 0
\circulararc 120 degrees from 14 3.2 center at 15 0
\circulararc 125 degrees from 17 6.4 center at 15 0
\plot -14.5 1.6 -14 2.4 -14 3.2 /
\plot -14.5 1.6 -14.5 2.4 -14 3.2 /
\plot -15.5 1.6 -16 2.4 -16 3.2 /
\plot -15.5 1.6 -15.5 2.4 -16 3.2 /
\plot -16 3.2 -16 4 -16.5 4.8 /
\plot -16 3.2 -16.5 4 -16.5 4.8 /
\plot -16.5 4.8 -17 5.6 -17 6.4 /
\plot -16.5 4.8 -16.5 5.6 -17 6.4 /
\plot -16.5 4.8 -14.7 4.2 -14 3.2 /
\plot -17 6.4 -15 5.2 -14 3.2 /
\plot -17.5 8 -17.5 8.8 -18 9.6 /
\plot -17.5 8 -18 8.8 -18 9.6 /
\plot -12 9.6 -12 8.8 -12.5 8 /
\plot -12 9.6 -12.5 8.8 -12.5 8 /
\plot -12.5 8 -12.5 7.2 -13 6.4 /
\plot -12.5 8 -13 7.2 -13 6.4 /
\plot -13 6.4 -13 5.6 -13.5 4.8 /
\plot -13 6.4 -13.5 5.6 -13.5 4.8 /
\plot -13.5 4.8 -15.5 6.8 -17.5 8 /
\plot -13 6.4 -15 7.4 -17.5 8 /
\plot 16 3.2 16 4.8 17 6.4 /
\plot 16 3.2 17 4.8 17 6.4 /
\plot 10.2 -4.8 11.7 -3.9 12.6 -2.4 /
\plot 10.2 -4.8 11.1 -3.3 12.6 -2.4 /
\plot 16 3.2 14.7 4.2 13 6.4 /
\plot 10.2 -4.8 10.8 -2 11.8 -1 /
\put{$C$} at -13 11
\put{$s_{1}C$} at -18 11
\put{$C$} at 18 10
\put{$s_{1}s_{2}s_{1}s_{2}s_{1}C$} at 7 -7
\put{$s_{1}s_{2}s_{1}C$} at 5 2
\put{$s_{1}s_{2}s_{1}s_{2}C$} at 5 -2
\put{$s_{1}C$} at 12 10
\put{$s_{1}s_{2}C$} at 8 7
\put{$M$} at 17 7.5
\put{$N$} at 9.2 -5
\put{$M$} at -11 10
\put{$N$} at -17.4 2.4
\endpicture\]

\[\beginpicture
\setcoordinatesystem units <0.3cm,0.3cm>
\setplotarea x from -25 to 25, y from -10 to 10
\linethickness=3pt
\plot -5 -8.5 5 8.5 / 
\setdots
\plot -10 0 10 0 / 
\plot -8.5 5 8.5 -5 / 
\plot -5 8.5 5 -8.5 / 
\put{$\bullet$} at 1 3.2
\put{$\bullet$} at -1 3.2
\put{$\bullet$} at 2 6.4
\put{$\bullet$} at -2 6.4
\put{$\bullet$} at -6.4 -2
\put{$\bullet$} at -3.2 1
\put{$\bullet$} at -3.2 -1
\put{$\bullet$} at -6.4 2
\put{$\bullet$} at -2.4 2.4
\put{$\bullet$} at -2.4 -2.4
\put{$\bullet$} at -4.8 4.8
\put{$\bullet$} at -4.8 -4.8
\setdashes
\plot -8.5 -5 8.5 5 / 
\plot 0 -10 0 10 / 
\put{$t_{q^{2},1}, q^{2}$ not a primitive third or fourth root of unity and $q^{2} \neq \pm 1$} at 0 -12
\setquadratic
\setsolid
\circulararc 120 degrees from -1 3.2 center at 0 0
\circulararc 125 degrees from 2 6.4 center at 0 0
\plot 1 3.2 1 4.8 2 6.4 /
\plot 1 3.2 2 4.8 2 6.4 /
\plot -4.8 -4.8 -3.3 -3.9 -2.4 -2.4 /
\plot -4.8 -4.8 -3.9 -3.3 -2.4 -2.4 /
\plot 1 3.2 -0.3 4.2 -2 6.4 /
\plot -4.8 -4.8 -4.2 -2 -3.2 -1 /
\put{$C$} at 3 10
\put{$s_{1}s_{2}s_{1}s_{2}s_{1}C$} at -8 -7
\put{$s_{1}s_{2}s_{1}C$} at -10 2
\put{$s_{1}s_{2}s_{1}s_{2}C$} at -10 -2
\put{$s_{1}C$} at -3 10
\put{$s_{1}s_{2}C$} at -7 7
\put{$M$} at 2 7.5
\put{$N$} at -5.8 -5
\endpicture\]

\begin{proposition}
If $q^{2}$ is a primitive third root of unity, then $M$ and $N$ are irreducible.
\end{proposition}

\noindent \textit{Proof.} If $q^{2}$ is a primitive third root of unity, then $Z(t) = \{ \alpha_{2}, 3\alpha_{1} + \alpha_{2}, 3\alpha_{1} + 2\alpha_{2}\}$ and $P(t) = \{ \alpha_{1}, \alpha_{1} + \alpha_{2}, 2\alpha_{1} + \alpha_{2}\}$.

Then dim $M_{t}^{\textrm{gen}} = 4$ and dim $M_{s_{1}t}^{\textrm{gen}} = 2$. By Lemma \ref{eq:lemma}, if $M' \subseteq M$ is a submodule of $M$, then dim $(M')_{t}^{\textrm{gen}} \geq 2$ and dim $(M')_{s_{1}t}^{\textrm{gen}} \geq 2$. Then dim $(M/M')_{t}^{\textrm{gen}} \leq 2$, but dim $(M/M')_{s_{1}t}^{\textrm{gen}} = 0$, so that Lemma $\ref{eq:lemma}$ implies that $(M/M')_{t} = 0$. Thus $M' = M$ and $M$ is irreducible. Lemma \ref{eq:lemma} similarly implies that $N$ is irreducible. $\square$

Then since $M$ and $N$ have different weight spaces, they are not isomorphic and are the only irreducibles with central character $t$.

\begin{proposition} Assume $q^{2} \neq \pm 1$ and that $q^{2}$ is not a primitive third root of unity.

(a) If $q^{2}$ is a primitive fourth root of unity then $M_{s_{2}s_{1}s_{2}s_{1}t}$ is a 1-dimensional submodule of $M$, and $M'$, the image of the weight spaces $M_{s_{1}s_{2}s_{1}t}$ and $M_{s_{2}s_{1}t}$ in $M/M_{s_{1}s_{2}s_{1}t}$, is an irreducible submodule of $M/M_{s_{1}s_{2}s_{1}t}$. The resulting quotient of $M$ is irreducible.

(b) If $q^{2}$ is a primitive fourth root of unity then $N_{s_{1}t}$ is a 1-dimensional submodule of $N$, and $N'$, the image of the weight spaces $N_{s_{2}s_{1}t}$ and $N_{s_{1}s_{2}s_{1}t}$ in $N/N_{s_{1}t}$, is an irreducible submodule of $N/N_{s_{1}t}$. The resulting quotient of $N$ is irreducible.

(c) The composition factors of $M$ and $N$ are the only composition factors of $M(t)$.

(d) If $q^{2}$ is not a primitive third or fourth root of unity then $M$ and $N$ are irreducible, and are the only irreducible modules with central character $t$.

\end{proposition}

\noindent \textit{Proof.} If $q^{2}$ is not $\pm 1$ or a primitive third root of unity, $Z(t) = \{\alpha_{2}\}$, so that $M$ has one 2-dimensional weight space $M_{t}^{\textrm{gen}}$ and four 1-dimensional weight spaces $M_{s_{1}t}$, $M_{s_{2}s_{1}t}$, $M_{s_{1}s_{2}s_{1}t}$, and $M_{s_{2}s_{1}s_{2}s_{1}t}$.

(a) If $q^{2}$ is a primitive fourth root of unity, then $P(t) = \{ \alpha_{1}, \alpha_{1} + \alpha_{2}, 3\alpha_{1} + \alpha_{2}, 3\alpha_{1} + 2\alpha_{2}\}.$ For $w \in \{ s_{1}, s_{2}s_{1}, s_{1}s_{2}s_{1}, s_{2}s_{1}s_{2}s_{1}\}$, let $m_{wt}$ be a non-zero vector in $M_{wt}$. By a calculation as in \eqref{weightbasis}, \[ m_{wt} = T_{w}T_{2}v_{t} + \sum_{w' < w} a_{w,w'}T_{w'}T_{2}v_{t},\] for $w \in \{ s_{1}, s_{2}s_{1}, s_{1}s_{2}s_{1}, s_{2}s_{1}s_{2}s_{1}\}$, where $a_{w,w'} \in \mathbb{C}$. Then if $s_{i}w > w$, \[ \tau_{i}m_{wt} \neq 0\] for $w \in \{ s_{1}, s_{2}s_{1}, s_{1}s_{2}s_{1}\}$, since the term $T_{i}T_{w}T_{2}$ cannot be canceled by any other term in $\tau_{i}m_{wt}$.

Thus $\tau_{2}: M_{s_{2}s_{1}s_{2}s_{1}t} \rightarrow M_{s_{1}s_{2}s_{1}t}$ is the zero map since, by Theorem \ref{tauthm}, $\tau_{2}^{2}: M_{s_{1}s_{2}s_{1}}t \rightarrow M_{s_{1}s_{2}s_{1}t}$ is the zero map. Hence $M_{s_{2}s_{1}s_{2}s_{1}t}$ is a submodule of $M$. Let $M_{1} = M/M_{s_{2}s_{1}s_{2}s_{1}t}$. Similarly, $\tau_{2}:M_{s_{2}s_{1}t} \rightarrow M_{s_{1}t}$ must be the zero map since, by Theorem \ref{tauthm}, $\tau_{2}^{2}:M_{s_{1}t} \rightarrow M_{s_{1}t}$ is the zero map. Then $M'_{1}$, the subspace spanned by $\overline{m_{s_{2}s_{1}t}}$ and $\overline{m_{s_{1}s_{2}s_{1}t}}$ in $M_{1}$, is a submodule of $M_{1}$. Since $\tau_{1}^{2}: (M'_{1})_{s_{2}s_{1}t} \rightarrow (M'_{1})_{s_{2}s_{1}t}$ is invertible, $M'_{1}$ is irreducible, and Lemma $\ref{eq:lemma}$ shows that $M_{2} = M_{1}/M'_{1}$ is irreducible.

(b) Replacing $t$ by $w_{0}t$ in this argument shows that $N$ also has three composition factors. The weight space $N_{s_{1}t}$ is a submodule of $N$, and $N_{1} = N/N_{s_{1}t}$ has an irreducible 2-dimensional submodule $N'_{1}$, consisting of the image of $N_{s_{2}s_{1}t}$ and $N_{s_{1}s_{2}s_{1}t}$ in $N_{1}$. Lemma \ref{eq:lemma} shows that $N_{1}/N'_{1}$ is irreducible.

(c) The composition factors of $M$ and $N$ are not distinct, since $M'_{1}$ and $N'_{1}$ are irreducible 2-dimensional modules with the same weight spaces, and Theorem \ref{2dims} shows that $M'_{1} \cong N'_{1}$. The 1-dimensional composition factors of $M$ and $N$ are not isomorphic since they have different weights, and the 3-dimensional modules are different since their weight space structures are different.

Counting multiplicities of weight spaces in $M$, $N$, and $M(t)$ shows that the remaining composition factor(s) of $M(t)$ must contain an $s_{2}s_{1}t$ weight space and an $s_{1}s_{2}s_{1}t$ weight space, each of dimension 1. But Theorem \eqref{calibsame}b shows that there must be one remaining composition factor, and Theorem \ref{2dims} shows that it is isomorphic to $M_{1}$. Then the composition factors of $M$ and $N$ are all the composition factors of $M(t)$.

(d) Lemma \ref{eq:lemma} and Theorem \ref{calibsame}b show that both $M$ and $N$ are irreducible if $q^{2}$ is not a primitive third or fourth root of unity. Since $M$ and $N$ are not isomorphic and are each 6-dimensional, they are the only composition factors of $M(t)$. $\square$

\[\beginpicture
\setcoordinatesystem units <0.269cm,0.269cm>
\setplotarea x from -25 to 25, y from -15 to 15
\linethickness=3pt
\plot -20 -8.5 -10 8.5 / 
\setdots
\plot -25 0 -5 0 / 
\plot -23.5 5 -6.5 -5 / 
\plot -20 8.5 -10 -8.5 / 
\put{$\bullet$} at -14 3.2
\put{$\bullet$} at -16 3.2
\put{$\bullet$} at -13 6.4
\put{$\bullet$} at -17 6.4
\put{$\bullet$} at -21.4 -2
\put{$\bullet$} at -18.2 1
\put{$\bullet$} at -18.2 -1
\put{$\bullet$} at -21.4 2
\put{$\bullet$} at -17.4 2.4
\put{$\bullet$} at -17.4 -2.4
\put{$\bullet$} at -19.8 4.8
\put{$\bullet$} at -19.8 -4.8
\setdashes
\plot -23.5 -5 -6.5 5 / 
\plot -15 -10 -15 10 / 
\put{$t_{q^{2},1}, q^{2}$ not a primitive third or fourth root of unity} at -15 -12
\setquadratic
\setsolid
\circulararc 120 degrees from -16 3.2 center at -15 0
\circulararc 125 degrees from -13 6.4 center at -15 0
\plot -14 3.2 -14 4.8 -13 6.4 /
\plot -14 3.2 -13 4.8 -13 6.4 /
\plot -19.8 -4.8 -18.3 -3.9 -17.4 -2.4 /
\plot -19.8 -4.8 -18.9 -3.3 -17.4 -2.4 /
\plot -14 3.2 -15.3 4.2 -17 6.4 /
\plot -19.8 -4.8 -19.2 -2 -18.2 -1 /
\put{$C$} at -12 10
\put{$s_{1}s_{2}s_{1}s_{2}s_{1}C$} at -23 -7
\put{$s_{1}s_{2}s_{1}C$} at -25 2
\put{$s_{1}s_{2}s_{1}s_{2}C$} at -25 -2
\put{$s_{1}C$} at -18 10
\put{$s_{1}s_{2}C$} at -22 7
\put{$M$} at -13 7.5
\put{$N$} at -20.8 -5
\setlinear
\setsolid
\plot 10 -8.5 20 8.5 / 
\setdots
\plot 6.5 5 23.5 -5 / 
\put{$\bullet$} at 16 3.2
\put{$\bullet$} at 14 3.2
\put{$\bullet$} at 17 6.4
\put{$\bullet$} at 13 6.4
\put{$\bullet$} at 8.6 -2
\put{$\bullet$} at 11.8 1
\put{$\bullet$} at 11.8 -1
\put{$\bullet$} at 8.6 2
\put{$\bullet$} at 12.6 2.4
\put{$\bullet$} at 12.6 -2.4
\put{$\bullet$} at 10.2 4.8
\put{$\bullet$} at 10.2 -4.8
\setdashes
\plot 10 8.5 20 -8.5 / 
\plot 6.5 -5 23.5 5 / 
\plot 15 -10 15 10 / 
\plot 5 0 25 0 / 
\put{$t_{q^{2},1}, q^{2}$ a primitive fourth root of unity} at 15 -12
\setquadratic
\setsolid
\circulararc 30 degrees from 11.8 -1 center at 15 0
\circulararc 35 degrees from 17 6.4 center at 15 0
\circulararc 30 degrees from 10.2 4.8 center at 15 0
\circulararc 30 degrees from 12.6 2.4 center at 15 0
\plot 16 3.2 16 4.8 17 6.4 /
\plot 16 3.2 17 4.8 17 6.4 /
\plot 10.2 -4.8 11.7 -3.9 12.6 -2.4 /
\plot 10.2 -4.8 11.1 -3.3 12.6 -2.4 /
\plot 16 3.2 14.7 4.2 13 6.4 /
\plot 10.2 -4.8 10.8 -2 11.8 -1 /
\put{$C$} at 18 10
\put{$s_{1}s_{2}s_{1}s_{2}s_{1}C$} at 7 -7
\put{$s_{1}s_{2}s_{1}C$} at 5 2
\put{$s_{1}s_{2}s_{1}s_{2}C$} at 5 -2
\put{$s_{1}C$} at 12 10
\put{$s_{1}s_{2}C$} at 8 7
\endpicture\]

\noindent \textbf{Case 3d: $t_{\pm q, 1}$}

If $\alpha_{2} \in Z(t')$ and $2\alpha_{1} + \alpha_{2} \in P(t')$, then $t'(X^{2\alpha_{1}}) = q^{ \pm 2}$ and $t'(X^{\alpha_{1}}) = \pm q^{\pm 1}$. By replacing $t'$ by $w_{0}t'$ if necessary, it suffices to assume that $t'(X^{\alpha_{1}}) = \pm q$. If $t'(X^{\alpha_{1}}) = q^{-2}$, then $t'$ was analyzed in case 3c. This occurs when $q^{3} = 1$ and $t' = t_{q,1}$, or when $q^{3} = -1$ and $t'=t_{-q,1}$. Thus the following analysis will apply to $t_{q,1}$ except if $q^{3} = 1$, and $t_{-q,1}$ except for when $q^{3} = -1$. (This is tantamount to assuming that $P(t)$ and $Z(t)$ each contain exactly one element for this $t$.)

Also, if $t'(X^{3\alpha_{1}}) = q^{-2}$, then $P(t)$ also contains $3\alpha_{1}+\alpha_{2}$ and $3\alpha_{1}+2\alpha_{2}$. This occurs when $q^{5} = 1$ and $t'(X^{\alpha_{1}}) = q$ or when $q^{5} = -1$ and $t'(X^{\alpha_{1}}) = -q$. When either of these hold, $t'$ is the same orbit as $t_{q^{2},q^{2}}$. This case (which was specifically not addressed in Case 2 above) will be treated separately below.

Let $w_{0} =s_{1}s_{2}s_{1}s_{2}s_{1}s_{2}$ and define \[\widetilde{H}_{\{1\}} = \mathbb{C}\textrm{-span}\{ T_{1}X^{\lambda}, X^{\lambda} \, | \, \lambda \in P \},\] the subalgebra of $\widetilde{H}$ generated by $T_{1}$ and $\mathbb{C}[X]$. Define $t = s_{2}s_{1}t'$ so that $t(X^{\alpha_{1}}) = q^{2}$ and $t(X^{\alpha_{2}}) = \pm q^{-3}$. Let $\mathbb{C}v_{t}$ and $\mathbb{C}v_{w_{0}t}$ be the 1-dimensional $\widetilde{H}_{\{1\}}$-modules spanned by $v_{t}$ and $v_{wt}$, respectively, and given by \[ \begin{matrix} T_{1}v_{t} = qv_{t}, & \quad & X^{\lambda}v_{t} = t(X^{\lambda})v_{t}, \\ T_{1}v_{w_{0}t} = -q^{-1}v_{w_{0}t}, & \textrm{ and } & X^{\lambda}v_{w_{0}t} = w_{0}t(X^{\lambda})v_{w_{0}t}. \end{matrix}\] Then define \[ M = \mathbb{C}v_{t} \otimes_{\widetilde{H}_{\{1\}}} \widetilde{H} \quad \textrm{ and } \] \[ N = \mathbb{C}v_{w_{0}t} \otimes_{\widetilde{H}_{\{1\}}} \widetilde{H}.\]

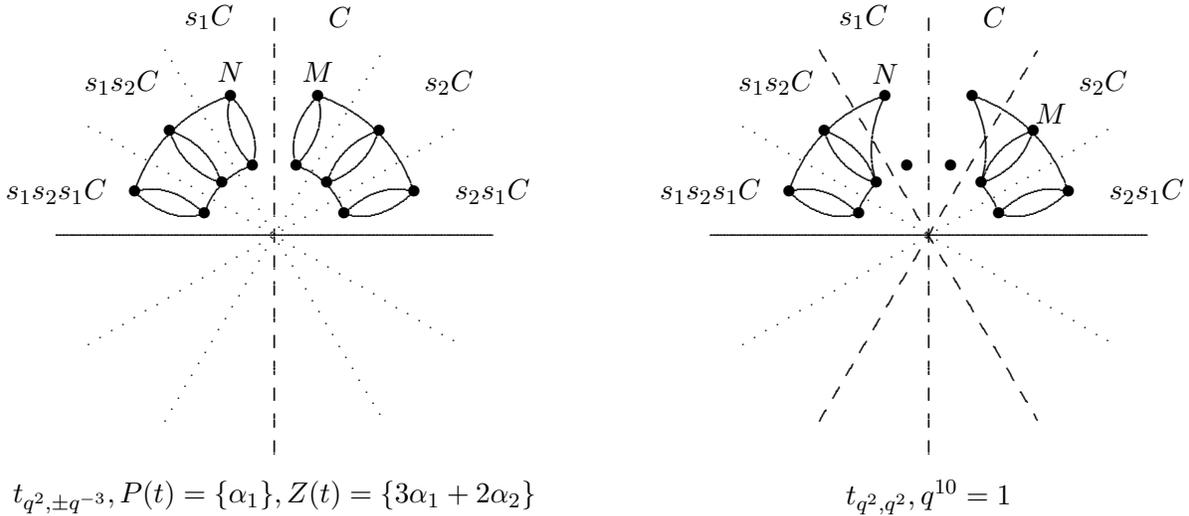
\[\beginpicture
\setcoordinatesystem units <0.29cm,0.29cm>
\setplotarea x from -25 to 25, y from -15 to 15
\linethickness=3pt
\plot -25 0 -5 0 / 
\setdots
\plot -20 8.5 -10 -8.5 / 
\plot -23.5 5 -6.5 -5 / 
\plot -23.5 -5 -6.5 5 / 
\plot -20 -8.5 -10 8.5 / 
\put{$\bullet$} at -14 3.2
\put{$\bullet$} at -16 3.2
\put{$\bullet$} at -13 6.4
\put{$\bullet$} at -17 6.4
\put{$\bullet$} at -8.6 2
\put{$\bullet$} at -18.2 1
\put{$\bullet$} at -11.8 1
\put{$\bullet$} at -21.4 2
\put{$\bullet$} at -17.4 2.4
\put{$\bullet$} at -12.6 2.4
\put{$\bullet$} at -19.8 4.8
\put{$\bullet$} at -10.2 4.8
\setdashes
\plot -15 -10 -15 10 / 
\put{$t_{q^{2},\pm q^{- 3}}, P(t) = \{ \alpha_{1} \} , Z(t) = \{ 3\alpha_{1} + 2\alpha_{2} \}$} at -15 -12
\setquadratic
\setsolid
\circulararc 55 degrees from -11.8 1 center at -15 0
\circulararc 55 degrees from -8.6 2 center at -15 0
\circulararc 55 degrees from -16 3.2 center at -15 0
\circulararc 55 degrees from -17 6.4 center at -15 0
\plot -14 3.2 -14 4.8 -13 6.4 /
\plot -14 3.2 -13 4.8 -13 6.4 /
\plot -16 3.2 -16 4.8 -17 6.4 /
\plot -16 3.2 -17 4.8 -17 6.4 /
\plot -10.2 4.8 -11.7 3.9 -12.6 2.4 /
\plot -10.2 4.8 -11.1 3.3 -12.6 2.4 /
\plot -19.8 4.8 -18.3 3.9 -17.4 2.4 /
\plot -19.8 4.8 -18.9 3.3 -17.4 2.4 /
\plot -18.2 1 -19.8 1 -21.4 2 /
\plot -18.2 1 -19.8 2 -21.4 2 /
\plot -11.8 1 -10.2 1 -8.6 2 /
\plot -11.8 1 -10.2 2 -8.6 2 /
\put{$C$} at -12 10
\put{$s_{2}s_{1}C$} at -5 2
\put{$s_{1}s_{2}s_{1}C$} at -25 2
\put{$s_{2}C$} at -7 7
\put{$s_{1}C$} at -18 10
\put{$s_{1}s_{2}C$} at -22 7
\put{$M$} at -13 7.5
\put{$N$} at -17 7.5
\setlinear
\setsolid
\plot 5 0 25 0 / 
\setdots
\plot 6.5 5 23.5 -5 / 
\plot 6.5 -5 23.5 5 / 
\put{$\bullet$} at 16 3.2
\put{$\bullet$} at 14 3.2
\put{$\bullet$} at 17 6.4
\put{$\bullet$} at 13 6.4
\put{$\bullet$} at 21.4 2
\put{$\bullet$} at 11.8 1
\put{$\bullet$} at 18.2 1
\put{$\bullet$} at 8.6 2
\put{$\bullet$} at 12.6 2.4
\put{$\bullet$} at 17.4 2.4
\put{$\bullet$} at 10.2 4.8
\put{$\bullet$} at 19.8 4.8
\setdashes
\plot 10 8.5 20 -8.5 / 
\plot 10 -8.5 20 8.5 / 
\plot 15 -10 15 10 / 
\put{$t_{q^{2},q^{2}}, q^{10} = 1$} at 15 -12
\setquadratic
\setsolid
\circulararc 55 degrees from 21.4 2 center at 15 0
\circulararc 30 degrees from 18.2 1 center at 15 0
\circulararc 30 degrees from 12.6 2.4 center at 15 0
\circulararc 55 degrees from 13 6.4 center at 15 0
\plot 12.6 2.4 12.4 4.4 13 6.4 /
\plot 17.4 2.4 17.6 4.4 17 6.4 /
\plot 10.2 4.8 11.7 3.9 12.6 2.4 /
\plot 10.2 4.8 11.1 3.3 12.6 2.4 /
\plot 11.8 1 10.2 1 8.6 2 /
\plot 11.8 1 10.2 2 8.6 2 /
\plot 19.8 4.8 18.3 3.9 17.4 2.4 /
\plot 19.8 4.8 18.9 3.3 17.4 2.4 /
\plot 18.2 1 19.8 1 21.4 2 /
\plot 18.2 1 19.8 2 21.4 2 /
\put{$C$} at 18 10
\put{$s_{2}s_{1}C$} at 25 2
\put{$s_{1}s_{2}s_{1}C$} at 5 2
\put{$s_{2}C$} at 23 7
\put{$s_{1}C$} at 12 10
\put{$s_{1}s_{2}C$} at 8 7
\put{$M$} at 20.6 5.6
\put{$N$} at 13 7.4
\endpicture\]

\begin{proposition} Let $t' = t_{\pm q,1}$. Assume that $t'(X^{\alpha_{1}}) \neq q^{-2}$ and $t'(X^{3\alpha_{1}}) \neq q^{-2}$. Then $M$ and $N$ are irreducible.

\end{proposition}

\noindent \textit{Proof.} Let $t = s_{2}s_{1}t'$. Under the assumptions, $Z(t) = \{3\alpha_{1} + 2\alpha_{2}\}$ and $P(t) = \{ \alpha_{1} \}$. Then dim $M_{t}^{\textrm{gen}} = $ dim $M_{s_{2}t}^{\textrm{gen}} = $ dim $M_{s_{1}s_{2}t}^{\textrm{gen}} = 2$. Then Lemma \ref{eq:lemma} and Theorem \ref{calibsame}b show that $M$ is irreducible. $N$ is also irreducible by the same reasoning. $\square$

\vspace{.1in}

Under the assumption of the theorem, since $M$ and $N$ are not isomorphic and are each 6-dimensional, they are the only composition factors of $M(t)$. Note that $t'(X^{\alpha_{1}}) = q^{-2}$ exactly if $q^{3} =  1$ or $-1$ and $t'(X^{\alpha_{1}}) = q$ or $-q$, respectively. In this case, the central character $t'$ has been analyzed above (Case 3c)
Also, $t'(X^{\alpha_{1}}) = q^{-4}$ exactly if $q^{5} = 1$ or $-1$ and $t'(X^{\alpha_{1}}) = q$ or $-q$, respectively. In this case, $t'$ is in the same orbit as $t_{q^{2},q^{2}}$.
\begin{proposition}
If $t = t_{q^{2},q^{2}}$ and $q^{2}$ is a fifth root of unity, then

(a) $M$ has a 5-dimensional irreducible submodule $M'$ and

(b) $N$ has a 5-dimensional irreducible submodule $N'$.

\end{proposition}

\noindent \textit{Proof.} Given these assumptions, $Z(t) = \{ 3\alpha_{1} + 2\alpha_{2} \}$ and $P(t) = \{\alpha_{1}, \alpha_{2}, 3\alpha_{1} + \alpha_{2}\}$. Then dim $M_{t}^{\textrm{gen}} = $ dim $M_{s_{2}t}^{\textrm{gen}} = $ dim $M_{s_{1}s_{2}t}^{\textrm{gen}} = 2$. Let $L_{q,q} = \mathbb{C}v$ be the 1-dimensional $\widetilde{H}$ module given by \[ T_{i}v = qv, \quad X^{\alpha_{i}} = q^{2}v, \quad \textrm{ for }  i =1,2.\] Since \[ \textrm{Hom}_{\widetilde{H}}(\widetilde{H} \otimes_{\widetilde{H}_{\{1\}}} \mathbb{C}v_{t}, L_{q,q}) = \textrm{Hom}_{\widetilde{H}_{\{1\}}}(\mathbb{C}v_{t}, L_{q,q}|_{\widetilde{H}_{\{1\}}})\] and
\renewcommand{\baselinestretch}{1} \normalsize \[ \begin{array}{rcl} \phi: \mathbb{C}v_{t} & \rightarrow & L_{q,q} \\  v_{t} & \mapsto & v\end{array}\]   \normalsize
is a map of $\widetilde{H}_{\{1\}}$-modules, there is a nonzero map $\theta: M \rightarrow L_{q,q}$. Then let $M_{1}$ be the kernel of $\theta$, which is 5-dimensional. Similarly, there is a map $\rho: N \rightarrow L_{q^{-1},q^{-1}}$, where $L_{q^{-1},q^{-1}} = \mathbb{C}v$ is given by \[ T_{i} = -q^{-1}v, X^{\alpha_{i}} = q^{-2}v, \quad \textrm{ for } i=1,2.\] Then if $N_{1}$ is the 5-dimensional kernel of $\rho$, Lemma \ref{eq:lemma} and Theorem \ref{calibsame}b show that $M_{1}$ and $N_{1}$ are both irreducible. $\square$

\vspace{.2in}

\noindent \textbf{Case 3e: $t_{q^{2/3},1}$}

If $\alpha_{2} \in Z(t)$ and $3\alpha_{1} + \alpha_{2} \in P(t)$, then $3\alpha_{1} + 2\alpha_{2} \in P(t)$ as well. If $t(X^{3\alpha_{1} + \alpha_{2}}) = q^{-2}$, then $w_{0}t(X^{\alpha_{2}}) = 1$ and $w_{0}t(X^{3\alpha_{1} + \alpha_{2}}) = q^{2}$, so by replacing $t$ with $w_{0}t$ if necessary, assume that $t(X^{\alpha_{1}})^{3} = q^{2}$. If $\alpha_{1} \in P(t)$, then this weight was analyzed in case 3c, and if $2\alpha_{1} + \alpha_{2} \in P(t)$, then this weight was analyzed in case 3d.

Then we assume $Z(t) = \{ \alpha_{2} \}$ and $P(t) = \{ 3\alpha_{1} + \alpha_{2}, 3\alpha_{1} + 2\alpha_{2}\}$. Let $t' = s_{1}t$ so that $Z(t') = \{ 3\alpha_{1} + \alpha_{2}\}$ and $P(t') = \{ \alpha_{2}, 3\alpha_{1}+2\alpha_{2}\}$. Let $w_{0} =s_{1}s_{2}s_{1}s_{2}s_{1}s_{2}$ and define \[\widetilde{H}_{\{1\}} = \mathbb{C}\textrm{-span}\{ T_{1}X^{\lambda}, X^{\lambda} \, | \, \lambda \in P \},\] the subalgebra of $\widetilde{H}$ generated by $T_{1}$ and $\mathbb{C}[X]$. Let $\mathbb{C}v_{t}$ and $\mathbb{C}v_{w_{0}t'}$ be the 1-dimensional $\widetilde{H}_{\{1\}}$-modules spanned by $v_{t'}$ and $v_{w_{0}t'}$, respectively, and given by \[ \begin{matrix} T_{1}v_{t'} = qv_{t'}, & \quad & X^{\lambda}v_{t'} = t'(X^{\lambda})v_{t'}, \\ T_{1}v_{w_{0}t'} = -q^{-1}v_{w_{0}t'}, & \textrm{ and } & X^{\lambda}v_{w_{0}t'} = w_{0}t'(X^{\lambda})v_{w_{0}t'}. \end{matrix}\] Then define \[ M = \mathbb{C}v_{t'} \otimes_{\widetilde{H}_{\{1\}}} \widetilde{H} \quad \textrm{ and } \] \[ N = \mathbb{C}v_{w_{0}t'} \otimes_{\widetilde{H}_{\{1\}}} \widetilde{H}.\]

\[\beginpicture
\setcoordinatesystem units <0.3cm,0.3cm>
\setplotarea x from -10 to 30, y from -15 to 15
\linethickness=3pt
\plot -10 -8.5 0 8.5 / 
\setdots
\plot -13.5 -5 3.5 5 / 
\plot -5 10 -5 -10 / 
\plot -13.5 5 3.5 -5 / 
\put{$\bullet$} at -4 3.2
\put{$\bullet$} at -6 3.2
\put{$\bullet$} at -3 6.4
\put{$\bullet$} at -7 6.4
\put{$\bullet$} at -11.4 -2
\put{$\bullet$} at -8.2 1
\put{$\bullet$} at -8.2 -1
\put{$\bullet$} at -11.4 2
\put{$\bullet$} at -7.4 2.4
\put{$\bullet$} at -7.4 -2.4
\put{$\bullet$} at -9.8 4.8
\put{$\bullet$} at -9.8 -4.8
\setdashes
\plot -10 8.5 0 -8.5 / 
\plot -15 0 5 0 / 
\put{$t_{q^{2/3},1}, P(t) = \{ 3\alpha_{1} + \alpha_{2}, 3\alpha_{1}+2\alpha_{2} \}, q^{2} \neq -1 $} at -7 -12
\setquadratic
\setsolid
\circulararc 30 degrees from -8.2 -1 center at -5 0
\circulararc 35 degrees from -3 6.4 center at -5 0
\circulararc 90 degrees from -4 3.2 center at -5 0
\circulararc 90 degrees from -9.8 4.8 center at -5 0
\plot -7 6.4 -7.5 4.4 -7.4 2.4 /
\plot -11.4 2 -9.5 0.8 -8.2 -1 /
\plot -6 3.2 -6 4.8 -7 6.4 /
\plot -6 3.2 -7 4.8 -7 6.4 /
\plot -4 3.2 -4 4.8 -3 6.4 /
\plot -4 3.2 -3 4.8 -3 6.4 /
\plot -9.8 -4.8 -8.3 -3.9 -7.4 -2.4 /
\plot -9.8 -4.8 -8.9 -3.3 -7.4 -2.4 /
\plot -8.2 -1 -9.8 -1 -11.4 -2 /
\plot -8.2 -1 -9.8 -2 -11.4 -2 /
\put{$C$} at -2 10
\put{$s_{1}s_{2}s_{1}s_{2}s_{1}C$} at -13 -7
\put{$s_{1}s_{2}s_{1}C$} at -15 2
\put{$s_{1}s_{2}s_{1}s_{2}C$} at -15 -2
\put{$s_{1}C$} at -8 10
\put{$s_{1}s_{2}C$} at -12 7
\setlinear
\setsolid
\plot 15 -8.5 25 8.5 / 
\setdots
\plot 11.5 -5 28.5 5 / 
\plot 20 10 20 -10 / 
\plot 11.5 5 28.5 -5 / 
\put{$\bullet$} at 21 3.2
\put{$\bullet$} at 19 3.2
\put{$\bullet$} at 22 6.4
\put{$\bullet$} at 18 6.4
\put{$\bullet$} at 13.6 -2
\put{$\bullet$} at 16.8 1
\put{$\bullet$} at 16.8 -1
\put{$\bullet$} at 13.6 2
\put{$\bullet$} at 17.6 2.4
\put{$\bullet$} at 17.6 -2.4
\put{$\bullet$} at 15.2 4.8
\put{$\bullet$} at 15.2 -4.8
\setdashes
\plot 15 8.5 25 -8.5 / 
\plot 10 0 30 0 / 
\put{$t_{q^{2/3},1}, P(t) = \{ 3\alpha_{1} + \alpha_{2}, 3\alpha_{1}+2\alpha_{2} \}, q^{2} = -1$} at 22 -12
\setquadratic
\setsolid
\circulararc 30 degrees from 16.8 -1 center at 20 0
\circulararc 35 degrees from 22 6.4 center at 20 0
\circulararc 30 degrees from 21 3.2 center at 20 0
\circulararc 30 degrees from 17.6 2.4 center at 20 0
\circulararc 30 degrees from 15.2 4.8 center at 20 0
\circulararc 28 degrees from 13.6 -2 center at 20 0
\plot 19 3.2 19 4.8 18 6.4 /
\plot 19 3.2 18 4.8 18 6.4 /
\plot 21 3.2 21 4.8 22 6.4 /
\plot 21 3.2 22 4.8 22 6.4 /
\plot 15.2 -4.8 16.7 -3.9 17.6 -2.4 /
\plot 15.2 -4.8 16.1 -3.3 17.6 -2.4 /
\plot 16.8 -1 15.2 -1 13.6 -2 /
\plot 16.8 -1 15.2 -2 13.6 -2 /
\put{$C$} at 23 10
\put{$s_{1}s_{2}s_{1}s_{2}s_{1}C$} at 12 -7
\put{$s_{1}s_{2}s_{1}C$} at 10 2
\put{$s_{1}s_{2}s_{1}s_{2}C$} at 10 -2
\put{$s_{1}C$} at 17 10
\put{$s_{1}s_{2}C$} at 13 7
\endpicture\]

\begin{proposition} Assume $t = t_{q^{2/3},1}$, where $q^{2/3}$ is a third root of $q^{2}$ not equal to $q^{\pm 2}$ or $\pm q^{\pm 1}$, and that $q^{2} \neq \pm 1$. Then $M$ and $N$ are irreducible.

\end{proposition}

\noindent \textit{Proof.} Under the assumptions, $Z(t') = \{ 3\alpha_{1} + \alpha_{2}\}$ and $P(t) = \{\alpha_{2}, 3\alpha_{1} + 2\alpha_{2} \}$, so that dim $M_{t}^{\textrm{gen}} = $ dim $M_{s_{1}t}^{\textrm{gen}} = 2$ while dim $M_{s_{2}s_{1}} =$ dim $M_{s_{1}s_{2}s_{1}t} = 1$. Then Lemma \ref{eq:lemma} and Theorem \ref{calibsame}b show that $M$ is irreducible. Similarly, $N$ is irreducible by the same reasoning, so that $M$ and $N$ are the only irreducible modules with central character $t$. $\square$

\begin{proposition} Assume $q^{2} = -1$ and that $t = t_{q^{2/3},1}$, where $q^{2/3}$ is a third root of $q^{2}$ not equal to $q^{\pm 2}$ or $\pm q^{\pm 1}$. Then $M$ and $N$ each have an irreducible 2-dimensional submodule consisting of their $s_{2}s_{1}t$ and $s_{1}s_{2}s_{1}t$ weight spaces. The resulting quotients are irreducible.

\end{proposition}

\noindent \textit{Proof.} Let $X = \{e, s_{1},s_{2}s_{1}, s_{1}s_{2}s_{1}\}$. By a calculation analogous to that in $\eqref{weightbasis}$, the generalized $t'$ weight space of $M$ is generated by $v_{t'}$ and a vector $v$ of the form $ \sum_{x \in X} a_{x} T_{x}v_{t'}$, where the $a_{x}$ are in $\mathbb{C}$ and $a_{s_{1}s_{2}s_{1}} \neq 0$. Then $\tau_{2}(v) \neq 0$, since it contains a non-zero $T_{2}T_{1}T_{2}T_{1}v_{t}$ term. But then $\tau_{2}\tau_{2}(v) = 0$ by Theorem $\ref{tauthm}$c, so that the space $M_{1} = M_{s_{2}t'} \oplus M_{s_{1}s_{2}t'}$ is actually a submodule of $M$. The resulting quotient $M/M_{1}$ is irreducible by Theorem $\ref{calibsame}$. A similar argument shows the same for $N$. $\square$

Note that Theorem $\ref{2dims}$ shows that the 2-dimensional composition factors of $M$ and $N$ are isomorphic, and this proposition implies that when $q^{2} = -1$, the composition factors of $M$ and $N$ are the only irreducibles with this central character. Counting dimensions of the weight spaces of these irreducibles shows that the final composition factor of $M(t)$ must be 2-dimensional with weights $s_{2}t'$ and $s_{1}s_{2}t'$, since there are no 1-dimensional modules with this central character. So the last composition factor of $M(t)$ must also be isomorphic to the 2-dimensional submodule of $M$.

\vspace{.2in}

\noindent \textbf{Summary.} We summarize the results of the previous theorems, including our choices of representatives for the various central characters, in the following tables. It should be noted that for any value of $q$ with $q^{2}$ not a root of unity of order 6 or less, the representation theory of $\widetilde{H}$ can be described in terms of $q$ only. If $q^{2}$ is a primitive root of unity of order 6 or less, then the representation theory of $\widetilde{H}$ does not fit that same description. This fact can be seen through a number of different lenses. It is a reflection of the fact that the sets $P(t)$ and $Z(t)$ for all possible central characters $t$ can be described solely in terms of $q$. In the local region pictures, this is reflected in the fact that the hyperplanes $H_{\alpha}$ and $H_{\alpha \pm \delta}$ are distinct \textit{unless} $q^{2}$ is a root of unity of order 6 or less. When these hyperplanes coincide, the sets $P(t)$ and $Z(t)$ change for characters on those hyperplanes.

Some notes are necessary about the table below. An entry of ``N/A'' means that the given central character is in the same orbit as a previous character for that particular value of $q$, as described after Proposition \ref{G2charsgen}.

If $q^{10} = 1$, we are assuming that $q^{5} = -1$, so that  $t_{\pm q,1} = t_{\pm q^{-4},1}$.

If $q^{8} = 1$, then only the central characters $t_{q^{2},1}$, $t_{q^{2},-q^{-2}}$, and $t_{q^{2},q^{2}}$ change from the generic case. All three of these characters are now in the same orbit. Also, we assume that for the central character $t_{q^{2/3},1}$, we choose a cube root of $q^{2}$ besides $q^{-2}$.

If $q^{6} = 1$, the entries for the central characters $t_{\pm q,1}$ and $t_{1,\pm q}$ only apply to the characters $t_{ - q^{-2},1}$ and $t_{1, -q^{-2}}$ (depending on whether $q^{3}$ is 1 or $-1$.) Then we note that $t_{1^{1/3},1} = t_{q^{\pm 2},1}$, and $t_{q^{-2},1}$ is in the same orbit as $t_{q^{2},1}$. Also, $t_{1,q^{-2}} = w_{0}t_{1,q^{2}}$, and $t_{q^{-2},1} = w_{0}t_{q^{2},1}$. Finally, $t_{1^{1/3},q^{2}} = t_{q^{\pm 2},q^{2}}$, but $s_{1}t_{q^{2},q^{2}} = t_{q^{-2},q^{2}} = s_{2}t_{1,q^{-2}}$ and so both are in the same orbit as $t_{1,q^{2}}$.

When $q^{2} = -1$, a number of characters change from the general case. Now, $t_{1,-1} = t_{1,q^{2}}$, which is in the same orbit as $t_{q^{2},-q^{-2}}$, $t_{q^{2},q^{2}}$ and $t_{q^{2},1}$. Similarly, $t_{1^{1/3},q^{2}}$ is in the same orbit as $t_{q^{2/3},1}$.

When $q^{2} = 1$, $Z(t) = P(t)$ for all $t \in T$.

\noindent \textbf{Acknowledgments.} The author would like to thank Arun Ram for his continued guidance and many helpful conversations regarding this work.

\begin{table}[htpb]
\centering
$\begin{array}{| c | c | c | c | c | c | c | c |}
\hline \, & \multicolumn{7}{|c|}{\textrm{Dims. of Irreds.}} \\
\hline t & q \textrm{ generic}& q^{12} = 1 & q^{10} = 1 & q^{8} = 1 & q^{6} = 1 & q^{4} = 1 & q = -1\\
\hline t_{1,1} & 12 & 12 & 12 & 12 & 12 & 12 & 1,1,1,1,2,2 \\
\hline t_{1,-1} & 12 & 12 & 12 & 12 & 12 & 1,2,2 & 3,3,3,3 \\
\hline t_{1^{1/3},1} & 12 & 12 & 12 & 12 & \textrm{N/A} & 12 & 3,3,6 \\
\hline t_{1,q^{2}} & 1,1,2,3,3 & 1,1,2,3,3 & 1,1,2,3,3 & 1,1,2,3,3 & 1,1,1,1,3,3 & \textrm{N/A} & \textrm{N/A} \\
\hline t_{1,\pm q} & 6,6 & 6,6 & 6,6 & 6,6 & 6,6 & 6,6 & \textrm{N/A} \\
\hline t_{1,z} & 12 & 12 & 12 & 12 & 12 & 12 & 6,6 \\
\hline t_{q^{2},1} & 6,6 & 6,6 & 6,6 & 1,1,2,3,3 & 6,6 & \textrm{N/A} & \textrm{N/A} \\
\hline t_{q,1} & 6,6 & 6,6 & 1,1,5,5 & 6,6 & 6,6 & 6,6 & \textrm{N/A} \\
\hline t_{- q,1} & 6,6 & 6,6 & 6,6 & 6,6 & 6,6 & \textrm{N/A} & \textrm{N/A} \\
\hline t_{q^{2/3},1} & 6,6 & 6,6 & 6,6 & 6,6 & 6,6 & 2,4,4 & \textrm{N/A} \\
\hline t_{z,1} & 12 & 12 & 12 & 12 & 12 & 12 & 6,6 \\
\hline t_{1^{1/3},q^{2}} & 3,3,3,3 & \textrm{N/A} & 3,3,3,3 & \textrm{N/A} & 3,3,3,3 & \textrm{N/A} & \textrm{N/A} \\
\hline t_{q^{2},-q^{-2}} & 2,2,4,4 & \textrm{N/A} & 2,2,4,4 & 2,2,4,4 & \textrm{N/A} & \textrm{N/A} & \textrm{N/A} \\
\hline t_{q^{2},q^{2}} & 1,1,5,5 & 1,1,2,2,3,3 & \textrm{N/A} & \textrm{N/A} & \textrm{N/A} & \textrm{N/A} & \textrm{N/A} \\
\hline t_{q^{2},z} & 6,6 & 6,6 & 6,6 & 6,6 & 6,6 & 6,6 & \textrm{N/A} \\
\hline t_{z,q^{2}} & 6,6 & 6,6 & 6,6 & 6,6 & 6,6 & 6,6 & \textrm{N/A} \\
\hline t_{z,w} & 12 & 12 & 12 & 12 & 12 & 12 & 12 \\ \hline \end{array}$
\end{table}

\end{document}